\pdfoutput=1
\documentclass[linksfromyear,preprint]{imsart}

\RequirePackage{amsthm,amsmath,amssymb}
\RequirePackage[authoryear,round]{natbib}
\RequirePackage[colorlinks,citecolor=blue,urlcolor=blue]{hyperref}
\RequirePackage{graphicx}
\RequirePackage{bbm}
\usepackage[dvipsnames,svgnames, x11names]{xcolor}
\usepackage{booktabs,threeparttable,multirow,makecell} 
\usepackage{subfig,float}
\usepackage{tikz,pgfplots} 
\pgfplotsset{compat=1.11}
\usetikzlibrary{arrows,intersections}
\usetikzlibrary{backgrounds,positioning,fit,petri}
\usetikzlibrary{mindmap,trees,calendar,external,shapes}
\usetikzlibrary{decorations.fractals,decorations.pathmorphing,shadows,shadings,patterns}
\usetikzlibrary{spy,calc} 

\startlocaldefs

\numberwithin{equation}{section}
\theoremstyle{plain}
\newtheorem{thm}{Theorem}[section]

\newtheorem{lem}{Lemma}[section]
\theoremstyle{definition}
\newtheorem{ass}{Assumption}[section]
\newtheorem{defn}{Definition}[section]
\newtheorem{ex}{Example}[section]
\theoremstyle{remark}

\newtheoremstyle{exctd}
{\topsep} {\topsep}%
{\upshape}
{}
{\bfseries\scshape}
{.}
{1em}
{\thmname{#1} \thmnumber{ #2}\thmnote{#3} (cont.)}
\theoremstyle{exctd}
\newtheorem*{exctd}{Example}

\newcommand{\oset}[3][0ex]{%
	\mathrel{\mathop{#3}\limits^{
			\vbox to#1{\kern0\ex@
				\hbox{$\scriptstyle#2$}\vss}}}}

\newcommand{\uset}[3][0ex]{%
	\mathrel{\mathop{#3}\limits^{
			\vbox to#1{\kern6\ex@
				\hbox{$\scriptstyle#2$}\vss}}}}

\endlocaldefs

\newcommand{\comment}[1]{}

\newcommand{\convl}{\xrightarrow{L}}

\newcommand{\convp}{\xrightarrow{p}}

\begin{document}

\begin{frontmatter}

\title{Refinements of the Kiefer-Wolfowitz Theorem and a Test of Concavity\thanksref{t1}}
\thankstext{t1}{I am indebted to the associated editor and an anonymous referee for constructive comments that have helped greatly improve the presentation. This paper grew out of a coauthored work with Brendan K.\ Beare. I thank him for encouragement and advice. Portions of this research were conducted with the advanced computing resources and consultation provided by Texas A\&M High Performance Research Computing.}
\runtitle{Refinements of the Kiefer-Wolfowitz Theorem and a Test of Concavity}

\begin{aug}

\author{\fnms{Zheng} \snm{Fang}\ead[label=e2]{zfang@tamu.edu}}
\address{Department of Economics\\
	Texas A\&M University\\
	College Station, TX 77843\\
	\printead{e2}}

\runauthor{Z.\ Fang}

\end{aug}

\begin{abstract}
This paper studies estimation of and inference on a distribution function $F$ that is concave on the nonnegative half line and admits a density function $f$ with potentially unbounded support. When $F$ is strictly concave, we show that the supremum distance between the Grenander distribution estimator and the empirical distribution may still be of order $O(n^{-2/3}(\log n)^{2/3})$ almost surely, which reduces to an existing result of Kiefer and Wolfowitz when $f$ has bounded support. We further refine this result by allowing $F$ to be not strictly concave or even non-concave and instead requiring it be ``asymptotically'' strictly concave. Building on these results, we then develop a test of concavity of $F$ or equivalently monotonicity of $f$, which is shown to have asymptotically pointwise level control under the entire null as well as consistency under any fixed alternative. In fact, we show that our test has local size control and nontrivial local power against any local alternatives that do not approach the null too fast, which may be of interest given the irregularity of the problem. Extensions to settings involving testing concavity/convexity/monotonicity are discussed.
\end{abstract}

\begin{keyword}[class=MSC]
\kwd{62G05}
\kwd{62G07}
\kwd{62G09}
\kwd{62G10}
\kwd{62G20}
\end{keyword}

\begin{keyword}
\kwd{Grenander estimator}
\kwd{Kiefer-Wolfowitz theorem}
\kwd{Least concave majorant}
\kwd{Test of concavity/monotonicity}
\end{keyword}



\end{frontmatter}

\section{Introduction}

Let $\{X_i\}_{i=1}^n$ be a sample of i.i.d.\ nonnegative random variables with common distribution function $F$ that is known to be concave. The seminal work by \citet{Grenander1956II} establish that the nonparametric maximum likelihood estimator of $F$ is given by the least concave majorant $\hat{\mathbb F}_n$ of the empirical distribution function $\mathbb F_n$. In this context, the present paper provides new results on two aspects of the problem.

The first set of results concerns the closeness between $\hat{\mathbb F}_n$ and $\mathbb F_n$ relative to the uniform norm $\|\cdot\|_\infty$. The celebrated Marshall's lemma indicates that $\hat{\mathbb F}_n$ is $\sqrt{n}$-consistent for $F$ relative to $\|\cdot\|_\infty$---specifically, for each $n\in\mathbf N$,
\[
\|\hat{\mathbb F}_n-F\|_\infty\le\|\mathbb F_n-F\|_\infty~.
\]
\citet{Kiefer_Wolfowitz1976minimax} show that $\hat{\mathbb F}_n$ in fact is asymptotically equivalent to $\mathbb F_n$ by essentially establishing the following result:
\begin{thm}\label{Thm: KW}
Let $\alpha_{F}\equiv\inf\{x\in\mathbf R: F(x)=1\}$ and $F$ be twice continuously differentiable on $[0,\alpha_{F}]$ with $\sup\{x\in\mathbf R: F(x)=0\}=0$. If $\alpha_{F}<\infty$ and
\begin{align}
\beta_{F} & \equiv\inf_{0<x<\alpha_{F}}\frac{-f'(x)}{f^2(x)}>0~,\label{Eqn: KW beta}\\
\gamma_{F} & \equiv\frac{\sup_{0<x<\alpha_{F}}-f'(x)}{\inf_{0<x<\alpha_{F}}f^2(x)}<\infty~,\label{Eqn: KW gamma}
\end{align}
then it follows that
\[
\sup_{0\le x\le\alpha_{F}}|\hat{\mathbb F}_n(x)-\mathbb F_n(x)|=O(n^{-2/3}(\log n)^{2/3})\text{ almost surely}~.
\]
\end{thm}
\noindent Since $F(0)=0$ by assumption, it follows from $\alpha_{F}<\infty$ that the distribution associated with $F$ has bounded support. In turn, the condition \eqref{Eqn: KW beta} implies that $f'(x)<0$ for all $x\in(0,\alpha_{F})$ and hence that $F$ is strictly concave on the support $[0,\alpha_{F}]$. Finally, the condition \eqref{Eqn: KW gamma} demands that the derivative $f'$ of the density be small relative to the minimum of $f$ on its support. This paper provides refinements of Theorem \ref{Thm: KW} along several dimensions. First, we consider situations where $f$ has potentially unbounded support, i.e., $\alpha_{F}=\infty$, an extension which appears to be nontrivial---see Remark 5 in \citet[p.82]{Kiefer_Wolfowitz1976minimax} and also  \citet{Kiefer_Wolfowitz1977minimax}. In doing so, we necessarily have $\gamma_{F}=\infty$ (so that \eqref{Eqn: KW gamma} is being violated) because the density $f$ must approach zero along the right tail. Our second refinement, which may or may not be surprising, allows $F$ to be non-strictly concave or even non-concave but requires that it be ``asymptotically'' strictly concave as the sample size tends to infinity. In any case, the above asymptotic order result is uncovered under weaker regularity conditions.

The second set of results concerns inference aspect of the problem. In particular, we develop a test for the hypothesis that $F$ is concave or equivalently $f$ is nonincreasing. The insight we exploit here is that the hypothesis can be equivalently formulated as:
\[
\mathrm H_0: \|\hat F-F\|_p=0 \qquad \text{v.s.} \qquad \mathrm H_1: \|\hat F-F\|_p>0 ~,
\]
where $\hat F$ is the least concave majorant of $F$, and $\|\cdot\|_p$ is an $L^p$ norm (with appropriate weighting) for $p\in[1,\infty]$. Thus, it is natural for us to employ the test statistic $\sqrt n\|\hat{\mathbb F}_n-\mathbb F_n\|_p$. There are several technical challenges, however, in establishing statistical properties of our test. First, as demonstrated in \citet{BeareandMoon2015} and \citet{Beare_Fang2016Grenander}, the least concave majorant operator $F\mapsto \hat F$ is not fully differentiable if $F$ is concave but not strictly concave, and hence the conventional Delta method is inapplicable in establishing the asymptotic distribution of $\sqrt n\|\hat{\mathbb F}_n-\mathbb F_n\|_p$. However, the same authors have also shown that it is Hadamard directionally differentiable, a weaker notion of differentiability under which the Delta method is in fact preserved \citep{Shapiro1991,Dumbgen1993}. Second, as shown by \citet{Dumbgen1993} and \citet{FangSantos2018HDD}, even though the Delta method generalizes to deriving asymptotic distributions under Hadamard directional differentiability, it does not generalize to bootstrap consistency. This can be remedied by appealing to the rescaled bootstrap of \citet{Dumbgen1993}. Lastly, the weak limit of $\sqrt n\|\hat{\mathbb F}_n-\mathbb F_n\|_p$ is degenerate at zero when $F$ is strictly concave, a consequence of the Kiefer-Wolfowitz theorem. To build up a level $\alpha$ test, we leverage the asymptotic order results on $\|\hat{\mathbb F}_n-\mathbb F_n\|_\infty$ and propose a selection procedure we call the KW-selection that determines whether $F$ is strictly concave or not under the null. The idea is that the convergence rates of $\hat{\mathbb F}_n$ are different in these two cases, and hence we may introduce a suitable tuning parameter to identify the truth. Our test is in fact shown to have local size control which may be of interest since the weak limit of the statistic exhibits ``discontinuity'' with respect to the true distribution. Unfortunately, the local power will be poor under local (contiguous) alternatives that approach a strictly concave distribution.


The literature on testing concavity of the cumulative distribution function (or equivalently monotonicity of the density function) is surprisingly limited. In Section \ref{Sec: comparison}, we review related studies in this regard and compare them to our test. Overall, there are roughly two types of tests in the density context. The first type is concerned with special classes of alternatives. This includes the work of \citet{Woodroofe_Sun1999} who test uniformity against monotonically increasing (but not uniform) densities, or uniformity against convex cdfs. Therefore, these tests are powerful in detecting these special classes by design, but may perform poorly for other alternatives, which is confirmed in our simulation studies. The second type remains agnostic about the particular natures of alternative hypothesis, but relies on critical values from the asymptotically least favorable distributions. This includes the tests of \citet{Durot2003testing} and \citet{Kulikov_Lopuhaa2004testing}. While the use of least favorable distributions is common practice in some other settings \citep{Carolan_Tebbs2005,Delgado_Escanciano2012monotonicity,DelgadoEscanciano2013ConditionalSD,DelgadoEscanciano2016ConditionalMI}, the resulting tests may be too conservative and cause substantial power loss. Finally, in similar contexts, there are recent studies such as \citet{BeareShi2019DRO} and \citet{Seo2018SM} that aim to improve power by bootstrap. However, due to the lack of Kiefer-Wolfowitz type results, it is unclear what happens to their tests when the limiting distributions of the test statistics are degenerate. We shall remove the dependence on least favorable distributions by bootstrap, and control the size of our test (both pointwise and locally) even when degeneracy occurs by utilizing the Kiefer-Wolfowitz's theorem.

We now introduce some notation. We denote by $\mathbf R_+$ the set of nonnegative real numbers. For a sequence $\{Z_n\}_{n=1}^\infty$ of random variables, we write $Z_n=O_{a.s.}(a_n)$ for a deterministic sequence $\{a_n\}_{n=1}^\infty$ if $Z_n/a_n=O(1)$ as $n\to\infty$ almost surely. Here and elsewhere in this paper, for a sequence $\{b_n\}$ of real numbers, $b_n=O(1)$ as $n\to\infty$ means that $\{b_n\}$ is bounded. For two sequences of real numbers $\{a_n\}$ and $\{b_n\}$, we write $a_n\lesssim b_n$ if $a_n/b_n =O(1)$ as $n\to\infty$, and write $a_n\gtrsim b_n$ if $b_n\lesssim a_n$. Moreover, the notation $a_n\asymp b_n$ means that $a_n/b_n=O(1)$ and $b_n/a_n=O(1)$ as $n\to\infty$. Analogously, for two functions $f,g:\mathbf R\to\mathbf R_+$, we write $f(\epsilon)\asymp g(\epsilon)$ as $\epsilon\downarrow 0$ if $f(\epsilon)/g(\epsilon)=O(1)$ and $g(\epsilon)/f(\epsilon)=O(1)$ as $\epsilon\downarrow 0$. For an arbitrary nonempty set $T$, $\ell^\infty(T)$ is the space of bounded real-valued functions on $\mathbf R_+$ equipped with the uniform norm $\|\cdot\|_\infty$, i.e., $\|f\|_\infty\equiv \sup_{x\in T}|f(x)|$. The space $C_0(\mathbf R_+)$ is the family of all real-valued continuous functions on $\mathbf R_+$ that vanish at infinity. We denote by $\convl$ the weak convergence in the sense of Hoffmann-J{\o}rgensen \citep{Vaart1996}.

The remainder of the paper is structured as follows. Section \ref{Sec: KW} presents refinements of Theorem \ref{Thm: KW}. Section \ref{Sec: tests} formally develops our test and then compare to some existing monotonicity tests from the literature. Section \ref{Sec: simulation} conducts simulation studies, with some of the results relegated to Appendix \ref{App: Simulation}. Section \ref{Sec: conclusion} concludes. All proofs are collected in Appendix \ref{Sec: proofs}, while extensions of the test results to a general setup are discussed in Appendix \ref{App: extension}.

\section{The Kiefer-Wolfowitz Theorems}\label{Sec: KW}

This section presents refinements of the Kiefer-Wolfowitz theorem that allow the support of the density function $f$ to be unbounded. Throughout, we think of the cdf $F$ as a function on $\mathbf R_+$ (rather than on $\mathbf R$). We proceed with the following assumption.

\begin{ass}\label{Ass: setup}
(i) $\{X_i\}_{i=1}^n$ is an i.i.d.\ sample with common distribution function $F$, and (ii) $F$ is strictly concave on $\mathbf R_+$.
\end{ass}

Assumption \ref{Ass: setup} simply formalizes the i.i.d.\ setup and strict concavity of the distribution function $F$. To establish the refinements, we need to impose further regularity conditions on $F$. In particular, we need to introduce an analog of $\beta_F$ in \eqref{Eqn: KW beta}. Specifically, for any $\epsilon\in[0,1]$, we define
\begin{align}\label{Eqn: beta}
\bar\beta(\epsilon)\equiv \inf_{0\le x\le F^{-1}(1-\epsilon)}\frac{-f'(x)}{f^2(x)}~,
\end{align}
where we suppress the dependence of $\bar\beta(\epsilon)$ on the underling distribution for simplicity. Thus, $\bar\beta(\epsilon)$ may be viewed as a truncated version of $\beta_F$ in \eqref{Eqn: KW beta}. As it turns out, while $\bar\beta(\epsilon)$ is allowed to approach zero as $\epsilon\to 0$, the speed at which it approaches zero determines the asymptotic order of $\|\hat{\mathbb F}_n-\mathbb F_n\|_\infty$. To formalize our discussions, we now impose:

\begin{ass}\label{Ass: equivalence}
(i) $F$ is twice continuously differentiable on $\mathbf R_+$, (ii) $\bar\beta(\epsilon) \asymp \epsilon^\tau$ as $\epsilon\to 0$ for some constant $\tau>-1$, and (iii)
\[
\|F''\|_{1/2}\equiv\left[\int_0^\infty \sqrt{|F''(x)|}\,\mathrm dx\right]^2<\infty~.
\]
\end{ass}

Assumption \ref{Ass: equivalence}(i) is the same smoothness condition required by Theorem \ref{Thm: KW}. Assumption \ref{Ass: equivalence}(ii) characterizes the exact rates at which $\bar\beta(\epsilon)$ is allowed to approach zero as $\epsilon\to 0$. The special case $\tau=0$ implies $\inf_{0\le x<\infty}\{-f'(x)/f^2(x)\}>0$, which is exactly the condition \eqref{Eqn: KW beta} in the case of unbounded support. If $\tau>0$, then $\bar\beta(\epsilon)\to 0$ as $\epsilon\to 0$; if $\tau\in(-1,0)$, then $\bar\beta(\epsilon)\to \infty$ as $\epsilon\to 0$ which still fulfills \eqref{Eqn: KW beta}. Assumption \ref{Ass: equivalence}(iii) is a technical condition that serves to control the interpolation error for the cdf $F$. When the support of the density function $f$ is bounded, Assumption \ref{Ass: equivalence}(iii) is automatically fulfilled because $F''$ as a continuous function is bounded. Thus, our assumptions are strictly weaker than those imposed in \citet{Kiefer_Wolfowitz1976minimax}.

Given Assumptions \ref{Ass: setup} and \ref{Ass: equivalence}, we now present the first refinement.

\begin{thm}\label{Thm: Equivalence}
If Assumptions \ref{Ass: setup} and \ref{Ass: equivalence} hold, then, as $n\to\infty$,
\[
\|\hat{\mathbb F}_n-\mathbb F_n\|_\infty=\begin{cases}
O_{a.s.}\big(((\log n)/n)^{\frac{2}{2\tau+3}}\big) & \text{ if }\tau\ge 0\\
O_{a.s.}\big(((\log n)/n)^{\frac{\tau+2}{2\tau+3}}\big) & \text{ if }\tau\in (-1,0)\\
\end{cases}~.
\]
\end{thm}

Theorem \ref{Thm: Equivalence} delivers the asymptotic order of $\|\hat{\mathbb F}_n-\mathbb F_n\|_\infty$, which crucially depends on the parameter $\tau$. The slower $\bar\beta(\epsilon)$ approaches zero (as $\epsilon\to 0$) or the smaller $\tau$ is, the closer $\hat{\mathbb F}_n$ and $\mathbb F_n$ is (asymptotically). If $\bar\beta(\epsilon)$ approaches zero too fast or $\tau$ is too large, then the asymptotic equivalence of $\hat{\mathbb F}_n$ and $\mathbb F_n$ are no longer implied. This happens precisely when $\tau\ge 1/2$, in which case
\begin{align}
\sqrt n ((\log n)/n)^{\frac{2}{2\tau+3}}= n^{\frac{2\tau-1}{2\tau+3}} (\log n)^{\frac{2}{2\tau+3}}\to\infty\text{ as } n\to\infty~.
\end{align}
The special case $\tau=0$ leads to $\|\hat{\mathbb F}_n-\mathbb F_n\|_\infty=O_{a.s.}\big(((\log n)/n)^{2/3}\big)$, which is exactly the result in Theorem \ref{Thm: KW}. There are two major new ingredients in establishing the theorem, compared to the proof of Theorem \ref{Thm: KW}. One is controlling the variation of $\mathbb F_n$ over a small upper quantile region---see the treatment of the $T_{n,k_n+1}$ term and the $B_n$ term in the proof. The other is bounding the interpolation error for $F$ by the $L_{1/2}$-integral of $F''$ based on a result from approximation theory \citep{Burchard1974splines}. The latter allows us to dispense with the regularity condition $\gamma_F<\infty$ \citep[p.82]{Kiefer_Wolfowitz1976minimax}.

Our second refinement allows the common distribution function $F$ that the finite sample $X_1,\ldots,X_n$ share to be not strictly concave or even non-concave. However, we do require that it vary with the sample size $n$ in such a way that it approaches (in a suitable sense to be specified) a strictly concave distribution function as $n\to\infty$, a notion we call ``local to concavity.'' In order to formalize our local analysis, we denote by $\mathbf P$ the set of probability measures with the support contained in $\mathbf R_+$ that possibly govern the data:
\[
\mathbf P\equiv\{P\in\mathbf M: P(\mathbf R_+)=1\}~,
\]
where $\mathbf M$ is the set of all Borel probability measures on $\mathbf R$ that are dominated by the Lebesgue measure on $\mathbf R$. Further, we think of the distribution function $F\equiv F(P)$ as a map $F:\mathbf P\to\ell^\infty(\mathbf R_+)$ defined by
\[
F(P)(x)\equiv P((-\infty,x])~,\,\forall\,x\in\mathbf R_+~.
\]
We may now formalize the precise meaning of ``local'' as follows.

\begin{defn}\label{Defn: differentiable path}
A function $t\mapsto P_t$ mapping a neighborhood $(-\epsilon,\epsilon)$ of zero into $\mathbf P$ is called a differentiable path passing through $P$ if, for $P_0=P$ and some function $h: \mathbf R_+\to\mathbf R$,
\begin{align}\label{Eqn: differentiable path}
\lim_{t\to 0}\int\left[\frac{\mathrm dP_t^{1/2}-\mathrm dP^{1/2}}{t}-\frac{1}{2}h\,\mathrm dP^{1/2}\right]^2= 0~.
\end{align}
\end{defn}

The notation $\mathrm dP_t$ and $\mathrm dP$ may be understood as the densities of $P_t$ and $P$ with respect to some dominating measure $\mu_t$ (for each $t$), though the integral in \eqref{Eqn: differentiable path} does not depend on the choice of $\mu_t$ \citep[p.362]{Vaart1998}. Loosely speaking, \eqref{Eqn: differentiable path} implies that $P_t$ gets closer and closer to $P_0$ ``on average'', as $t\to 0$. The function $h$ is referred to as the score function of $P$ and satisfies $\int h\,\mathrm dP=0$ and $h\in L^2(P)$---see, for example, Lemma 25.14 in \citet{Vaart1998}. The term ``score function'' makes sense in view of the relation $h(x)=\frac{\mathrm d}{\mathrm dt}\log \mathrm dP_t(x)\big|_{t=0}$ (which is the usual definition of score function) under regularity conditions---see, for example, Lemma 7.6 in \citet{Vaart1998}.

If $\{P_t: |t|<\epsilon\}$ is a differentiable path, then we obtain by Example 5.3.1 in \citet{BKRW993Efficient} that
\begin{align}\label{Eqn: F diff}
\lim_{t\to 0}\|\frac{F(P_t)-F(P)}{t}-\dot F(h)\|_\infty=0~.
\end{align}
where
\[
\dot F(h)(x)\equiv\int_{-\infty}^x h(u)\,P(\mathrm du)~.
\]
Intuitively, \eqref{Eqn: F diff} means that the distribution function $F_t\equiv F(P_t)$ smoothly passes through $F\equiv F(P_0)$ at the same speed as the underlying probability measure $P_t$ passes through $P_0$. For a differentiable path $\{P_t\}$ satisfying \eqref{Eqn: differentiable path} and $F_t\equiv F(P_t)$, we say $\{F_t\}$ is {\it local to concavity} if $F$ is concave and is {\it local to strict concavity} if $F$ is strictly concave.

We next show that the conclusion in Theorem \ref{Thm: KW} is preserved under local to strict concavity. To this end, we impose the following assumption.

\begin{ass}\label{Ass: local}
(i) $\{X_i\}_{i=1}^n$ is an i.i.d.\ sample with common probability measure $P_{1/\sqrt{n}}$, and (ii) $\{P_{1/\sqrt{n}}\}\subset\mathbf P$ corresponds to the differentiable path $\{P_t\}$ passing through $P_0\equiv P$ in the sense of Definition \ref{Defn: differentiable path}.
\end{ass}

\begin{thm}\label{Thm: Equivalence, local}
If Assumptions \ref{Ass: setup}(ii), \ref{Ass: equivalence} and \ref{Ass: local} hold, then, as $n\to\infty$,
\[
\|\hat{\mathbb F}_n-\mathbb F_n\|_\infty=\begin{cases}
O_{a.s.}\big(((\log n)/n)^{\frac{2}{2\tau+3}}\big) & \text{ if }\tau\ge 0\\
O_{a.s.}\big(((\log n)/n)^{\frac{\tau+2}{2\tau+3}}\big) & \text{ if }\tau\in (-1,0)\\
\end{cases}~.
\]
\end{thm}

Theorem \ref{Thm: Equivalence, local} is not much deeper than Theorem \ref{Thm: KW}, but it has implications for the problem of testing concavity as we shall elaborate in Section \ref{Sec: tests}.

\section{Testing Concavity}\label{Sec: tests}

In this section, we develop a test of concavity that controls size regardless of whether concavity is strict or not. This is accomplished by building on the asymptotic order results established previously. We shall also compare our concavity test with existing monotonicity tests. For simplicity, in what follows, we focus on the canonical case $\tau=0$ in Assumption \ref{Ass: equivalence}(ii). The general case makes no essential differences---see Appendix \ref{App: extension} for extensions.

\subsection{The Test Statistic}

First, we introduce our test statistic and derive its asymptotic distribution. To this end, we introduce the least concave majorant (LCM) operator following \citet{Beare_Fang2016Grenander}---see also \citet{BeareandMoon2015}.

\begin{defn}
Given a convex set $T\subseteq \mathbf R_+$, the LCM over $T$ is the operator $\mathcal M_T: \ell^\infty(\mathbf R_+)\to\ell^\infty(T)$ that maps each $\theta\in\ell^\infty(\mathbf R_+)$ to the function
\[
\mathcal M_T\theta(x)\equiv\inf\{g(x): g\in \ell^\infty(T)\text{, }g\text{ is concave, and }\theta\le g\text{ on }T\},\quad \forall \, x\in T ~.
\]
If $T=\mathbf R_+$, then we write $\mathcal M\equiv\mathcal M_{\mathbf R_+}$.
\end{defn}

The hypothesis of our interest can now be formulated as:
\[
\mathrm H_0: \phi(F)=0\quad \text{v.s.}\quad \mathrm H_1: \phi(F)>0~,
\]
where $\phi: \ell^\infty(\mathbf R_+)\to\mathbf R$ is defined by
\[
\phi(F)=\|\mathcal MF-F\|_p=\begin{cases}
[\int_{\mathbf R_+}(\mathcal MF(x)-F(x))^p g(x)\,\mathrm dx]^{1/p} & \text{ if }p\in[1,\infty)\\
\sup_{x\in\mathbf R_+}|\mathcal MF(x)-F(x)| & \text{ if }p=\infty
\end{cases}~,
\]
with $g\in L^1(\mathbf R_+)$ some known positive weighting function. Here and throughout, we work with an arbitrarily fixed $p\in[1,\infty]$, and thus suppress the dependence of $\phi$ on $p$ for notational simplicity. Given the above formulation, we employ the statistic $\sqrt n \phi(\mathbb F_n)$. To derive the weak limit of $\sqrt n \phi(\mathbb F_n)$, note that under the null hypothesis, $\mathcal MF=F$ and hence we may rewrite:
\begin{align}\label{Eqn: test statistic using D}
\sqrt n \phi(\mathbb F_n)=\sqrt n\|\mathcal M\mathbb F_n-\mathbb F_n\|_p=\|\sqrt n\{\mathcal D\mathbb F_n-\mathcal D F\}\|_p~,
\end{align}
where $\mathcal D \equiv\mathcal M-\mathcal I$ with $\mathcal I$ the identity operator on $\ell^\infty(\mathbf R_+)$. Thus, the asymptotic distribution of $\sqrt n \phi(\mathbb F_n)$ would be an implication of the continuous mapping theorem and the Delta method, if we could show that $\mathcal M$ (and hence $\mathcal D$) is Hadamard differentiable \citep[p.372-374]{Vaart1996}. Unfortunately, as demonstrated by \citet{BeareandMoon2015} and \citet{Beare_Fang2016Grenander}, the LCM operator fails to be fully differentiable but only Hadamard directionally differentiable in general \citep{Shapiro1990}. Nonetheless, this type of directional differentiability suffices for applying a generalized version of the Delta method \citep{Shapiro1991,Dumbgen1993}.

\begin{defn}\label{Defn: HDD}
Let $\mathbb D$ and $\mathbb E$ be normed spaces equipped with norms $\|\cdot\|_{\mathbb D}$ and $\|\cdot\|_{\mathbb E}$ respectively, and $\phi:\mathbb D_\phi\subseteq \mathbb D\to\mathbb E$. The map $\phi$ is said to be {\it Hadamard directionally differentiable} at $\theta \in\mathbb D_\phi$ {\it tangentially} to $\mathbb D_0\subseteq\mathbb D$, if there is a map $\phi_\theta':\mathbb D_0\to\mathbb E$ such that:
\begin{equation}\label{Eqn: HDD}
\lim_{n\rightarrow \infty}\|\frac{\phi(\theta +t_n h_n)-\phi(\theta)}{t_n} -\phi_\theta'(h)\|_{\mathbb E} = 0 ~,
\end{equation}
for all sequences $\{h_n\}\subset\mathbb D$ and $\{t_n\}\subset\mathbf R_+$ such that $t_n\downarrow 0$, $h_n\to h\in\mathbb D_0$ as $n\to\infty$ and $\theta+t_nh_n\in\mathbb D_\phi$ for all $n$.
\end{defn}

The defining feature of Hadamard directional differentiability is that, unlike Hadamard (full) differentiability, the directional derivative is in general nonlinear though necessarily continuous and positively homogeneous of degree one \citep{Shapiro1990}.\footnote{Let $\mathbb D$ be a vector space. A function $f: \mathbb D\to\mathbf R$ is said to be positively homogeneous of degree one if and only if $f(ax) = af(x)$ for all $x\in\mathbb D$ and all $a\in\mathbf R_+$.} We refer the reader to \citet{Shapiro1990} and a more recent review by \citet{FangSantos2018HDD} for additional discussions. Proposition 2.1 in \citet{Beare_Fang2016Grenander} implies that $\mathcal M$ is Hadamard directionally differentiable at any concave $F\in\ell^\infty(\mathbf R_+)$ tangentially to the set $C_0(\mathbf R_+)$ with the derivative $\mathcal M_F': C_0(\mathbf R_+)\to\ell^\infty(\mathbf R_+)$ given by: for any $h\in C_0(\mathbf R_+)$ and $x\in \mathbf R_+$,
\[
(\mathcal M_F'h)(x)=(\mathcal M_{T_{F,x}}h)(x)~,
\]
where $T_{F,x}=\{x\}\cup U_{F,x}$ with $U_{F,x}$ the union of all open intervals $A\subseteq \mathbf R_+$ such that (i) $x\in A$, and (ii) $F$ is affine over $A$. We emphasize that $\mathcal M_F'$ is equal to the identity operator on $C_0(\mathbf R_+)$ if and only if $F$ is strictly concave, in which case $\mathcal M$ is Hadamard differentiable at $F$ tangentially to $C_0(\mathbf R_+)$---see Proposition 2.2 in \citet{Beare_Fang2016Grenander} for more details.

We are now in a position to state the asymptotic distribution of our statistic $\sqrt n \phi(\mathbb F_n)$ by invoking the generalized Delta method.

\begin{lem}\label{Lem: weak limit}
If Assumption \ref{Ass: setup}(i) holds, then $\sqrt n \phi(\mathbb F_n)\convl \|\mathcal M_F'\mathbb G-\mathbb G\|_p$ under $\mathrm H_0$, where $\mathbb G\equiv\mathbb B\circ F$ and $\mathbb B$ is the standard Brownian bridge on $[0,1]$.
\end{lem}

Lemma \ref{Lem: weak limit} establishes the weak limit of the test statistic $\sqrt n \phi(\mathbb F_n)$ under the null. The Brownian bridge $t\mapsto\mathbb G(t)\equiv\mathbb B(F(t))$ is a zero mean Gaussian process with covariance function: for all $s,t\in[0,1]$,
\begin{align}
\mathrm{Cov}(\mathbb G(t),\mathbb G(s))=F(\min(s,t))-F(t)F(s)~.
\end{align}
The limiting distribution in Lemma \ref{Lem: weak limit} is not pivotal because it depends on $F$ through the process $\mathbb G$ and, critically, the derivative $\mathcal M_F'$.

\subsection{The Critical Values}

Towards constructing critical values for our test, we next aim to estimate the law of $\|\mathcal M_F'\mathbb G-\mathbb G\|_p$ through bootstrap. There are, however, two complications involved, as we now elaborate.

First, when $F$ is non-strictly concave in which case $\mathcal M$ is only Hadamard directionally differentiable, the standard bootstrap (compare to \eqref{Eqn: test statistic using D}),
\begin{align}\label{Eqn: standard bootstrap}
\|\sqrt n\{\mathcal D\mathbb F_n^*-\mathcal D\mathbb F_n\}\|_p~,
\end{align}
is necessarily inconsistent. This is a consequence of Proposition 1 in \citet{Dumbgen1993}, which has been formalized by Theorem A.1 in \citet{FangSantos2018HDD}. Second, when $F$ is strictly concave, the weak limit of $\sqrt n \phi(\mathbb F_n)$ is degenerate at zero since $\mathcal M_F'=\mathcal I$, which is not surprising in view of Theorem \ref{Thm: Equivalence}. The first issue can be resolved by appealing to the rescaled bootstrap in \citet{Dumbgen1993}. As shall be seen shortly, the rescaled bootstrap is connected in a subtle way to the modified bootstrap in \citet{FangSantos2018HDD}; namely, it amounts to composing a suitable derivative estimator (see \eqref{Eqn: rescaled bootsrap2}) with some bootstrap process. The second issue is more challenging, which we shall fix by leveraging Theorems \ref{Thm: Equivalence} and \ref{Thm: Equivalence, local} as we describe now.

The key to the level control of our test is a selection procedure that determines whether the concavity of $F$ is strict or not by exploiting the fact that the convergence rates in these two cases are different (see Theorem \ref{Thm: Equivalence} and Lemma \ref{Lem: weak limit}). Specifically, let $\{\kappa_n\}$ be a sequence of positive scalars such that $\kappa_n=o(1)$ and $(\log n)^{2/3}n^{-1/6}/\kappa_n=o(1)$ as $n\to\infty$. For example, we may take $\kappa_n=n^{-1/7}$ or $(\log n)^{-1}$. Define for each $n\in\mathbf N$,
\begin{align}\label{Eqn: KW selection}
\xi_n\equiv \frac{\sqrt{n}\|\hat{\mathbb F}_n-\mathbb F_n\|_p}{\kappa_n}~.
\end{align}
Then $\xi_n\convp 0$ if $F$ is strictly concave by Theorem \ref{Thm: Equivalence} and $\xi_n\convp\infty$ if $F$ is non-strictly concave by Lemma \ref{Lem: weak limit}. Thus, if $\sqrt n\phi(\mathbb F_n)\le \kappa_n$, we may take $F$ to be strictly concave and non-strictly concave otherwise. We refer to such a selection procedure as the KW-selection. As a result, the asymptotic level of our test can be controlled in this case by choosing the critical value to be $\kappa_n$.

On the other hand, when $F$ is non-strictly concave, the nondegenerate limit $\|\mathcal M_F'(\mathbb G)-\mathbb G\|_p$ in Lemma \ref{Lem: weak limit} can be estimated by composing a suitable estimator $\hat{\mathcal M}_n'$ of $\mathcal M_F'$ with the nonparametric bootstrap estimator $\sqrt{n}\{\mathbb F_n^*-\mathbb F_n\}$ of the Gaussian process $\mathbb G$, i.e., we employ the bootstrap estimator
\begin{align}\label{Eqn: rescaled bootsrap}
\|\hat{\mathcal M}_n'(\sqrt{n}\{\mathbb F_n^*-\mathbb F_n\})-\sqrt{n}\{\mathbb F_n^*-\mathbb F_n\}\|_p~,
\end{align}
where $\mathbb F_n^*(x)=\frac{1}{n}\sum_{i=1}^{n}W_{ni}1\{X_i\le x\}$ for $x\in\mathbf R_+$, $W_n\equiv(W_{n1},\ldots,W_{nn})$ a multinomial random vector with $n$ categories and probabilities $(1/n,\ldots,1/n)$, and $\hat{\mathcal M}_n': \ell^\infty(\mathbf R_+)\to\ell^\infty(\mathbf R_+)$ is some appropriate derivative estimator. Such a bootstrap procedure is justified by Theorem 3.2 in \citet{FangSantos2018HDD}. It is in fact precisely the rescaled bootstrap proposed by \citet{Dumbgen1993}, which amounts to estimating $\mathcal M_F'$ by: for any $h\in\ell^\infty(\mathbf R_+)$,
\begin{align}\label{Eqn: rescaled bootsrap2}
\hat{\mathcal M}_n'(h)=\frac{\mathcal M(\mathbb F_n+t_nh)-\mathcal M(\mathbb F_n)}{t_n}~,
\end{align}
where $t_n\downarrow 0$ such that $t_n\sqrt n\to\infty$---see \citet[p.390-1]{FangSantos2018HDD}. In turn, we may choose the critical value in this case as: for $\alpha\in(0,1)$,
\[
\hat c_{1-\alpha}^*\equiv \inf\{c\in\mathbf R: P_W(\|\hat{\mathcal M}_n'(\sqrt{n}\{\mathbb F_n^*-\mathbb F_n\})-\sqrt{n}\{\mathbb F_n^*-\mathbb F_n\}\|_p\le c)\ge 1-\alpha\}~,
\]
where $P_W$ denotes the probability with respect to the bootstrap weights $\{W_{ni}\}_{i=1}^n$ holding the data $\{X_i\}_{i=1}^n$ fixed. In practice, $\hat c_{1-\alpha}^*$ can be estimated by Monte Carlo simulations by drawing a large number of bootstrap samples $\{X_{b,i}^*\}_{i=1}^n$ from $\mathbb F_n$ for $b=1,\ldots,B$ (with the data $\{X_i\}_{i=1}^n$ fixed).

Now, for each fixed $\alpha\in(0,1)$ we set the critical value of our test as
\[
\hat c_{1-\alpha}\equiv\max\{\kappa_n,\hat c_{1-\alpha}^*\}~,
\]
which provides pointwise asymptotic level control as confirmed by the following theorem. Before stating the theorem, we formalize requirements on $t_n$, $\kappa_n$, and the $(1-\alpha)$th quantile $c_{1-\alpha}^*$ of the cdf of the weak limit $\|\mathcal M_F'(\mathbb G)-\mathbb G\|_p$.

\begin{ass}\label{Ass: kn}
(i) $\{t_n\}$ is a sequence of positive scalars such that $t_n\downarrow 0$ such that $t_n\sqrt n\to\infty$ as $n\to\infty$; (ii) $\{\kappa_n\}$ is a sequence of positive scalars such that $\kappa_n\downarrow 0$ and $(\log n)^{2/3}n^{-1/6}/\kappa_n\to 0$ as $n\to\infty$.
\end{ass}

\begin{ass}\label{Ass: quantile regularity}
The cdf of $\|\mathcal M_F'(\mathbb G)-\mathbb G\|_p$ is continuous and strictly increasing at its $(1-\alpha)$th quantile $c_{1-\alpha}^*$ when $F$ is non-strictly concave.
\end{ass}

Assumption \ref{Ass: kn} imposes the convergence rates on the tuning parameters $\kappa_n$ and $t_n$. We do not touch the challenging issue of optimal choices in this paper. Assumption \ref{Ass: quantile regularity} is a standard technical condition (often implicitly imposed) to ensure that consistent bootstrap can produce consistent critical values---see Lemma 11.2.1 in \citet{TSH2005}. In verifying Assumption \ref{Ass: quantile regularity}, we note that, by the proof of Theorem \ref{Thm: test concavity, local}, the map $g\mapsto \phi_F'(g)\equiv\|\mathcal M_F' g-g\|_p$ is subadditive\footnote{Let $\mathbb D$ be a vector space. A function $f: \mathbb D\to\mathbf R$ is said to be subadditive if and only if $f(x+y)\le f(x)+f(y)$ for all $x,y\in\mathbb D$.} and hence convex since it is obviously positively homogeneous of degree one. Hence, by Theorem 11.1 in \citet{Davydov1998local}, the distribution function $H$ of $\|\mathcal M_F'(\mathbb G)-\mathbb G\|_p$ is absolutely continuous and strictly increasing on $(r_0,\infty)$, with $r_0\equiv\inf\{r\in\mathbf R: H(r)>0\}$. Therefore, Assumption \ref{Ass: quantile regularity} holds whenever $c_{1-\alpha}^*>r_0$.

We are now in a position to state the first main result of this paper.

\begin{thm}\label{Thm: test concavity, pointwise}
If Assumptions \ref{Ass: setup}(i), \ref{Ass: equivalence} with $\tau=0$, \ref{Ass: kn} and \ref{Ass: quantile regularity} hold, then it follows that, under $\mathrm H_0$,
\begin{align}
\limsup_{n\to\infty}P(\sqrt{n}\phi(\mathbb F_n)> \hat c_{1-\alpha})\le\alpha~,
\end{align}
and under $\mathrm H_1$,
\begin{align}
\lim_{n\to\infty}P(\sqrt{n}\phi(\mathbb F_n)> \hat c_{1-\alpha})=1~.
\end{align}
\end{thm}

Theorem \ref{Thm: test concavity, pointwise} shows that our test is pointwise (in $P$) asymptotically level $\alpha$ and consistent under any fixed alternative. However, in view of the irregularity of the problem and as argued in the literature \citep{TSH2005,AndrewsandGuggen2010ET,Romano_Shaikh2012uniform}, pointwise asymptotics can be unreliable in ``nonstandard'' settings. It is therefore of interest to investigate the uniform or at least local properties of our test.

To study both local size control and local power of our proposed test, we follow \citet[p.384-6]{Vaart1998} and consider differentiable paths in $\mathbf P$ passing through $P_0\equiv P$ that also belong to the set
\[
\mathbf H\equiv \{\{P_t\}: \text{(i)} \phi(F(P_t))=0 \text{ if }t\le 0\text{, and (ii)} \phi(F(P_t))>0 \text{ if }t> 0\}~.
\]
Thus, a differentiable path $\{P_t\}$ in $\mathbf H$ is such that if it satisfies the null hypothesis whenever $t\le 0$ but otherwise the alternative for all $t>0$. For a differentiable path $\{P_{\eta/\sqrt{n}}\}$ where $\eta\in\mathbf R$, we set $P^n\equiv \prod_{i=1}^n P$, $P_n^n\equiv\prod_{i=1}^n P_{\eta/\sqrt{n}}$ and define the power function of our test for sample size $n$ as
\[
\pi_n(P_{\eta/\sqrt{n}})\equiv P_n^n(\sqrt{n}\phi(\mathbb F_n)>\hat c_{1-\alpha})~.
\]

Our next theorem establishes local properties of our test.

\begin{thm}\label{Thm: test concavity, local}
 Let $\{P_{t}\}$ be a differentiable path in $\mathbf H$ and let Assumptions \ref{Ass: equivalence}, \ref{Ass: local}, \ref{Ass: kn}, and \ref{Ass: quantile regularity} hold with $\tau=0$ and $F\equiv F(P_0)$. Then it follows that
\begin{enumerate}
\item For any $\eta\in\mathbf R$, (i) if $F$ is strictly concave, then $\liminf_{n\to\infty}\pi_n(P_{\eta/\sqrt{n}})=0$, and (ii) if $F$ is non-strictly concave, then
\begin{multline}\label{Thm: test concavity, limit local}
\liminf_{n\to\infty}\pi_n(P_{\eta/\sqrt{n}})\\ \ge P(\|\mathcal M_F' (\mathbb G+\eta\dot F(h))-\{\mathbb G+\eta\dot F(h)\}\|_p>c_{1-\alpha}^*)~,
\end{multline}
where \eqref{Thm: test concavity, limit local} holds with equality if, in addition, the cdf of $\|\mathcal M_F' (\mathbb G+\eta\dot F(h))-\{\mathbb G+\eta\dot F(h)\}\|_p$ is continuous at $c_{1-\alpha}^*$.
\item For any $\eta\le 0$, we have
\begin{align}\label{Thm: test concavity, local size}
\limsup_{n\to\infty}\pi_n(P_{\eta/\sqrt{n}})\le\alpha~.
\end{align}
\item If $\|\mathcal M_F' (\dot F(h))-\dot F(h)\|_p>0$, then
\begin{align}\label{Thm: test concavity, local power}
\liminf_{\eta\uparrow \infty} \liminf_{n\to\infty}\pi_n(P_{\eta/\sqrt{n}}) = 1~.
\end{align}
\end{enumerate}
\end{thm}

The first part of Theorem \ref{Thm: test concavity, local} delivers a lower bound of the local limiting power function. If the true distribution function $F_n\equiv F(P_{\eta/\sqrt n})$ is local to a strictly concave function $F$, then the limiting local power along the path $\{P_{1/\sqrt n}\}$ is zero. This is unfortunate and in fact is an implication of Theorem \ref{Thm: Equivalence, local}. The second part shows that the asymptotic null rejection rate along $\{P_{1/\sqrt n}\}$ is no larger than $\alpha$, establishing the asymptotic local size control.

The first part of Theorem \ref{Thm: test concavity, local} might leave one the impression that the test has poor local power against any sequence of local alternatives $\{P_{\eta/\sqrt n}\}$ (with $\eta>0$). The third part is intended to reconcile such a misconception. Heuristically, it says that our test has nontrivial local power if the sequence $\{P_{\eta/\sqrt n}\}$ does not approach the null too fast. To appreciate the condition $\|\mathcal M_F' (\dot F(h))-\dot F(h)\|_p>0$, we first note that it prevents $F$ from being strictly concave in which case $\mathcal M_F'$ is the identity operator and the limiting local power is zero by Part 1-(i). Second, by \eqref{Eqn: F diff}, we have that, as $n\to\infty$,
\begin{align}
\sqrt n \{F_n- F\}=\sqrt n \{F(P_{h/\sqrt n})- F(P)\}\to \eta\dot F(h)~.
\end{align}
In turn, this implies by Proposition 2.1 in \citet{Beare_Fang2016Grenander} and $\phi(F)=0$ (by the definition of $\mathbf H$) that, as $n\to\infty$,
\begin{align}
\sqrt n \phi(F_n) =\sqrt n \{\phi(F_n)- \phi(F)\} \to \|\mathcal M_F'(\eta\dot F(h))-\eta\dot F(h)\|_p~.
\end{align}
By the positive homogeneity of degree one of $\mathcal M_F'$ (as a Hadamard directional derivative), we thus have: for $\Delta\equiv \|\mathcal M_F'(\dot F(h))-\dot F(h)\|_p$,
\begin{align}
\phi(F_n)= \frac{\eta}{\sqrt n} (\Delta+o(1))~,
\end{align}
for all $\eta>0$. Therefore, $\Delta\equiv \|\mathcal M_F'(\dot F(h))-\dot F(h)\|_p>0$ implies that the third part is concerned with local power along a nontrivial Pitman drift, a canonical device for local power analysis. When $\Delta=0$ (as when $F$ is strictly concave), the speed that $\{P_{\eta/\sqrt n}\}$ approaches the null is too fast in the sense that $\sqrt n \phi(F_n)\to 0$, thereby making the test hard to reject.

\subsection{Comparisons with Existing Tests}\label{Sec: comparison}

We now compare our concavity test with some existing tests. While there is a rich literature on testing monotonicity in settings such as nonparametric regression---see \citet{Chetverikov2018Monotonicity} for a recent study with a brief survey, it is somewhat surprising that results for the same problem in the density context seem rather limited. This is just one piece of evidence consistent with Jon Wellner's view that ``the (shape-constrained) community also needs to do more work to provide inferential methods beyond estimation'' \citep{BanerjeeSamworth2018Interview}.

\citet{Banerjee_Wellner2001} propose likelihood ratio tests for equality at a fixed point, rather than the global shape as we consider, in the setting of nonparametric estimation of a monotone function---see also \citet{Banerjee2005LR}. More related to our setting is the work by \citet{Woodroofe_Sun1999}, who are concerned with testing uniformity versus a monotone (but not uniform) density. Translated to our setup, they study the simple null $F$ being the uniform distribution against the composite alternative that $F$ is convex (but is not a uniform distribution). Therefore, the parameter spaces under the null and the alternative are smaller than the corresponding spaces in our setup. \citet{Woodroofe_Sun1999} propose two tests, namely, the $P$-test and the $D$-test, based on a penalized nonparametric maximum likelihood estimator of the density. The $P$-test statistic is simply the penalized likelihood ratio, while the $D$-test statistic is $\sqrt n$---with $n$ being the sample size---times the uniform distance between the penalized cdf estimator and the uniform cdf.

The simple null hypothesis in \citet{Woodroofe_Sun1999}, though greatly simplifies the asymptotic analysis, may be restrictive from a practical point of view. This is a reflection of the fact that the pointwise asymptotic distributions of the test statistic under a composite null such as ours are highly nonstandard. In a nonparametric regression setting, \citet{Durot2003testing} studies the composite null that the regression function is nonincreasing against the alternative that it is not. The test statistic proposed by Durot is a rescaled version of the uniform distance between a cumulative regression estimator and its least concave majorant. \citet{Durot2003testing} shows that the constant regression functions are asymptotically least favorable, base on which critical values are constructed. \citet{Kulikov_Lopuhaa2004testing} adapt \citet{Durot2003testing}'s test to the density setup based on the same critical values---see their Theorem 3.1, and also Proposition 3.1 in \citet{Kulikov_Lopuhaa2008}. \citet{Kulikov_Lopuhaa2004testing} also propose a test based on the $L_p$-norm (against the empirical distribution) of the empirical cdf and its least concave majorant, also relying on critical values constructed from the least favorable distribution. We note that the test based on the statistic $T_n$ in \citet{Kulikov_Lopuhaa2004testing} is invalid in the sense that it may over-reject under the null, as pointed out by the same authors and confirmed by their simulations.

Tests based on the least favorable distributions, though control the size in a very simple way, may be too conservative and result in power loss. This includes the tests of \citet{Durot2003testing} and \citet{Kulikov_Lopuhaa2004testing}. On the other hand, the $P$-test and the $D$-test explicitly take into account the alternative hypothesis of nondecreasing densities. Therefore, they are powerful in detecting nondecreasing densities, but may perform poorly when the underlying distribution is neither concave nor convex. Our simulations confirm these predictions. Finally, we reiterate that the need of bootstrap for testing concavity as well as the tuning parameters $\kappa_n$ and $t_n$ is in line with the nonstandard nature of the problem, rather than a special attribute of our inferential framework. The implementation of our bootstrap is as simple as calculating the test statistic, or as simple as computing the least favorable majorant.

\section{Simulation Experiments}\label{Sec: simulation}

We next evaluate the finite sample performance of our test through Monte Carlo simulations. Special attention shall be paid to the tuning parameters $\kappa_n$ (for the KW selection) and $t_n$ (for the rescaled bootstrap). For this, we consider $\kappa_n\in\{n^{-1/7},n^{-1/8}, \log^{-1}n\}$ and $t_n\in\{n^{-1/\varpi}: \varpi=3,4,\ldots,7\}$. We run two sets of simulations in order to compare with the standard bootstrap and some existing monotonicity tests. Throughout, we let the significance level be $5\%$, and all results are based on 1000 Monte Carlo simulations and 500 bootstrap repetitions for each simulation replication (to implement our test).

First, we consider the distributions defined by the following densities supported on $\mathbf R_+$: for any $x\in\mathbf R_+$,
\begin{align}
f_1(x) & =
\frac{1}{\sqrt{\pi/2}}\exp\{-\frac{x^2}{2}\} ~,\label{Eqn: sim1 f1} \\
f_2(x) & =\begin{cases}
3/2 & \text{ if }0\le x<1/2\\
1/2 & \text{ if }1/2\le x<9/10\\
\exp\{9/10-x\}/20 &\text{ if }x\ge 9/10
\end{cases}~.\label{Eqn: sim1 f2}
\end{align}
The function $f_1$ is strictly decreasing, while $f_2$ is weakly decreasing though strictly decreasing on the region $[9/10,\infty)$. Hence, the corresponding distribution functions, denoted $F_1$ and $F_2$, are strictly concave and non-strictly concave respectively. To examine local behaviors of our test, we construct differentiable paths passing through $F_1$ and $F_2$ as follows: for $t\in\mathbf R_+$ and $x\in\mathbf R_+$,
\begin{align}
f_{1,t}(x)&= \frac{\exp\{-tx\}f_1(x)}{\int_{\mathbf R_+} \exp\{-tx\}f_1(x) dx} =
\frac{1}{a(t)}\exp\{-\frac{(x+t)^2}{2}\} ~,\label{Eqn: sim1 f1t} \\
f_{2,t}(x)&=\begin{cases}
e^{-t}3/2 & \text{ if }0\le x<1/2\\
e^t/2 & \text{ if }1/2\le x<9/10\\
b(t)\exp\{9/10-x\}/20 &\text{ if }x\ge 9/10
\end{cases}~,\label{Eqn: sim1 f2t}
\end{align}
where $a(t)\equiv\int_t^\infty \exp\{-u^2/2\}du$ and $b(t)\equiv 20-15e^{-t}-4e^t$. By Example 3.2.1 and Proposition 2.1.1 in \citet{BKRW993Efficient}, these two paths are indeed differentiable (under the null). We then generate i.i.d.\ samples $\{X_i\}_{i=1}^n$ from \eqref{Eqn: sim1 f1}, \eqref{Eqn: sim1 f2}, \eqref{Eqn: sim1 f1t} and \eqref{Eqn: sim1 f2t}, with $n\in\{100,200,300,400,1000\}$. For simplicity, we only consider the sup statistic $\sqrt n\|\hat{\mathbb F}_n-\mathbb F_n\|_\infty$, and report results based on the rescaled and the standard bootstrap.

{\renewcommand{\arraystretch}{1.45}
\begin{table}[!h]
\caption{Pointwise Size Control with Strict Concavity Based on \eqref{Eqn: sim1 f1}}\label{Tab: size pointwise, sc}
\begin{footnotesize}
\begin{center}
\begin{tabular}{cccccccc}
\hline\hline
 &   & \multicolumn{5}{c}{Rescaled} & \multirow{3}{*}[-0.5em]{Standard}\\
\cmidrule{3-7}
 &   & \multicolumn{5}{c}{$t_n$}    & \\
\cmidrule{3-7}
$n$ & $\kappa_n$                     & $n^{-1/3}$ & $n^{-1/4}$ & $n^{-1/5}$ & $n^{-1/6}$ & $n^{-1/7}$ & \\
 \hline
 \multirow{3}{*}{ 100 } & $n^{-1/7}$ & $0.017$    & $0.008$    & $0.006$    & $0.005$    & $0.004$ & \multirow{3}{*}{ $0.053$ }\\
                        & $n^{-1/8}$ & $0.017$    & $0.008$    & $0.006$    & $0.005$    & $0.005$ &   \\
                  & $(\log n)^{-1}$  & $0.016$    & $0.008$    & $0.006$    & $0.005$    & $0.004$ &   \\ 
                  \cmidrule{3-8}
 \multirow{3}{*}{ 200 } & $n^{-1/7}$ & $0.022$    & $0.015$    & $0.011$    & $0.009$    & $0.009$ & \multirow{3}{*}{ $0.071$ }\\
                        & $n^{-1/8}$ & $0.021$    & $0.015$    & $0.011$    & $0.009$    & $0.009$ &   \\
                  & $(\log n)^{-1}$  & $0.021$    & $0.015$    & $0.011$    & $0.009$    & $0.009$ &   \\ 
                  \cmidrule{3-8}
 \multirow{3}{*}{ 300 } & $n^{-1/7}$ & $0.010$    & $0.002$    & $0.002$    & $0.003$    & $0.003$ & \multirow{3}{*}{ $0.064$ }\\
                        & $n^{-1/8}$ & $0.011$    & $0.002$    & $0.002$    & $0.003$    & $0.002$ &   \\
                  & $(\log n)^{-1}$  & $0.010$    & $0.002$    & $0.002$    & $0.003$    & $0.003$ &   \\ 
                  \cmidrule{3-8}
 \multirow{3}{*}{ 400 } & $n^{-1/7}$ & $0.013$    & $0.006$    & $0.000$    & $0.000$    & $0.000$ & \multirow{3}{*}{ $0.068$ }\\
                        & $n^{-1/8}$ & $0.014$    & $0.005$    & $0.000$    & $0.000$    & $0.000$ &   \\
                  & $(\log n)^{-1}$  & $0.014$    & $0.005$    & $0.000$    & $0.000$    & $0.000$ &   \\ 
                  \cmidrule{3-8}
 \multirow{3}{*}{ 1000 }& $n^{-1/7}$ & $0.009$    & $0.002$    & $0.003$    & $0.000$    & $0.000$ & \multirow{3}{*}{ $0.078$ }\\
                        & $n^{-1/8}$ & $0.009$    & $0.002$    & $0.003$    & $0.000$    & $0.000$ &   \\
                  & $(\log n)^{-1}$  & $0.009$    & $0.002$    & $0.003$    & $0.000$    & $0.000$ &   \\ 
\hline\hline
\end{tabular}
\end{center}
\end{footnotesize}
\end{table}
}

Tables \ref{Tab: size pointwise, sc}, \ref{Tab: size pointwise, nsc}, \ref{Tab: size local, sc} and \ref{Tab: size local, nsc} present the numerical results. Inconsistency of the standard bootstrap is prominently evidenced in last columns of Tables \ref{Tab: size pointwise, nsc} and \ref{Tab: size local, nsc} which are based on non-strictly concave functions (so that the operator $\mathcal M$ are only Hadamard directionally differentiable at these distribution functions). On the contrary, our test alleviates the size distortion both pointwise and locally, though the results vary with the choices of the tuning parameters. Overall, the choices $t_n\in\{n^{-1/3},n^{-1/4}\}$ tend to be too small, while $t_n\in\{n^{-1/5},n^{-1/6},n^{-1/7}\}$ are adequate in that they control the size well. Tables \ref{Tab: size pointwise, sc} and \ref{Tab: size local, sc} also show that, when the distribution function is strictly concave, our test by design controls the size both pointwise and locally, though the rejection rates are often close to zero, which is consistent with Theorem \ref{Thm: test concavity, local}. Moreover, as in Tables \ref{Tab: size pointwise, sc} and \ref{Tab: size local, sc}, the standard bootstrap also exhibits size control, due to the facts that the LCM operator is fully (Hadamard) differentiable at strictly concave distribution functions and that the distance between the Grenander distribution estimator and the empirical cdf converges faster than $\sqrt n$ by the Kiefer-Wolfowitz theorems. Finally, we note that the results are quite insensitive to the choice of $\kappa_n$.

{\renewcommand{\arraystretch}{1.5}
\begin{table}[!h]
\caption{Pointwise Size Control with Non-Strict Concavity Based on \eqref{Eqn: sim1 f2}}\label{Tab: size pointwise, nsc}
\begin{footnotesize}
\begin{center}
\begin{tabular}{cccccccc}
\hline\hline
 &   & \multicolumn{5}{c}{Rescaled} & \multirow{3}{*}[-0.5em]{Standard}\\
\cmidrule{3-7}
 &   & \multicolumn{5}{c}{$t_n$}    & \\
\cmidrule{3-7}
$n$ & $\kappa_n$                     & $n^{-1/3}$ & $n^{-1/4}$ & $n^{-1/5}$ & $n^{-1/6}$ & $n^{-1/7}$ & \\
 \hline
 \multirow{3}{*}{ 100 } & $n^{-1/7}$ & $0.068$    & $0.051$    & $0.046$    & $0.041$    & $0.036$ & \multirow{3}{*}{ $0.132$ }\\
                        & $n^{-1/8}$ & $0.069$    & $0.049$    & $0.044$    & $0.041$    & $0.036$ &   \\
                  & $(\log n)^{-1}$  & $0.070$    & $0.052$    & $0.045$    & $0.041$    & $0.036$ &   \\ 
                  \cmidrule{3-8}
 \multirow{3}{*}{ 200 } & $n^{-1/7}$ & $0.082$    & $0.061$    & $0.052$    & $0.046$    & $0.045$ & \multirow{3}{*}{ $0.151$ }\\
                        & $n^{-1/8}$ & $0.084$    & $0.063$    & $0.054$    & $0.045$    & $0.043$ &   \\
                  & $(\log n)^{-1}$  & $0.080$    & $0.063$    & $0.053$    & $0.046$    & $0.044$ &   \\ 
                  \cmidrule{3-8}
 \multirow{3}{*}{ 300 } & $n^{-1/7}$ & $0.086$    & $0.066$    & $0.056$    & $0.046$    & $0.040$ & \multirow{3}{*}{ $0.169$ }\\
                        & $n^{-1/8}$ & $0.085$    & $0.068$    & $0.056$    & $0.047$    & $0.040$ &   \\
                  & $(\log n)^{-1}$  & $0.086$    & $0.067$    & $0.056$    & $0.046$    & $0.040$ &   \\ 
                  \cmidrule{3-8}
 \multirow{3}{*}{ 400 } & $n^{-1/7}$ & $0.109$    & $0.081$    & $0.064$    & $0.057$    & $0.053$ & \multirow{3}{*}{ $0.191$ }\\
                        & $n^{-1/8}$ & $0.106$    & $0.079$    & $0.064$    & $0.057$    & $0.050$ &   \\
                  & $(\log n)^{-1}$  & $0.108$    & $0.079$    & $0.063$    & $0.058$    & $0.051$ &   \\ 
                  \cmidrule{3-8}
 \multirow{3}{*}{ 1000 }& $n^{-1/7}$ & $0.096$    & $0.071$    & $0.063$    & $0.055$    & $0.046$ & \multirow{3}{*}{ $0.205$ }\\
                        & $n^{-1/8}$ & $0.095$    & $0.072$    & $0.063$    & $0.054$    & $0.046$ &   \\
                  & $(\log n)^{-1}$  & $0.093$    & $0.073$    & $0.062$    & $0.056$    & $0.046$ &   \\ 
\hline\hline
\end{tabular}
\end{center}
\end{footnotesize}
\end{table}
}

{\renewcommand{\arraystretch}{1.5}
\begin{table}[!h]
\caption{Local Size Control with Strict Concavity Based on $f_{1,1/\sqrt n}$ from \eqref{Eqn: sim1 f1t}}\label{Tab: size local, sc}
\begin{footnotesize}
\begin{center}
\begin{tabular}{cccccccc}
\hline\hline
 &   & \multicolumn{5}{c}{Rescaled} & \multirow{3}{*}[-0.5em]{Standard}\\
\cmidrule{3-7}
 &   & \multicolumn{5}{c}{$t_n$}    & \\
\cmidrule{3-7}
$n$ & $\kappa_n$                     & $n^{-1/3}$ & $n^{-1/4}$ & $n^{-1/5}$ & $n^{-1/6}$ & $n^{-1/7}$ & \\
 \hline
 \multirow{3}{*}{ 100 } & $n^{-1/7}$ & $0.017$    & $0.009$    & $0.007$    & $0.008$    & $0.006$ & \multirow{3}{*}{ $0.050$ }\\
                        & $n^{-1/8}$ & $0.016$    & $0.010$    & $0.007$    & $0.007$    & $0.006$ &   \\
                  & $(\log n)^{-1}$  & $0.017$    & $0.008$    & $0.008$    & $0.007$    & $0.006$ &   \\ 
                  \cmidrule{3-8}
 \multirow{3}{*}{ 200 } & $n^{-1/7}$ & $0.018$    & $0.006$    & $0.004$    & $0.003$    & $0.003$ & \multirow{3}{*}{ $0.060$ }\\
                        & $n^{-1/8}$ & $0.018$    & $0.006$    & $0.004$    & $0.003$    & $0.003$ &   \\
                  & $(\log n)^{-1}$  & $0.018$    & $0.006$    & $0.004$    & $0.003$    & $0.003$ &   \\ 
                  \cmidrule{3-8}
 \multirow{3}{*}{ 300 } & $n^{-1/7}$ & $0.015$    & $0.009$    & $0.006$    & $0.002$    & $0.002$ & \multirow{3}{*}{ $0.059$ }\\
                        & $n^{-1/8}$ & $0.014$    & $0.009$    & $0.006$    & $0.002$    & $0.002$ &   \\
                  & $(\log n)^{-1}$  & $0.016$    & $0.009$    & $0.006$    & $0.002$    & $0.002$ &   \\ 
                  \cmidrule{3-8}
 \multirow{3}{*}{ 400 } & $n^{-1/7}$ & $0.015$    & $0.002$    & $0.003$    & $0.000$    & $0.000$ & \multirow{3}{*}{ $0.070$ }\\
                        & $n^{-1/8}$ & $0.015$    & $0.002$    & $0.003$    & $0.000$    & $0.000$ &   \\
                  & $(\log n)^{-1}$  & $0.015$    & $0.002$    & $0.003$    & $0.000$    & $0.000$ &   \\ 
                  \cmidrule{3-8}
 \multirow{3}{*}{ 1000 }& $n^{-1/7}$ & $0.009$    & $0.002$    & $0.002$    & $0.002$    & $0.002$ & \multirow{3}{*}{ $0.070$ }\\
                        & $n^{-1/8}$ & $0.010$    & $0.002$    & $0.002$    & $0.002$    & $0.002$ &   \\
                  & $(\log n)^{-1}$  & $0.011$    & $0.002$    & $0.000$    & $0.002$    & $0.002$ &   \\ 
\hline\hline
\end{tabular}
\end{center}
\end{footnotesize}
\end{table}
}

{\renewcommand{\arraystretch}{1.5}
\begin{table}[!h]
\caption{Local Size Control with Non-Strict Concavity Based on $f_{2,1/\sqrt n}$ from \eqref{Eqn: sim1 f2t}}\label{Tab: size local, nsc}
\begin{footnotesize}
\begin{center}
\begin{tabular}{cccccccc}
\hline\hline
 &   & \multicolumn{5}{c}{Rescaled} & \multirow{3}{*}[-0.5em]{Standard}\\
\cmidrule{3-7}
 &   & \multicolumn{5}{c}{$t_n$}    & \\
\cmidrule{3-7}
$n$ & $\kappa_n$                     & $n^{-1/3}$ & $n^{-1/4}$ & $n^{-1/5}$ & $n^{-1/6}$ & $n^{-1/7}$ & \\
 \hline
 \multirow{3}{*}{ 100 } & $n^{-1/7}$ & $0.070$    & $0.052$    & $0.043$    & $0.041$    & $0.036$ & \multirow{3}{*}{ $0.128$ }\\
                        & $n^{-1/8}$ & $0.070$    & $0.051$    & $0.042$    & $0.038$    & $0.038$ &   \\
                  & $(\log n)^{-1}$  & $0.073$    & $0.051$    & $0.043$    & $0.038$    & $0.034$ &   \\ 
                  \cmidrule{3-8}
 \multirow{3}{*}{ 200 } & $n^{-1/7}$ & $0.112$    & $0.074$    & $0.063$    & $0.051$    & $0.045$ & \multirow{3}{*}{ $0.193$ }\\
                        & $n^{-1/8}$ & $0.111$    & $0.076$    & $0.063$    & $0.049$    & $0.045$ &   \\
                  & $(\log n)^{-1}$  & $0.113$    & $0.075$    & $0.062$    & $0.049$    & $0.042$ &   \\ 
                  \cmidrule{3-8}
 \multirow{3}{*}{ 300 } & $n^{-1/7}$ & $0.098$    & $0.066$    & $0.055$    & $0.050$    & $0.044$ & \multirow{3}{*}{ $0.182$ }\\
                        & $n^{-1/8}$ & $0.099$    & $0.067$    & $0.058$    & $0.049$    & $0.044$ &   \\
                  & $(\log n)^{-1}$  & $0.100$    & $0.070$    & $0.057$    & $0.050$    & $0.044$ &   \\ 
                  \cmidrule{3-8}
 \multirow{3}{*}{ 400 } & $n^{-1/7}$ & $0.113$    & $0.084$    & $0.071$    & $0.061$    & $0.055$ & \multirow{3}{*}{ $0.197$ }\\
                        & $n^{-1/8}$ & $0.113$    & $0.085$    & $0.070$    & $0.062$    & $0.055$ &   \\
                  & $(\log n)^{-1}$  & $0.112$    & $0.083$    & $0.069$    & $0.061$    & $0.056$ &   \\ 
                  \cmidrule{3-8}
\multirow{3}{*}{ 1000 } & $n^{-1/7}$ & $0.092$    & $0.064$    & $0.055$    & $0.046$    & $0.036$ & \multirow{3}{*}{ $0.196$ }\\
                        & $n^{-1/8}$ & $0.093$    & $0.061$    & $0.055$    & $0.044$    & $0.036$ &   \\
                  & $(\log n)^{-1}$  & $0.094$    & $0.066$    & $0.055$    & $0.045$    & $0.037$ &   \\ 
\hline\hline
\end{tabular}
\end{center}
\end{footnotesize}
\end{table}
}

\begin{figure}
\centering
\subfloat[Starting from the top are cdfs in \eqref{Eqn: cdf compare1} corresponding to $\lambda= -5, -1,-0.2,0.2,1,5$]{
\begin{tikzpicture}[scale=0.75,baseline]
\begin{axis}[domain=0:1, samples=100, restrict y to domain=0:1,
cycle list={%
{NavyBlue,solid},
{Azure4,densely dashdotdotted},
{Green3, densely dotted},
{Green3, densely dotted},
{Azure4,densely dashdotdotted},
{NavyBlue,solid}
}
]
\foreach \theta in {-5,-1, -0.2, 0.2, 1, 5}{%
\addplot+[domain=0:1] {(exp(\theta*\x)-1)/(exp(\theta)-1)};
}
\end{axis}
\end{tikzpicture}}
\quad
\subfloat[Starting from the top are cdfs in \eqref{Eqn: cdf compare2} corresponding to $\lambda= 0.9, 0.6, 0.3$]{
\begin{tikzpicture}[scale=0.75,baseline]
\begin{axis}[domain=0:1, samples=100,restrict y to domain=0:1,
cycle list={%
{Green3, densely dotted},
{Azure4,densely dashdotdotted},
{NavyBlue,solid}
}
]
\foreach \theta in {0.3, 0.6, 0.9}{%
\addplot+[domain=0:1] {(\x<\theta)*(-(\x-\theta)^2/\theta+\theta)+(\x>=\theta)*\x}; 
}
\end{axis}
\end{tikzpicture}}
\caption{Graphical representations of the cdfs in \eqref{Eqn: cdf compare1} and \eqref{Eqn: cdf compare2}}\label{Fig: cdfs in simulations}
\end{figure}
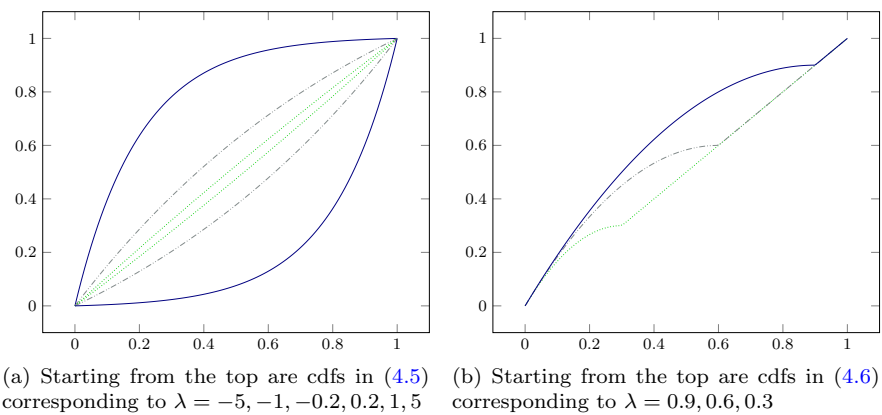

Next, we compare our test with some existing tests. As before, the null hypothesis is that the true cdf is conave (or, equivalently, the true pdf is nonincreasing), while the alternative is that it is not. We shall consider the $D$ and the $P$ tests in \citet{Woodroofe_Sun1999}, the Kolmogorov-Smirnov type test in \citet{Durot2003testing} and \citet{Kulikov_Lopuhaa2004testing}, labelled as the KS test, and the test based on the $R_n$ statistic with $k=1$ in \citet{Kulikov_Lopuhaa2004testing}, labelled as the KL test. We design two families of distributions. The first family consists of distributions defined by: for $\lambda\in\mathbf R$,
\begin{align}\label{Eqn: cdf compare1}
F_\lambda(x) = \begin{cases}
0 & \text{ if }x<0\\
\frac{e^{\lambda x} -1 }{e^{\lambda}- 1} & \text{ if }0\le x\le 1 \\
1 & \text{ if }x>1
\end{cases}~,
\end{align}
where $F_\lambda$ is understood to be the uniform cdf if $\lambda=0$. Clearly, $F_\lambda$ is concave and hence in the null whenever $\lambda\le 0$, but is convex and hence in the alternative if $\lambda>0$---see Figure \ref{Fig: cdfs in simulations}-(a). Thus, this is a setup well suited for the $D$ test and the $P$ test proposed by \citet{Woodroofe_Sun1999}, because the alternative consists of only convex cdfs (or nondecreasing pdfs). As mentioned in Section \ref{Sec: comparison}, the $D$ and the $P$ tests may perform poorly in detecting distributions that are neither convex nor concave. To verify this point numerically, we consider a family of distributions (that belongs to the alternative) defined by: for $\lambda \in (0,1)$,
\begin{align}\label{Eqn: cdf compare2}
F_\lambda(x) = \begin{cases}
0 & \text{ if }x<0\\
-(x-\lambda)^2/\lambda+\lambda & \text{ if }0\le x<\lambda \\
x & \text{ if }\lambda\le x<1\\
1 & \text{ if }\lambda\ge 1
\end{cases}~.
\end{align}
As shown in Figure \ref{Fig: cdfs in simulations}-(b), $F_\lambda$ defined in \eqref{Eqn: cdf compare2} is neither concave nor convex on $[0,1]$ for any $\lambda\in(0,1)$, though it is concave on $[0,\lambda]$ and $[\lambda,1]$.

As suggested by \citet{Woodroofe_Sun1999}, we choose the penalty parameter $c$ for the $D$ and the $P$ tests to be $0.2$, and the corresponding critical values may be found in their Table 1. The critical values for the KS and the KL tests are given in Table 1 in \citet{Kulikov_Lopuhaa2004testing}---see also Table 1 in \citet{Durot2003testing}. For convenience, if $(t_n,\kappa_n)=(n^{-1/3},n^{1-/7})$, then we label our test as R37---``R'' stands for ``the rescaled bootstrap''; if $(t_n,\kappa_n)=(n^{-1/3},\log^{-1}n)$, then the resulting test is labelled as R3L---``L'' stands for ``logarithm''; so on and so forth. Below, we try more sample sizes with $n\in\{100,200,300,400,600,800,1000\}$, but only report results for $n\in\{100,400\}$ in the main text and relegate the rest to Appendix \ref{App: Simulation}.

{\renewcommand{\arraystretch}{1.5}
\setlength{\tabcolsep}{5pt}
\begin{table}[!h]
\caption{Comparisons with Existing Tests Based on \eqref{Eqn: cdf compare1}: $n=100$}\label{Tab: comparison, n100}
\begin{footnotesize}
\begin{center}
\begin{tabular}{ccccccccccc}
\hline\hline
  & \multicolumn{10}{c}{$\lambda$}  \\
\cmidrule{2-11}
 & $-1$ & $-0.1$ & $-0.05$ & $0$ & $0.05$ & $0.1$ & $0.15$ & $0.2$ & $0.25$ & $1$ \\
\hline
R37    & $0.010$ & $0.084$ & $0.106$ & $0.095$ & $0.124$ &  $0.132$ & $0.181$ & $0.194$ & $0.235$ & $0.819$ \\
R38    & $0.009$ & $0.085$ & $0.107$ & $0.094$ & $0.123$ & $0.133$ & $0.180$ & $0.187$ & $0.235$ & $0.818$ \\
R3L    & $0.009$ & $0.085$ & $0.107$ & $0.096$ & $0.123$ & $0.135$ & $0.180$ & $0.189$ & $0.236$ & $0.818$ \\
R47    & $0.005$ & $0.062$ & $0.084$ & $0.079$ & $0.103$ & $0.109$ & $0.157$ & $0.165$ & $0.206$ & $0.799$ \\
R48    & $0.004$ & $0.064$ & $0.086$ & $0.080$ & $0.104$ & $0.111$ & $0.156$ & $0.165$ & $0.208$ & $0.799$ \\
R4L    & $0.004$ & $0.061$ & $0.088$ & $0.079$ & $0.102$ & $0.110$ & $0.157$ & $0.165$ & $0.209$ & $0.798$ \\
R57    & $0.003$ & $0.049$ & $0.079$ & $0.076$ & $0.094$ & $0.097$ & $0.149$ & $0.152$ & $0.194$ & $0.789$ \\
R58    & $0.003$ & $0.051$ & $0.078$ & $0.075$ & $0.094$ & $0.099$ & $0.146$ & $0.153$ & $0.195$ & $0.789$ \\
R5L    & $0.003$ & $0.050$ & $0.079$ & $0.075$ & $0.095$ & $0.097$ & $0.145$ & $0.154$ & $0.195$ & $0.789$ \\
R67    & $0.002$ & $0.049$ & $0.075$ & $0.075$ & $0.091$ & $0.091$ & $0.137$ & $0.144$ & $0.187$ & $0.784$ \\
R68    & $0.002$ & $0.049$ & $0.076$ & $0.075$ & $0.091$ & $0.091$ & $0.137$ & $0.145$ & $0.187$ & $0.784$ \\
R6L    & $0.002$ & $0.047$ & $0.075$ & $0.073$ & $0.091$ & $0.090$ & $0.144$ & $0.146$ & $0.188$ & $0.787$ \\
R77    & $0.002$ & $0.047$ & $0.073$ & $0.071$ & $0.087$ & $0.088$ & $0.137$ & $0.143$ & $0.180$ & $0.782$ \\
R78    & $0.003$ & $0.045$ & $0.074$ & $0.071$ & $0.088$ & $0.090$ & $0.139$ & $0.145$ & $0.182$ & $0.783$ \\
R7L    & $0.002$ & $0.045$ & $0.074$ & $0.072$ & $0.089$ & $0.088$ & $0.141$ & $0.143$ & $0.182$ & $0.784$ \\ \hline
P-Test & $0.003$ & $0.032$ & $0.030$ & $0.037$ & $0.064$ & $0.068$ & $0.072$ & $0.079$ & $0.119$ & $0.764$ \\
D-Test & $0.000$ & $0.018$ & $0.039$ & $0.036$ & $0.061$ & $0.071$ & $0.090$ & $0.108$ & $0.170$ & $0.809$ \\
KS-Test & $0.000$ & $0.010$ & $0.017$ & $0.018$ & $0.031$ & $0.035$ & $0.048$ & $0.057$ & $0.083$ & $0.642$ \\
KL-Test & $0.000$ & $0.012$ & $0.020$ & $0.024$ & $0.038$ & $0.042$ & $0.052$ & $0.074$ & $0.109$ & $0.765$ \\
\hline\hline
\end{tabular}
\end{center}
\end{footnotesize}
\end{table}
}

{\renewcommand{\arraystretch}{1.5}
\setlength{\tabcolsep}{5pt}
\begin{table}[!ht]
\caption{Comparisons with Existing Tests Based on \eqref{Eqn: cdf compare1}: $n=400$}\label{Tab: comparison, n400}
\begin{footnotesize}
\begin{center}
\begin{tabular}{ccccccccccc}
\hline\hline
  & \multicolumn{10}{c}{$\lambda$}  \\
\cmidrule{2-11}
 & $-1$ & $-0.1$ & $-0.05$ & $0$ & $0.05$ & $0.1$ & $0.15$ & $0.2$ & $0.25$ & $1$ \\
\hline
R37  & $0.004$ & $0.052$ & $0.077$ & $0.086$ & $0.134$ & $0.189$ & $0.286$ & $0.327$ & $0.420$ & $0.999$ \\
R38  & $0.004$ & $0.053$ & $0.074$ & $0.088$ & $0.134$ & $0.191$ & $0.287$ & $0.326$ & $0.420$ & $0.999$ \\
R3L  & $0.004$ & $0.052$ & $0.075$ & $0.087$ & $0.133$ & $0.191$ & $0.287$ & $0.324$ & $0.422$ & $0.999$ \\
R47  & $0.003$ & $0.033$ & $0.058$ & $0.068$ & $0.107$ & $0.163$ & $0.252$ & $0.296$ & $0.382$ & $0.999$ \\
R48  & $0.003$ & $0.036$ & $0.059$ & $0.065$ & $0.107$ & $0.166$ & $0.252$ & $0.292$ & $0.384$ & $0.999$ \\
R4L  & $0.002$ & $0.035$ & $0.059$ & $0.067$ & $0.107$ & $0.164$ & $0.251$ & $0.296$ & $0.385$ & $0.999$ \\
R57  & $0.003$ & $0.027$ & $0.050$ & $0.061$ & $0.096$ & $0.150$ & $0.237$ & $0.282$ & $0.364$ & $0.999$ \\
R58  & $0.003$ & $0.029$ & $0.053$ & $0.062$ & $0.097$ & $0.150$ & $0.237$ & $0.279$ & $0.365$ & $0.999$ \\
R5L  & $0.003$ & $0.029$ & $0.051$ & $0.061$ & $0.095$ & $0.153$ & $0.237$ & $0.278$ & $0.364$ & $0.999$ \\
R67  & $0.000$ & $0.027$ & $0.048$ & $0.060$ & $0.092$ & $0.144$ & $0.227$ & $0.272$ & $0.349$ & $0.999$ \\
R68  & $0.000$ & $0.024$ & $0.046$ & $0.052$ & $0.080$ & $0.151$ & $0.226$ & $0.276$ & $0.368$ & $0.999$ \\
R6L  & $0.000$ & $0.027$ & $0.050$ & $0.058$ & $0.093$ & $0.146$ & $0.230$ & $0.273$ & $0.354$ & $0.999$ \\
R77  & $0.000$ & $0.023$ & $0.048$ & $0.055$ & $0.087$ & $0.140$ & $0.224$ & $0.268$ & $0.347$ & $0.999$ \\
R78  & $0.000$ & $0.026$ & $0.050$ & $0.053$ & $0.087$ & $0.143$ & $0.223$ & $0.266$ & $0.345$ & $0.999$ \\
R7L  & $0.003$ & $0.024$ & $0.048$ & $0.056$ & $0.088$ & $0.140$ & $0.224$ & $0.265$ & $0.346$ & $0.999$ \\
 \hline
P-Test & $0.000$ & $0.024$ & $0.035$ & $0.041$ & $0.078$ & $0.105$ & $0.189$ & $0.234$ & $0.343$ & $0.999$ \\
D-Test & $0.000$ & $0.013$ & $0.023$ & $0.037$ & $0.065$ & $0.135$ & $0.219$ & $0.294$ & $0.382$ & $0.999$ \\
KS-Test& $0.000$ & $0.010$ & $0.020$ & $0.025$ & $0.046$ & $0.080$ & $0.145$ & $0.198$ & $0.259$ & $0.999$ \\
KL-Test& $0.000$ & $0.009$ & $0.018$ & $0.037$ & $0.060$ & $0.102$ & $0.184$ & $0.258$ & $0.340$ & $1.000$ \\
\hline\hline
\end{tabular}
\end{center}
\end{footnotesize}
\end{table}
}

Tables \ref{Tab: comparison, n100} and \ref{Tab: comparison, n400} record the results for $n\in\{100, 400\}$. Across the choices of $(t_n,\kappa_n)$, when $F_\lambda$ is close to being uniform, our test tends to over-reject (under the null) with smaller $t_n$, especially with $t_n=n^{-1/3}$ and in small samples. All other tests control size well in large samples, but tend to under-reject in small samples. In terms of power, the KS test appears to be the least powerful against local alternatives (i.e., distributions with positive but small $\lambda$). The $D$-test overall appears to be the most powerful one, though the power is comparable to ours and that of the $P$ and the KL tests. For a fair comparison, we also compute the size-adjusted power by subtracting the empirical rejection rates for $\lambda=0$ (the least favorable case) from the corresponding rejections rates for $\lambda>0$. The results are recorded in Table \ref{Tab: comparison}, from which we see that the power patterns more or less remain. Those results for $n\in\{200,300,600,800,1000\}$ share similar patterns---see Appendix \ref{App: Simulation}.

{\renewcommand{\arraystretch}{1.5}
\setlength{\tabcolsep}{4pt}
\begin{table}[!h]
\caption{Size-Adjusted Power Comparisons with Existing Tests Based on \eqref{Eqn: cdf compare1}}\label{Tab: comparison}
\begin{scriptsize}
\begin{center}
\begin{tabular}{cccccccccccccc}
\hline\hline
  & \multicolumn{6}{c}{$\lambda$ for $n=100$} & & \multicolumn{6}{c}{$\lambda$ for $n=400$}\\
\cmidrule{2-7}  \cmidrule{9-14}
 & $0.05$ & $0.1$ & $0.15$ & $0.2$ & $0.25$ & $1$ & & $0.05$ & $0.1$ & $0.15$ & $0.2$ & $0.25$ & $1$  \\
\hline
R37                            & $0.029$ & $0.037$ & $0.086$ & $0.099$ & $0.140$ & $0.724$ && $0.048$ & $0.103$ & $0.200$ & $0.241$ & $0.334$ & $0.913$ \\
R38                            & $0.029$ & $0.039$ & $0.086$ & $0.093$ & $0.141$ & $0.724$ && $0.046$ & $0.103$ & $0.199$ & $0.238$ & $0.332$ & $0.911$ \\
R3L                            & $0.027$ & $0.039$ & $0.084$ & $0.093$ & $0.140$ & $0.722$ && $0.046$ & $0.104$ & $0.200$ & $0.237$ & $0.335$ & $0.912$ \\
R47                            & $0.024$ & $0.030$ & $0.078$ & $0.086$ & $0.127$ & $0.720$ && $0.039$ & $0.095$ & $0.184$ & $0.228$ & $0.314$ & $0.931$ \\
R48                            & $0.024$ & $0.031$ & $0.076$ & $0.085$ & $0.128$ & $0.719$ && $0.042$ & $0.101$ & $0.187$ & $0.227$ & $0.319$ & $0.934$ \\
R4L                            & $0.023$ & $0.031$ & $0.078$ & $0.086$ & $0.130$ & $0.719$ && $0.040$ & $0.097$ & $0.184$ & $0.229$ & $0.318$ & $0.932$ \\
R57                            & $0.018$ & $0.021$ & $0.073$ & $0.076$ & $0.118$ & $0.713$ && $0.035$ & $0.089$ & $0.176$ & $0.221$ & $0.303$ & $0.938$ \\
R58                            & $0.019$ & $0.024$ & $0.071$ & $0.078$ & $0.120$ & $0.714$ && $0.035$ & $0.088$ & $0.175$ & $0.217$ & $0.303$ & $0.937$ \\
R5L                            & $0.020$ & $0.022$ & $0.070$ & $0.079$ & $0.120$ & $0.714$ && $0.034$ & $0.092$ & $0.176$ & $0.217$ & $0.303$ & $0.938$ \\
R67                            & $0.016$ & $0.016$ & $0.062$ & $0.069$ & $0.112$ & $0.709$ && $0.032$ & $0.084$ & $0.167$ & $0.212$ & $0.289$ & $0.939$ \\
R68                            & $0.016$ & $0.016$ & $0.062$ & $0.070$ & $0.112$ & $0.709$ && $0.028$ & $0.099$ & $0.174$ & $0.224$ & $0.316$ & $0.947$ \\
R6L                            & $0.018$ & $0.017$ & $0.071$ & $0.073$ & $0.115$ & $0.714$ && $0.035$ & $0.088$ & $0.172$ & $0.215$ & $0.296$ & $0.941$ \\
R77                            & $0.016$ & $0.017$ & $0.066$ & $0.072$ & $0.109$ & $0.711$ && $0.032$ & $0.085$ & $0.169$ & $0.213$ & $0.292$ & $0.944$ \\
R78                            & $0.017$ & $0.019$ & $0.068$ & $0.074$ & $0.111$ & $0.712$ && $0.034$ & $0.090$ & $0.170$ & $0.213$ & $0.292$ & $0.946$ \\
R7L                            & $0.017$ & $0.016$ & $0.069$ & $0.071$ & $0.110$ & $0.712$ && $0.032$ & $0.084$ & $0.168$ & $0.209$ & $0.290$ & $0.943$ \\\hline
P-Test                         & $0.027$ & $0.031$ & $0.035$ & $0.042$ & $0.082$ & $0.727$ && $0.037$ & $0.064$ & $0.148$ & $0.193$ & $0.302$ & $0.958$ \\
D-Test                         & $0.025$ & $0.035$ & $0.054$ & $0.072$ & $0.134$ & $0.773$ && $0.028$ & $0.098$ & $0.182$ & $0.257$ & $0.345$ & $0.962$ \\
KS-Test                        & $0.013$ & $0.017$ & $0.030$ & $0.039$ & $0.065$ & $0.624$ && $0.021$ & $0.055$ & $0.120$ & $0.173$ & $0.234$ & $0.974$ \\
KL-Test                        & $0.014$ & $0.018$ & $0.028$ & $0.050$ & $0.085$ & $0.741$ && $0.023$ & $0.065$ & $0.147$ & $0.221$ & $0.303$ & $0.963$ \\
\hline\hline
\end{tabular}
\end{center}
\end{scriptsize}
\end{table}
}

Tables \ref{Tab: comparison2, n100} and \ref{Tab: comparison2, n400} record the rejection rates for $n\in\{100,400\}$ based on the design \eqref{Eqn: cdf compare2}. As expected, the $D$ and the $P$ tests, which are designed to be powerful against convex alternatives, now perform poorly and in most cases have power less than $5\%$---this is the case even with large sample sizes. The KS and the KL tests appear more powerful than the $D$ and the $P$ tests, but their performance is strictly dominated by our tests across all sample sizes and pairs $(t_n,\kappa_n)$, with substantial power discrepancies in many cases. We remind the reader that in practice, one can rarely rule out a priori distributions such as those in \eqref{Eqn: cdf compare2}. Thus, the simulations suggest that our test is more robust in the sense that they are powerful against a larger class of alternatives. Those results for $n\in\{200,300,600,800,1000\}$ share similar patterns---see Appendix \ref{App: Simulation}.

Admittedly, the simulation results reinforce the importance of the choice of $\kappa_n$ for our test---the choice of $t_n$ does not appear as important. Overall, they provide comforting numerical evidence that our test is a useful addition to the literature. The associated computation cost is reasonable: with an Intel Xeon CPU E5-1650 v4 3.60GHz, a single replication based on 1000 samples from the design in \eqref{Eqn: cdf compare1} and 500 bootstrap repetitions is completed in 15 seconds.

{\renewcommand{\arraystretch}{1.5}
\setlength{\tabcolsep}{5pt}
\begin{table}[!h]
\caption{Comparisons with Existing Tests Based on \eqref{Eqn: cdf compare2}: $n=100$}\label{Tab: comparison2, n100}
\begin{footnotesize}
\begin{center}
\begin{tabular}{cccccccccc}
\hline\hline
  & \multicolumn{9}{c}{$\lambda$}  \\
\cmidrule{2-10}
 & $0.1$ & $0.2$ & $0.3$ & $0.4$ & $0.5$ & $0.6$ & $0.7$ & $0.8$ & $0.9$ \\
\hline
R37    & $0.127$ & $0.221$ & $0.339$ & $0.482$ & $0.652$ & $0.664$ & $0.607$ & $0.386$ & $0.056$  \\
R38    & $0.125$ & $0.223$ & $0.336$ & $0.483$ & $0.656$ & $0.669$ & $0.607$ & $0.387$ & $0.053$     \\
R3L    & $0.127$ & $0.222$ & $0.338$ & $0.483$ & $0.654$ & $0.669$ & $0.608$ & $0.386$ & $0.055$  \\
R47    & $0.102$ & $0.184$ & $0.273$ & $0.406$ & $0.566$ & $0.586$ & $0.503$ & $0.296$ & $0.030$     \\
R48    & $0.103$ & $0.182$ & $0.273$ & $0.407$ & $0.569$ & $0.584$ & $0.499$ & $0.299$ & $0.029$     \\
R4L    & $0.102$ & $0.182$ & $0.275$ & $0.406$ & $0.563$ & $0.582$ & $0.504$ & $0.297$ & $0.031$  \\
R57    & $0.097$ & $0.164$ & $0.255$ & $0.374$ & $0.531$ & $0.537$ & $0.447$ & $0.248$ & $0.024$     \\
R58    & $0.095$ & $0.160$ & $0.252$ & $0.377$ & $0.529$ & $0.537$ & $0.448$ & $0.250$ & $0.024$     \\
R5L    & $0.096$ & $0.162$ & $0.256$ & $0.374$ & $0.527$ & $0.539$ & $0.449$ & $0.247$ & $0.024$  \\
R67    & $0.086$ & $0.153$ & $0.238$ & $0.358$ & $0.494$ & $0.508$ & $0.411$ & $0.213$ & $0.020$     \\
R68    & $0.084$ & $0.152$ & $0.239$ & $0.358$ & $0.497$ & $0.506$ & $0.412$ & $0.212$ & $0.020$     \\
R6L    & $0.086$ & $0.153$ & $0.237$ & $0.357$ & $0.499$ & $0.512$ & $0.412$ & $0.214$ & $0.020$     \\
R77    & $0.087$ & $0.147$ & $0.234$ & $0.346$ & $0.479$ & $0.484$ & $0.389$ & $0.193$ & $0.019$     \\
R78    & $0.085$ & $0.148$ & $0.233$ & $0.344$ & $0.475$ & $0.486$ & $0.389$ & $0.196$ & $0.019$     \\
R7L    & $0.085$ & $0.146$ & $0.234$ & $0.342$ & $0.480$ & $0.483$ & $0.394$ & $0.196$ & $0.019$     \\ \hline
P-Test & $0.015$ & $0.018$ & $0.008$ & $0.009$ & $0.010$ & $0.004$ & $0.003$ & $0.006$ & $0.003$  \\
D-Test & $0.030$ & $0.048$ & $0.039$ & $0.042$ & $0.032$ & $0.017$ & $0.006$ & $0.004$ & $0.000$  \\
KS-Test& $0.017$ & $0.055$ & $0.090$ & $0.121$ & $0.154$ & $0.126$ & $0.086$ & $0.021$ & $0.000$ \\
KL-Test& $0.025$ & $0.050$ & $0.043$ & $0.039$ & $0.015$ & $0.003$ & $0.000$ & $0.000$ & $0.000$ \\
\hline\hline
\end{tabular}
\end{center}
\end{footnotesize}
\end{table}
}

{\renewcommand{\arraystretch}{1.5}
\setlength{\tabcolsep}{5pt}
\begin{table}[!h]
\caption{Comparisons with Existing Tests Based on \eqref{Eqn: cdf compare2}: $n=400$}\label{Tab: comparison2, n400}
\begin{footnotesize}
\begin{center}
\begin{tabular}{cccccccccc}
\hline\hline
  & \multicolumn{9}{c}{$\lambda$}  \\
\cmidrule{2-10}
 & $0.1$ & $0.2$ & $0.3$ & $0.4$ & $0.5$ & $0.6$ & $0.7$ & $0.8$ & $0.9$ \\
\hline
R37    & $0.207$ & $0.595$ & $0.938$ & $0.995$ & $1.000$ & $1.000$ & $1.000$ & $1.000$ & $0.718$  \\
R38    & $0.207$ & $0.598$ & $0.935$ & $0.995$ & $1.000$ & $1.000$ & $1.000$ & $1.000$ & $0.721$  \\
R3L    & $0.207$ & $0.596$ & $0.937$ & $0.995$ & $1.000$ & $1.000$ & $1.000$ & $1.000$ & $0.723$  \\
R47    & $0.173$ & $0.516$ & $0.908$ & $0.991$ & $1.000$ & $1.000$ & $1.000$ & $0.995$ & $0.563$     \\
R48    & $0.176$ & $0.517$ & $0.908$ & $0.991$ & $1.000$ & $1.000$ & $1.000$ & $0.995$ & $0.564$     \\
R4L    & $0.175$ & $0.514$ & $0.908$ & $0.991$ & $1.000$ & $1.000$ & $1.000$ & $0.995$ & $0.561$     \\
R57    & $0.157$ & $0.485$ & $0.886$ & $0.987$ & $0.998$ & $1.000$ & $1.000$ & $0.985$ & $0.463$     \\
R58    & $0.155$ & $0.487$ & $0.886$ & $0.989$ & $0.998$ & $0.999$ & $1.000$ & $0.986$ & $0.470$     \\
R5L    & $0.158$ & $0.485$ & $0.887$ & $0.988$ & $0.998$ & $1.000$ & $1.000$ & $0.985$ & $0.467$     \\
R67    & $0.143$ & $0.469$ & $0.876$ & $0.985$ & $0.998$ & $0.999$ & $1.000$ & $0.975$ & $0.394$     \\
R68    & $0.143$ & $0.471$ & $0.875$ & $0.985$ & $0.998$ & $0.999$ & $1.000$ & $0.976$ & $0.393$     \\
R6L    & $0.144$ & $0.468$ & $0.875$ & $0.986$ & $0.998$ & $0.999$ & $1.000$ & $0.977$ & $0.391$     \\
R77    & $0.137$ & $0.462$ & $0.864$ & $0.980$ & $0.997$ & $0.999$ & $1.000$ & $0.966$ & $0.346$     \\
R78    & $0.134$ & $0.461$ & $0.864$ & $0.979$ & $0.997$ & $0.999$ & $1.000$ & $0.967$ & $0.350$     \\
R7L    & $0.136$ & $0.456$ & $0.866$ & $0.981$ & $0.997$ & $0.999$ & $1.000$ & $0.967$ & $0.348$     \\ \hline
P-Test & $0.020$ & $0.014$ & $0.010$ & $0.015$ & $0.010$ & $0.007$ & $0.007$ & $0.005$ & $0.003$  \\
D-Test & $0.050$ & $0.054$ & $0.042$ & $0.038$ & $0.032$ & $0.012$ & $0.011$ & $0.004$ & $0.000$  \\
KS-Test& $0.076$ & $0.295$ & $0.689$ & $0.907$ & $0.966$ & $0.977$ & $0.962$ & $0.750$ & $0.056$ \\
KL-Test& $0.107$ & $0.232$ & $0.356$ & $0.362$ & $0.298$ & $0.083$ & $0.003$ & $0.000$ & $0.000$ \\
\hline\hline
\end{tabular}
\end{center}
\end{footnotesize}
\end{table}
}

\section{Conclusion}\label{Sec: conclusion}

This paper studies estimation of and inference on a cumulative distribution function with concavity constraint. The estimation results generalize the asymptotic order results in \citet{Kiefer_Wolfowitz1976minimax} to settings with unbounded support and contiguous distributions. These results are not only of interest in their own right, but also useful for conducting inference on the concavity. In particular, in conjunction with the rescaled bootstrap of \citet{Dumbgen1993} and the recent work of \citet{FangSantos2018HDD}, they allow us to build up a test that controls size, pointwise and locally, even when the distribution function is strictly concave (in which case the test statistic is asymptotically degenerate). Through simulation studies, we find that our test is powerful against a larger class of alternatives, compared to some existing tests such as those in \citet{Woodroofe_Sun1999} that are designed to work against a specific class of alternatives, or those in \citet{Durot2003testing} and \citet{Kulikov_Lopuhaa2004testing} that rely on critical values from the least favorable asymptotic distributions. In Appendix \ref{App: extension}, we show how the testing results may be extended to a general setup that includes regression and hazard rate problems as special cases.

In our limited simulation studies, we also find that the choice of $\kappa_n$ for the KW selection does not appear as important as the choice of $t_n$ for implementing the rescaled bootstrap. A formal investigation of both choices is, however, beyond the scope of this paper, and therefore left for future study.


\clearpage

\newpage

\appendix

\section{Proofs of Main Results}\label{Sec: proofs}
\renewcommand{\theequation}{A.\arabic{equation}}
\setcounter{equation}{0}

\begin{figure}
\centering
\begin{tikzpicture} [domain=0:3.2]
\begin{axis}[width=0.8\textwidth,height=0.6\textwidth,
             axis y line=middle, axis x line=bottom,
             xmin=-0.2, xmax=3.2,
             ymin=0, ymax=1.1,
             xlabel=$x$,
             every axis y label/.style={at=(current axis.above origin),anchor=south,font=\footnotesize},
             every axis x label/.style={at=(current axis.right of origin),anchor=north,font=\footnotesize},
             ticklabel style = {font=\tiny},
             xtick={0,0.10536051565, 0.22314355131, 2.30258509299,2.99573227355},
             xticklabels={$a_0$, $a_1$, $a_2$, $a_{k_n-1}$, $a_{k_n}$},
             ytick={0.1,0.2, 0.9,0.95,1},
             yticklabels={$1/k_n$,$2/k_n$,  $1-1/k_n$,$1-1/(2k_n)$,1},
            ]
\addplot [NavyBlue,line width=0.2pt] {{1-exp(-\x)}}  node [pos=0.4, below right] {$F(x)$};
\draw[very thin,dashed] (0,0.1) -- (0.10536051565,0.1) -- (0.10536051565,0);
\draw[very thin,dashed] (0,0.2) -- (0.22314355131,0.2) -- (0.22314355131,0);
\draw[very thin,dashed] (0,0.9) -- (2.30258509299,0.9) -- (2.30258509299,0);
\draw[very thin,dashed] (0,0.95) -- (2.99573227355,0.95) -- (2.99573227355,0);
\draw[very thin,dashed] (0,1) -- (5,1);
\end{axis}
\end{tikzpicture}
\caption{Construction of the grid points $a_0, a_1, a_2, \ldots, a_{k_n}$}
\label{Fig: grid}
\end{figure}
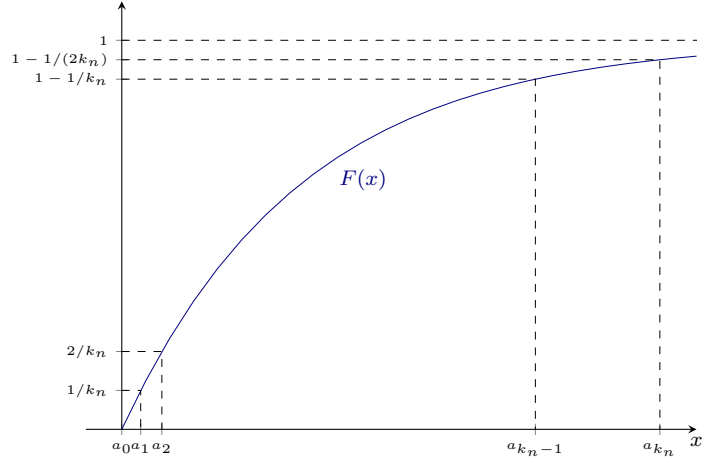

\noindent{\sc Proof of Theorem \ref{Thm: Equivalence}:} Our proof is a mixture of the original one in \citet{Kiefer_Wolfowitz1976minimax} and the ``modernized'' one given in \citet{Balabdaoui_Wellner2007Wolfowitz} but takes into account that the support of the density function $f$ is potentially unbounded. Following \citet{Kiefer_Wolfowitz1976minimax}, we first define interpolating processes for $F$ and $\mathbb F_n$. Let $k_n\uparrow\infty$ be a sequence of positive integers (to be chosen). Define $a_j\equiv a_j^{(k_n)}\equiv F^{-1}(j/k_n)$ for $j=1,\ldots,k_n-1$, and $a_{k_n}\equiv a_{k_n}^{(k_n)}\equiv F^{-1}(1-\frac{1}{(2k_n)})$. Moreover, set $a_0=0$ and $a_{k_n+1} \equiv a_{k_n+1}^{(k_n)} =\infty$. Let $L^{(k_n)}$ be the piecewise linear on $[a_{j-1},a_{j}]$ for $j=1,\ldots,k_n$ satisfying
\[
L^{(k_n)}(a_j^{(k_n)})=F(a_j^{(k_n)})~,\,j=0,\ldots,k_n~,
\]
and $L^{(k_n)}(x)=F(x)$ for $x\in[a_{k_n},a_{k_n+1}]$. Thus, in Figure \ref{Fig: grid}, $L^{k_n}$ is the function that connects $\{a_j,j/k_n\}_{j=0}^{k_n}$ on the interval $[0,a_{k_n}]$ in a piecewise linear way and is identical to $F(x)$ on the tail $[a_{k_n},a_{k_n+1}]$. Clearly, $L^{(k_n)}$ inherits concavity from $F$. Using the notation in \citep[p.31]{deBoor2001}, we may write $L^{(k_n)}= I_2 F$ even though $L^{(k_n)}$ is nonlinear on $[a_{k_n},a_{k_n+1}]$.

Next, we define $L_n^{(k_n)}$ by: for $x\in[a_j,a_{j+1}]$ and $j=0,\ldots,k_n$,
\[
L_n^{(k_n)}(x)\equiv \mathbb F_n(a_j)+\frac{\mathbb F_n(a_{j+1})-\mathbb F_n(a_j)}{F(a_{j+1})-F(a_j)}\{L^{(k_n)}(x)-F(a_j)\}~.
\]
Similarly, write $L_n^{(k_n)}= I_2\mathbb F_n$ and note, for $j=0,\ldots,k_n+1$,
\[
L_n^{(k_n)}(a_j^{(k_n)})=\mathbb F_n(a_j^{(k_n)}) ~.
\]
Thus, $L_n^{k_n}$ and $\mathbb F_n$ intersect at the grid points $\{a_j\}$, just as $L^{k_n}$ and $F$. Moreover, since $L_n^{(k_n)}$ is an affine transformation of $L^{(k_n)}$ on each of the intervals $[a_j,a_{j+1}]$, the former inherits the piecewise linearity from the latter on $[a_0,a_{k_n}]$. Heuristically, one may thus think of $L_n^{(k_n)}$ as a finite sample analog of $L^{(k_n)}$. We next show that $L_n^{(k_n)}$ also inherits concavity from $L^{(k_n)}$ for a suitable $k_n$.

\noindent\underline{\sc Step 1:} For $A_n\equiv \{L_n^{(k_n)} \text{ is concave on }\mathbf R_+\}$ and a suitably chosen sequence $\{k_n\}$, show that
\begin{align}\label{Eqn: Ln concave, aux}
\lim_{n\to\infty}P(A_n)=1~.
\end{align}
In fact, we shall show that $A_n$ occurs for all $n$ large with probability one. We follow the arguments given by \citet[p.81-2]{Kiefer_Wolfowitz1976minimax}. For $j=1,\ldots,k_n$, define the increments:
\[
T_{n,j}\equiv\mathbb F_n(a_j)-\mathbb F_n(a_{j-1})~,\,\Delta a_j\equiv a_j-a_{j-1}~.
\]
Note that the definition of $\Delta a_j$ is in line with \citet{Balabdaoui_Wellner2007Wolfowitz} but differs slightly from \citet{Kiefer_Wolfowitz1976minimax}. Since $L_n^{(k_n)}$ is piecewise linear on each $[a_j,a_{j-1}]$ for all $j=1,\ldots,k_n$, the slope of $L_n^{(k_n)}$ on each such $[a_j,a_{j-1}]$ is precisely $T_{n,j}/\Delta a_j$. In turn, it follows that establishing ``asymptotic concavity'' amounts to showing a string of inequalities comparing consecutive slopes---we remind the reader that $L_n^{(k_n)}$ is concave on the interval $[a_{k_n},a_{k_n+1}]$. To formalize the idea, write, for $j=1,\ldots,k_n-1$,
\[
B_{n,j}\equiv \left\{\frac{T_{n,j+1}}{\Delta a_{j+1}}\le \frac{T_{n,j}}{\Delta a_{j}}\right\}~,\,B_{n,k_n}\equiv \left\{\frac{T_{n,k_n+1}}{F(a_{k_n+1})-F(a_{k_n})}f(a_{k_n})\le \frac{T_{n,a_{k_n}}}{\Delta a_{k_n}}\right\}~,
\]
where $\frac{T_{n,k_n+1}}{F(a_{k_n+1})-F(a_{k_n})}f(a_{k_n})$ is the slope of $L_n^{(k_n)}$ at $a_{k_n}$. For all $j=1,\ldots,k_n$, if $B_{n,j}$ holds, then $L_n^{(k_n)}$ stays concave as it moves from $[a_{j-1},a_j]$ into $[a_j,a_{j+1}]$. Therefore, we obtain the following representation of the event $A_n$:
\begin{align}\label{Eqn: Ln concave, aux0}
A_n=\bigcap_{j=1}^{k_n} B_{n,j}~.
\end{align}
We next consider the sets $\{B_{n,j}\}$ one by one. The goal is to establish analytically tractable sufficient conditions for each $B_{n,j}$ to occur. Then we shall bound $P(A_n)$ from below by computing $\sum_{j=1}^{k_n} P(B_{n,j}^c)$ (in view of \eqref{Eqn: Ln concave, aux0}), which in turn may be controlled by utilizing those sufficient conditions.

\noindent\underline{\sc For $j=1,\ldots,k_n-2$:} It is simple to verify that $B_{n,j}$ occurs if
\begin{align}\label{Eqn: Ln concave, aux1}
|T_{n,i}-\frac{1}{k_n}|\le\frac{\delta_n}{k_n}~,\,i=j,j+1~,\text{ and }\frac{\Delta a_{j+1}}{\Delta a_j}\ge 1+3\delta_n~,
\end{align}
provided $0\le\delta_n\le 1/3$. Indeed, if \eqref{Eqn: Ln concave, aux1} holds, then
\begin{align}\label{Eqn: Ln concave, aux1a}
\frac{T_{n,j+1}}{T_{n,j}} \le \frac{1/k_n  +  \delta_n/k_n}{1/k_n - \delta_n/k_n} \le 1+ 3\delta_n\le \frac{\Delta a_{j+1}}{\Delta a_j}~,
\end{align}
so that $B_{n,j}$ occurs, where the second inequality holds whenever $0\le\delta_n\le 1/3$. Next, we show that $\frac{\Delta a_{j+1}}{\Delta a_j}\ge 1+3\delta_n$ holds for all $n$ large under an additional rate restriction on $\delta_n$ (to be specified)---we shall control the probability of $|T_{n,i}-\frac{1}{k_n}|\le\frac{\delta_n}{k_n}$ shortly. By Assumption \ref{Ass: equivalence}(i), we may employ Taylor's theorem to conclude that, for some $\xi_{j+1}\in[a_j,a_{j+1}]$,
\begin{align}\label{Eqn: Ln concave, aux2}
\Delta a_{j+1}=F^{-1}(\frac{j+1}{k_n})-F^{-1}(\frac{j}{k_n})=\frac{1}{f(a_j)}\frac{1}{k_n}+\frac{1}{2k_n^2}\frac{-f'(\xi_{j+1})}{f^3(\xi_{j+1})}~.
\end{align}
Moreover, since $F^{-1}$ is strictly convex (because $\frac{d}{dx}F^{-1}(x)=1/f(F^{-1}(x))$ is strictly increasing) and $\Delta a_j\equiv F^{-1}(\frac{j}{k_n})-F^{-1}(\frac{j-1}{k_n})$ by definition, the ratio $\Delta a_j/k_n^{-1}$ is no larger than the slope of $F^{-1}$ at the right end point of the interval $[(j-1)/k_n,j/k_n]$ (i.e., $j/k_n$). Mathematically, this means that
\begin{align}\label{Eqn: Ln concave, aux3}
\Delta a_j\le \frac{1}{f(a_j)}\frac{1}{k_n}~.
\end{align}
Combining previous results \eqref{Eqn: Ln concave, aux2} and \eqref{Eqn: Ln concave, aux3} yields
\begin{multline}\label{Eqn: Ln concave, aux3a0}
\frac{\Delta a_{j+1}}{\Delta a_j}\ge 1+\frac{1}{2k_n}\frac{-f'(\xi_{j+1})}{f^3(\xi_{j+1})} f(a_j)\\
\ge 1+\frac{1}{2k_n}\frac{-f'(\xi_{j+1})}{f^2(\xi_{j+1})}\ge 1+\frac{1}{2k_n}\bar\beta(\frac{1}{2k_n})~,
\end{multline}
where the second inequality holds since $f(\xi_{j+1})\le f(a_j)$ because $f$ is decreasing and $\xi_{j+1}\ge a_j$. By setting $\delta_n\le \frac{1}{6k_n}\bar\beta(\frac{1}{k_n})$, the second part of display \eqref{Eqn: Ln concave, aux1} thus holds. Note that $\bar\beta(\epsilon)$ is clearly nondecreasing in $\epsilon$ (because the domain of the infimum in \eqref{Eqn: beta} shrinks as $\epsilon$ increases), so $\delta_n\le 1/3$ for all large $n$. Moreover, by Assumption \ref{Ass: equivalence}(ii), $\bar\beta(\epsilon)>0$ for all small $\epsilon>0$ and so $\delta_n>0$ for large $n$.

\noindent\underline{\sc For $j=k_n-1$:} By identical arguments in \eqref{Eqn: Ln concave, aux1a}, we see $B_{n,k_n-1}$ occurs if
\begin{multline}\label{Eqn: Ln concave, aux4}
|T_{n,k_n-1}-\frac{1}{k_n}|\le\frac{\delta_n}{k_n}~,\,|T_{n,k_n}-\frac{1}{2k_n}|\le \frac{\delta_n}{2k_n}~,\\
\text{ and }2\frac{\Delta a_{k_n}}{\Delta a_{k_n-1}}\ge 1+3\delta_n~,
\end{multline}
provided $0\le\delta_n\le 1/3$. Indeed, the first two inequalities in \eqref{Eqn: Ln concave, aux4} imply
\begin{align*}
\frac{T_{n,k_n}}{T_{n,k_n-1}}\le \frac{1/(2k_n)+\delta_n/(2k_n)}{1/k_n-\delta_n/k_n}=\frac{1}{2}\cdot\frac{1+\delta_n}{1-\delta_n}\le \frac{1+3\delta_n}{2}\le \frac{\Delta a_{k_n}}{\Delta a_{k_n-1}}~,
\end{align*}
as claimed, where the second inequality holds provided $0\le \delta_n\le 1/3$. By arguments analogous to those leading to \eqref{Eqn: Ln concave, aux3a0}, the last displayed part of \eqref{Eqn: Ln concave, aux4} holds if $\delta_n\le \frac{1}{12k_n}\bar\beta(\frac{1}{2k_n})\le 1/3$. Indeed, by Taylor's theorem (as in \eqref{Eqn: Ln concave, aux2}) we have, for some $\xi_{k_n}\in [a_{k_n-1},a_{k_n}]$,
\begin{multline}\label{Eqn: Ln concave, aux4a}
\Delta a_{k_n}=F^{-1}(1-\frac{1}{2k_n})-F^{-1}(1-\frac{1}{k_n})\\
=\frac{1}{f(a_{k_n-1})}\frac{1}{2k_n}+\frac{1}{8k_n^2}\frac{-f'(\xi_{k_n})}{f^3(\xi_{k_n})}~,
\end{multline}
Moreover, \eqref{Eqn: Ln concave, aux3} holds for $j=k_n-1$ as well, i.e.,
\begin{align}\label{Eqn: Ln concave, aux4b}
\Delta a_{k_n-1}\le \frac{1}{f(a_{k_n-1})}\frac{1}{k_n}~.
\end{align}
A combination of \eqref{Eqn: Ln concave, aux4a} and \eqref{Eqn: Ln concave, aux4b}, together with the definition of $\bar\beta$, $f$ being decreasing, then implies that
\begin{multline}
\frac{\Delta a_{k_n}}{\Delta a_{k_n-1}}\ge \frac{1}{2}+\frac{1}{8k_n}\frac{-f'(\xi_{k_n})}{f^3(\xi_{k_n})}f(a_{k_n-1}) \\
\ge \frac{1}{2}+\frac{1}{8k_n}\bar\beta(\frac{1}{2k_n})\ge \frac{1}{2}(1+3\delta_n)~,
\end{multline}
where the last inequality holds if $\delta_n\le \frac{1}{12k_n}\bar\beta(\frac{1}{2k_n})$. We shall control the probabilities of the first two events in \eqref{Eqn: Ln concave, aux4} shortly.

\noindent\underline{\sc For $j=k_n$:} This case is easier because $B_{n,k_n}$ occurs whenever
\begin{align}\label{Eqn: Ln concave, aux5}
|T_{n,k_n}-\frac{1}{2k_n}|\le \frac{\delta_n}{2k_n} ~,\,T_{n,k_n+1} \le \frac{1-\delta_n}{2k_n}~.
\end{align}
To see this, note that, as in result \eqref{Eqn: Ln concave, aux3}, convexity of $F^{-1}$ and $\Delta a_{k_n}\equiv F^{-1}(1-1/(2k_n))-F^{-1}(1-1/k_n)$ together imply that the ratio $\Delta a_j/(1/(2k_n))$ is no larger than the slope of $F^{-1}$ at $1-1/(2k_n)$, i.e.,
\begin{align*}
\Delta a_{k_n}\le \frac{1}{f(a_{k_n})}\frac{1}{2k_n}~.
\end{align*}
Hence, by construction of the grid points, we have
\begin{align}\label{Eqn: Ln concave, aux5a}
\frac{F(a_{k_n+1})-F(a_{k_n})}{\Delta a_{k_n}}\frac{1}{f(a_{k_n})}= \frac{1/(2k_n)}{\Delta a_{k_n} f(a_{k_n})}\ge 1~.
\end{align}
It follows from \eqref{Eqn: Ln concave, aux5a} that, whenever \eqref{Eqn: Ln concave, aux5} occurs,
\begin{align*}
\frac{T_{n,k_n}}{\Delta a_{k_n}}\frac{F(a_{k_n+1})-F(a_{k_n})}{f(a_{k_n})}\ge \frac{1-\delta_n}{2k_n} \ge T_{n,k_n+1}~.
\end{align*}

At this point, we note that $\delta_n$ should be selected to satisfy (i) $0<\delta_n\le 1/3$, (ii) $\delta_n\le \bar\beta(1/k_n)/(6k_n)$, and (iii) $\delta_n\le\bar\beta(1/(2k_n))/(12k_n)$. Since $\bar\beta(\epsilon)$ decreases as $\epsilon\downarrow 0$, we may thus set $\delta_n=\bar\beta(1/(2k_n))/(12k_n)$ in what follows, so that (i)(ii)(iii) hold simultaneously for all $n$ large. Note that $\delta_n\downarrow 0$ since $k_n\uparrow\infty$ as required (and to be constructed). We may now obtain from \eqref{Eqn: Ln concave, aux0}, and results \eqref{Eqn: Ln concave, aux1}, \eqref{Eqn: Ln concave, aux4} and \eqref{Eqn: Ln concave, aux5} that,  for all $n$ large enough so that $\delta_n\le1/3$,
\begin{multline}\label{Eqn: Ln concave, aux5b}
P(A_n^c) \le\sum_{j=1}^{k_n} P(B_{n,j}^c)
\le\sum_{j=1}^{k_n-1} 2P(|T_{n,j}-\frac{1}{k_n}|>\frac{\delta_n}{k_n}) \\
+2P(|T_{n,k_n}-\frac{1}{2k_n}|>\frac{\delta_n}{2k_n})+ P(T_{n,k_n+1} \le \frac{1-\delta_n}{2k_n})~.
\end{multline}
Following \citet{Balabdaoui_Wellner2007Wolfowitz}, we may apply Theorem 1.1.2 in \citet{Shorack_Wellner1986empirical} and Lemma 5.2 in \citet{Balabdaoui_Wellner2007Wolfowitz} to conclude that, for $p_n\equiv k_n^{-1}$ and each $j=1,\ldots,k_n-1$,
\begin{align}\label{Eqn: Ln concave, aux5c}
P(|T_{n,j}-\frac{1}{k_n}|>\frac{\delta_n}{k_n})&\le 2\exp\{-\frac{1}{2}np_n\delta_n^2(1+o(1))~,
\end{align}
where the $o(1)$ term depends only on $\{\delta_n\}$ (and nothing else, including $j$). The same arguments also imply
\begin{align}\label{Eqn: Ln concave, aux5d}
P(|T_{n,k_n}-\frac{1}{2k_n})|>\frac{\delta_n}{2k_n})\le 2\exp\{-\frac{1}{4}np_n\delta_n^2(1+o(1))~,
\end{align}
where the $o(1)$ term depends only on $\delta_n$. For the last term in \eqref{Eqn: Ln concave, aux5b}, we need somewhat different arguments, which is important in handling the unbounded support. Specifically, by Theorem 1.1.2 in \citet{Shorack_Wellner1986empirical},
\begin{align}\label{Eqn: Ln concave, aux5e}
P(T_{n,k_n+1} \le \frac{1-\delta_n}{2k_n})=P(\frac{1/(2k_n)}{\mathbb G_n(1/(2k_n))} \ge \frac{1}{1-\delta_n})~,
\end{align}
where $\mathbb G_n$ is the empirical cdf defined by $n$ i.i.d.\ uniform random variables. Inequality 10.3.2-(6) in \citet[p.415-6]{Shorack_Wellner1986empirical} then implies
\begin{multline}\label{Eqn: Ln concave, aux5f}
P(\frac{1/(2k_n)}{\mathbb G_n(1/(2k_n))} \ge \frac{1}{1-\delta_n})\le P(\sup_{t=1/(2k_n)}^1\frac{t}{\mathbb G_n(t)} \ge \frac{1}{1-\delta_n})\\
\le \exp\{-\frac{1}{2}np_n h(1-\delta_n)\}~,
\end{multline}
where $h(x)\equiv x(\log x-1)+1$ for $x>0$. By L'H\^{o}pital's rule, $h(1+x)=\frac{1}{2}x^2(1+o(1))$ as $x\to 0$. Since $\delta_n\downarrow 0$ by our choice, we may thus conclude from results \eqref{Eqn: Ln concave, aux5e} and \eqref{Eqn: Ln concave, aux5f} that
\begin{align}\label{Eqn: Ln concave, aux5g}
P(T_{n,k_n+1} \le \frac{1-\delta_n}{2k_n})\le \exp\{-\frac{1}{4}np_n \delta_n^2(1+o(1))\}~,
\end{align}
where the $o(1)$ term again depends only on $\delta_n$. Finally, we conclude from results \eqref{Eqn: Ln concave, aux5b}, \eqref{Eqn: Ln concave, aux5c}, \eqref{Eqn: Ln concave, aux5d} and \eqref{Eqn: Ln concave, aux5g}, together with the choice of $\delta_n$, that
\begin{multline}\label{Eqn: Ln concave, aux5j}
P(A_n^c)\le 5k_n \exp\{-\frac{1}{4}np_n\delta_n^2(1+o(1)) \\
\le  5k_n \exp\{-\frac{1}{720}nk_n^{-3}\bar\beta^2(\frac{1}{2k_n}))~,
\end{multline}
for all $n$ large so that the $1+o(1)$ term is larger than $4/5$.

Now, we may choose $k_n$ to be the integer part of $\nu$ that satisfies
\begin{align}\label{Eqn: Ln concave, aux6a}
\frac{\nu^3}{\bar\beta^2(1/(2\nu))}=\frac{3}{7}\frac{1}{720}\frac{n}{\log n}~,
\end{align}
where $\nu\ge 1$. To see that a solution to \eqref{Eqn: Ln concave, aux6a} exists, we note several facts. First, $\epsilon\mapsto\bar\beta(\epsilon)$ being nondecreasing implies that $\nu\mapsto\nu^3/\bar\beta^2(1/(2\nu))$ is strictly increasing and diverges to infinity as $\nu\uparrow\infty$. Second, $\epsilon\mapsto\bar\beta(\epsilon)$ in \eqref{Eqn: beta} is continuous and hence so is $\nu\mapsto\nu^3/\bar\beta^2(1/(2\nu))$. Indeed, Assumption \ref{Ass: equivalence}(i) and $f(x)>0$ for each $x\in\mathbf R_+$ (by Assumption \ref{Ass: setup}(ii)) implies that $x\mapsto -f'(x)/f^2(x)$ is continuous on $\mathbf R_+$. Moreover, the correspondence $\varphi:(0,1) \twoheadrightarrow\mathbf R_+$ defined by $\varphi(\epsilon)=[0,F^{-1}(1-\epsilon)]$ is compact-valued, and continuous (justified by first writing $\varphi$ as the composition of the singleton-valued continuous map $\epsilon\mapsto F^{-1}(1-\epsilon)$ and the continuous correspondence $x\twoheadrightarrow[0,x]$, and then invoking Theorem 17.23 in \citet{AliprantisandBorder2006}). Berge's maximum theorem---see, for example, Theorem 17.31 in \citet{AliprantisandBorder2006}---then implies that $\bar\beta(\epsilon)$ is continuous in $\epsilon$. Third, the right hand side of \eqref{Eqn: Ln concave, aux6a} diverges to infinity as $n\to\infty$. Thus, for each fixed $n$ that is sufficiently large, the right hand side of \eqref{Eqn: Ln concave, aux6a} falls within the range $\{\nu^3/\bar\beta^2(1/(2\nu)): \nu\ge 1\}$.

The above three facts together imply that there is a unique $\nu=\nu_n$ that satisfies equation \eqref{Eqn: Ln concave, aux6a} for all $n$ large and diverges to infinity as $n\to\infty$. To evaluate the order of $\nu_n$, we note that, by Assumption \ref{Ass: equivalence}(ii),
\begin{align}\label{Eqn: Ln concave, aux6b}
\nu_n^3= \frac{3}{7}\frac{1}{720}\frac{n}{\log n}\bar\beta^2(\frac{1}{2\nu_n})\asymp \frac{n}{\log n} \nu_n^{-2\tau}~,
\end{align}
implying that $\nu_n\asymp (n/\log n)^{1/(2\tau+3)}$. Since $k_n$ is the integer part of $\nu_n$, we have $k_n\asymp (n/\log n)^{1/(2\tau+3)}$, which diverges to infinity as $n\uparrow\infty$ since $2\tau+3>0$ by Assumption \ref{Ass: equivalence}(iii). Since the bound in \eqref{Eqn: Ln concave, aux5j} is weakly increasing in $k_n$, we may thus conclude by results \eqref{Eqn: Ln concave, aux5j} and \eqref{Eqn: Ln concave, aux6b} that, for all $n$ large,
\begin{multline}
P(A_n^c)\le 5\nu_n \exp\{-\frac{1}{720}n\nu_n^{-3}\bar\beta^2(\frac{1}{2\nu_n})) \\
= 5\nu_n \exp\{-\frac{7}{3} \log n\}\lesssim n^{1/3} n^{-7/3} = n^{-2}\to 0~.
\end{multline}
Thus, $P(A_n)\to 1$. Since $\sum_{i=1}^{n}P(A_n^c)=\sum_{i=1}^{n}n^{-2}<\infty$, the first Borel-Cantelli lemma implies that, with probability one, $A_n$ is concave for all $n$ large.

\noindent\underline{\sc Step 2:} Relate $\|\hat{\mathbb F}_n-\mathbb F_n\|_\infty$ to the interpolating process $L^{(k_n)}$ and $L_n^{(k_n)}$. For notational simplicity, define
\begin{align*}
B_n & \equiv \sup_{x\in[a_{k_n},a_{k_n+1}]}|\mathbb F_n(x)- L_n^{(k_n)}(x)|~,\\
D_n & \equiv \sup_{x\in[0,a_{k_n}]}|\mathbb F_n(x)-I_2\mathbb F_n(x)-\{F(x)-I_2F(x)\}|~,\\
E_n & \equiv \|F-L^{(k_n)}\|_\infty~.
\end{align*}
By the triangle inequality and Lemma 2.2 in \citet{Durot_Tocquet2003}, we have
\begin{multline}\label{Eqn: decomposition, BDEa}
\|\hat{\mathbb F}_n-\mathbb F_n\|_\infty  \le \|\hat{\mathbb F}_n-L_n^{(k_n)}-\{\mathbb F_n-L_n^{(k_n)}\}\|_\infty\\
 \le \|\hat{\mathbb F}_n-L_n^{(k_n)}\|_\infty+\|\mathbb F_n-L_n^{(k_n)}\|_\infty\le 2\|\mathbb F_n-L_n^{(k_n)}\|_\infty~.
\end{multline}
Invoking the triangle inequality once again we in turn have from \eqref{Eqn: decomposition, BDEa} that
\begin{align}\label{Eqn: decomposition, BDE}
\|\hat{\mathbb F}_n-\mathbb F_n\|_\infty & \le  2\|\mathbb F_n-L_n^{(k_n)}-\{F-L^{(k_n)}\}\|_\infty+2\|F-L^{(k_n)}\|_\infty\notag\\
&\le 2(B_n+D_n+E_n)~,
\end{align}
where the second inequality exploited $F(x)=L^{(k_n)}(x)$ for all $x\in[a_{k_n},a_{k_n+1}]$.

\noindent\underline{\sc Step 3:} Conclude by controlling the orders of $B_n$, $D_n$ and $E_n$ separately. To this end, we introduce the concept of modulus of continuity following \citet[p.25]{deBoor2001}. For a generic function $g:[a,b]\to\mathbf R$, let
\[
\omega(g;h)\equiv\sup\{|g(x)-g(y)|: x,y\in[a,b], |x-y|\le h\}~.
\]
The treatment of $D_n$ follows closely \citet{Balabdaoui_Wellner2007Wolfowitz}, which we include here for completeness. Note that
\begin{align}\label{Eqn: Dn}
D_n\le\omega(\mathbb F_n-F;|a|)\overset{d}{=}n^{-1/2}\omega(\mathbb U_n;p_n)~,
\end{align}
where the inequality is due to result (18) in \citet[p.36]{deBoor2001} with $|a|\equiv\max_{j=1}^k(a_j-a_{j-1})$, and $\overset{d}{=}$ denotes equality in distribution and is due to Theorem 1.1.2 in \citet{Shorack_Wellner1986empirical}, with $\mathbb U_n\equiv\sqrt{n}\{\mathbb G_n-I\}$ the empirical process of $n$ i.i.d.\ Uniform$(0,1)$ random variables as defined in \citet{Shorack_Wellner1986empirical}. By Theorem 14.2.1 in \citet{Shorack_Wellner1986empirical}, we have
\begin{align}\label{Eqn: Dn, aux1}
\lim_{n\to\infty}\frac{\omega(\mathbb U_n;p_n)}{\sqrt{2p_n\log(1/p_n)}}=1\text{ a.s.,}
\end{align}
provided $p_n\to 0$, $np_n\to \infty$, $\log(1/p_n)/\log(\log n)\to\infty$, and $\log(1/p_n)/(np_n)\to 0$. To see that these rate conditions on $p_n\equiv k_n^{-1}$ are indeed met, we first note that $p_n\asymp ((\log n)/n)^{1/(2\tau+3)}$. Thus, trivially, $p_n\to 0$ and
\begin{align}
np_n= n^{\frac{2\tau+2}{2\tau+3}}(\log n)^{\frac{1}{2\tau+3}}\to\infty~,
\end{align}
as $n\to\infty$, since $\frac{2\tau+2}{2\tau+3}>0$ by Assumption \ref{Ass: equivalence}(ii). The third rate condition is also trivial because $\log(1/p_n)$ is of order $\log(n)-\log(\log n)$, diverging to infinity faster than $\log(\log n)$. For the fourth condition, $\log(1/p_n)$ grows at a logarithmic rate while $np_n$ grows at a polynomial rate---note $\tau+1>0$ and so $\frac{2\tau+2}{2\tau+3}>0$. Now, by simple algebra, we may conclude from result \eqref{Eqn: Dn, aux1} that
\begin{align}\label{Eqn: Dn order}
D_n=O_{a.s.}(n^{-1/2}\sqrt{p_n\log(1/p_n)})
=O_{a.s.}(((\log n)/n)^{\frac{\tau+2}{2\tau+3}})~.
\end{align}

We now turn to $B_n$ which is important in taking care of the unbounded support. Heuristically, $B_n$ measures the interpolating error for $\mathbb F_n$ on the (unbounded) right tail $[a_{k_n},a_{k_n+1}]$ of the support.

Plugging the expression of $L_n^{(k_n)}$ into $B_n$ we have,
\begin{align}\label{Eqn: Bn}
B_n&= \sup_{x\in[a_{k_n},a_{k_n+1}]} |\mathbb F_n(a_{k_n})+\frac{\mathbb F_n(a_{k_n+1})-\mathbb F_n(a_{k_n})}{F(a_{k_n+1})-F(a_{k_n})}\{L^{(k_n)}(x)-F(a_{k_n})\}-\mathbb F_n(x)|\notag\\
&\le \sup_{x\in[a_{k_n},a_{k_n+1}]}|\mathbb F_n(x)-\mathbb F_n(a_{k_n})|+|\mathbb F_n(a_{k_n+1})-\mathbb F_n(a_{k_n})|\notag\\
&=2\{\mathbb F_n(a_{k_n+1})-\mathbb F_n(a_{k_n})\}~,
\end{align}
where the inequality follows by $L^{(k_n)}=F$ on $[a_{k_n},a_{k_n+1}]$ so that $|L^{(k_n)}(x)-F(a_{k_n})|/|F(a_{k_n+1})-F(a_{k_n})|\le 1$, and the second equality follows by monotonicity of $\mathbb F_n$. Recall that $T_{n,k_n+1}\equiv\mathbb F_n(a_{k_n+1})-\mathbb F_n(a_{k_n})$, which, based on our previous results, satisfies: for all $n$ large,
\begin{align}\label{Eqn: Bn, aux1}
P(T_{n,k_n+1} \le \frac{1-\delta_n}{2k_n}) \lesssim n^{1/3} n^{-7/3} = n^{-2}~.
\end{align}
By $\delta_n\downarrow 0$ and the first Borel-Cantelli lemma, we obtain from \eqref{Eqn: Bn, aux1} that
\begin{align}\label{Eqn: Bn, aux2}
T_{n,k_n+1}= O_{a.s.}(k_n^{-2}) = O_{a.s.}((\log n)/n)^{2/(2\tau+3)})~,
\end{align}
where the second equality follows by $k_n\asymp (n/\log n)^{1/(2\tau+3)}$. It follows from results \eqref{Eqn: Bn} and \eqref{Eqn: Bn, aux2} that
\begin{align}\label{Eqn: Bn order}
B_n = O_{a.s.}((\log n)/n)^{2/(2\tau+3)})~.
\end{align}

Finally, consider $E_n$ which controls the interpolation error for $F$ using $L^{(k_n)}$. Recall that $L^{(k_n)}$ is piecewise linear on $[0,a_{k_n}]$, and equal to $F$ on $[a_{k_n},a_{k_n+1}]$. It follows by result (5) in \citet[p.32]{deBoor2001}, Assumption \ref{Ass: equivalence}, and Theorem 1 in \citet{Burchard1974splines} (see also Theorem 1 in \citet{deBoor1973}) that
\begin{align}\label{Eqn: En order}
E_n=\|F-L^{(k_n)}\|_\infty &\le C k_n^{-2}\left[\int_0^{a_{k_n}}|f'(x)|^{1/2}\,\mathrm dx\right]^2\notag\\
&\le C k_n^{-2}\|f'\|_{1/2}=O((\log n)/n)^{2/(2\tau+3)})~,
\end{align}
where $C$ is a constant. If $\tau\ge 0$, then $\frac{2}{2\tau+3}\le\frac{\tau+2}{2\tau+3}$; if $-1<\tau<0$, then $\frac{2}{2\tau+3}>\frac{\tau+2}{2\tau+3}$. Therefore, combining results \eqref{Eqn: Dn order}, \eqref{Eqn: Bn order}, and \eqref{Eqn: En order}, together with {\sc Step 1}, we may establish the theorem.\qed

\noindent{\sc Proof of Theorem \ref{Thm: Equivalence, local}:} Let $\{k_n\}$ be the same sequence that is determined by \eqref{Eqn: Ln concave, aux6b}. Define the grid points $\{a_j\}_{j=0}^{k_n+1}$ and the interpolation processes $L^{(k_n)}$ and $L_n^{(k_n)}$ in the same fashion as before (based on $F$). The proof analogously consists of three steps as in the proof of Theorem \ref{Thm: Equivalence}. One of the key observations is that Theorem 1.1.2 in \citet{Shorack_Wellner1986empirical} holds for any i.i.d.\ $\{X_i\}_{i=1}^n$ whose common distribution function is allowed to depend on $n$.

\noindent\underline{\sc Step 1:} For $A_n\equiv \{L_n^{(k_n)} \text{ is concave on }\mathbf R_+\}$ (as before), show that \eqref{Eqn: Ln concave, aux} holds. This is accomplished by the same arguments as before. In particular, the arguments preceding \eqref{Eqn: Ln concave, aux5j} are purely algebraic, Theorem 1.1.2 in \citet{Shorack_Wellner1986empirical} (for result \eqref{Eqn: Ln concave, aux5j}) holds for any $n$, and Lemma 5.2 in \citet{Balabdaoui_Wellner2007Wolfowitz} and Inequality 10.3.2 in \citet{Shorack_Wellner1986empirical} are concerned with the uniform empirical process whose relations to $\mathbb F_n$ only depend on the i.i.d.\ assumption (for each $n$) as characterized by Theorem 1.1.2 in \citet{Shorack_Wellner1986empirical}. This last observation means that, again, the common cdf shared by $X_1,\ldots,X_n$ can depend on $n$, as in Assumption \ref{Ass: local}.

\noindent\underline{\sc Step 2:} For $B_n$, $D_n$ and $E_n$ defined as before, show that \eqref{Eqn: decomposition, BDE} holds. This is a consequence of the triangle inequality and the fact $F(x)=L^{(k_n)}(x)$ for all $x\in[a_{k_n},a_{k_n+1}]$. The arguments are thus the same as before.

\noindent\underline{\sc Step 3:} Control $B_n$, $D_n$ and $E_n$. The treatment of $B_n$ is the same as before, because Theorem 1.1.2 in \citet{Shorack_Wellner1986empirical} is still applicable which then allows one to invoke Inequality 10.3.2 in \citet{Shorack_Wellner1986empirical}. The treatment of $E_n$ is merely a consequence of the same approximation result, namely, Theorem 1 in \citet{Burchard1974splines}. It thus remains to consider $D_n$.

Let $F_n\equiv F(P_{1/\sqrt n})$ (the true underlying cdf). We then have:
\begin{align}\label{Eqn: equivalence local, aux0}
D_n&\equiv \sup_{x\in[0,a_{k_n}]}|\mathbb F_n(x)-I_2\mathbb F_n(x)-\{F(x)-I_2F(x)\}|\notag\\
&\le\omega(\mathbb F_n-F;|a|)\le \omega(\mathbb F_n-F_n;|a|)+\omega( F_n-F;|a|)\notag\\
&\overset{d}{=}n^{-1/2}\omega(\mathbb U_n;p_n)+\omega(F_n-F;|a|)~,
\end{align}
where the suprema in the $\omega$'s above are all taken over $[0,a_{k_n}]$, $|a|\equiv\max_{j=1}^k(a_j-a_{j-1})$, the second inequality follows by the subadditivity of $g\mapsto\omega(g,|a|)$ (as an implication of the supremum operator), and $\overset{d}{=}$ is due to Theorem 1.1.2 in \citet{Shorack_Wellner1986empirical}. Since $n^{-1/2}\omega(\mathbb U_n;p_n)=O_{a.s.}(((\log n)/n)^{\frac{\tau+2}{2\tau+3}})$ as before (see \eqref{Eqn: Dn order}), it remains to analyze $\omega( F_n-F;|a|)$.

Let $x_n,y_n\in[0,a_{k_n}]$ with $|x_n- y_n|\le a$ and $|F(x_n)-F(y_n)|\le p_n$. Define $f_n\equiv \mathrm dP_{1/\sqrt{n}}/\mathrm d\mu$ and $f\equiv dP/d\mu$ with $\mu$ the Lebesgue measure, which exist by the definition of $\mathbf P$. Since $\{P_{1/\sqrt{n}}\}\subset\mathbf P\subset\mathbf M$, we then obtain by the fundamental theorem of calculus and Jensen's inequality that:
\begin{align}\label{Eqn: equivalence local, aux1}
|\{&F_n(x_n)-F(x_n)\}-\{F_n(y_n)-F(y_n)\}|\le \int_{\mathbf R_+}1\{[x_n,y_n]\}|f_n-f|\,\mathrm d\mu\notag\\
&= \int_{\mathbf R_+}1\{[x_n,y_n]\}|f_n-f|1\{f=0\}\,\mathrm d\mu + \int_{\mathbf R_+}1\{[x_n,y_n]\}|f_n-f|1\{f\neq 0\}\,\mathrm d\mu\notag\\
&\le \int_{\mathbf R_+}f_n1\{f=0\}\,\mathrm d\mu + \int_{f\neq 0}1\{[x_n,y_n]\}|\frac{f_n}{f}-1|\,\mathrm dP~,
\end{align}
where the equality follows by the linearity of integrals, and the second inequality exploited $f\cdot 1\{f=0\}=0$ in the first integral. Here, $1\{f=0\}$ and $1\{f\neq 0\}$ are the indicator functions of the sets $\{f=0\}$ and $\{f\neq 0\}$ respectively. By Proposition A.5.3(E) in \citet[p.459]{BKRW993Efficient}, we have
\begin{align}\label{Eqn: equivalence local, aux2}
\int_{\mathbf R_+}f_n1\{f=0\}\,\mathrm d\mu=o(n^{-1})=O(((\log n)/n)^{\frac{\tau+2}{2\tau+3}})~.
\end{align}
For each $n\in\mathbf N$, let $r_n: \mathbf R\to\mathbf R$ be defined by (on $\{f\neq 0\}$):
\[
\frac{f_n}{f}=1+n^{-1/2}h+n^{-1/2}r_n~.
\]
By the triangle inequality and the definition of $r_n$, we then have
\begin{align}\label{Eqn: equivalence local, aux3}
\int_{f\neq 0}1\{[x_n,&y_n]\}|\frac{f_n}{f}-1|\,\mathrm dP \le \int_{f\neq 0}1\{[x_n,y_n]\} |n^{-1/2}r_n|\,\mathrm dP+\int_{\mathbf R_+} 1\{[x_n,y_n]\}| n^{-1/2} h|\,\mathrm dP\notag\\
&\le   \int_{f\neq 0}1\{[x_n,y_n]\} |n^{-1/2}r_n|\,\mathrm dP+ n^{-1/2}\{P([x_n,y_n])\}^{1/2} \{\int |h|^2\,\mathrm dP\}^{1/2}\notag\\
&\le  \int_{f\neq 0}1\{[x_n,y_n]\} |n^{-1/2}r_n|\,\mathrm dP+ n^{-1/2} p_n^{1/2}  \{\int h^2\,\mathrm dP\}^{1/2}\notag\\
&=  \int_{f\neq 0}1\{[x_n,y_n]\} |n^{-1/2}r_n|\,\mathrm dP+ O(((\log n)/n)^{\frac{\tau+2}{2\tau+3}})~,
\end{align}
where the second inequality follows by the Cauchy-Schwarz inequality, the third inequality by $P([x_n,y_n])=F(y_n)-F(x_n)\le p_n$, and the last step by $p_n\asymp ((\log n)/n)^{1/(2\tau+3)}$ and $Ph^2<\infty$. Next, note
\begin{multline}\label{Eqn: equivalence local, aux4}
\int_{f\neq 0}1\{[x_n,y_n]\} |n^{-1/2}r_n| \,  \mathrm dP\le \int_{f\neq 0}1\{[x_n,y_n]\} |n^{-1/2}r_n| 1\{|n^{-1/2}r_n|>1\}\,\mathrm dP \\
 + \int_{f\neq 0}1\{[x_n,y_n]\} |n^{-1/2}r_n| 1\{|n^{-1/2}r_n|\le 1\}\,\mathrm dP~.
\end{multline}
For the first term on the right side above, we have
\begin{align}\label{Eqn: equivalence local, aux5}
\int_{f\neq 0}1\{[x_n,y_n]\} |n^{-1/2}r_n| & 1\{|n^{-1/2}r_n|>1\}\,\mathrm dP\notag\\
&\le \int_{f\neq 0}|n^{-1/2}r_n| 1\{|n^{-1/2}r_n|>1\}\,\mathrm dP\notag\\
&=o(n^{-1})=O(((\log n)/n)^{\frac{\tau+2}{2\tau+3}})~,
\end{align}
where the first equality follows from \citet[p.27, lines 5-8]{Pfanzagl1985expansion}, the equivalence between pathwise/Hellinger differentiability and a notion of weak differentiability (which implies the condition (1.2.11') on p.25 in \citet{Pfanzagl1985expansion}). For the second term, by the Cauchy-Schwarz inequality, we have
\begin{align}\label{Eqn: equivalence local, aux6}
\int_{f\neq 0}1\{[x_n,y_n]\}& |n^{-1/2}r_n| 1\{|n^{-1/2}r_n|\le 1\}\,\mathrm dP\notag\\
&\le p_n^{1/2} \{\int_{f\neq 0} |n^{-1/2}r_n|^2 1\{|n^{-1/2}r_n|\le 1\}\,\mathrm dP\}^{1/2}\notag\\
&= p_n^{1/2}n^{-1/2} \{\int_{f\neq 0} |r_n|^2 1\{|n^{-1/2}r_n|\le 1\}\,\mathrm dP\}^{1/2}\notag\\
&=p_n^{1/2} n^{-1/2}o(1)=O(((\log n)/n)^{\frac{\tau+2}{2\tau+3}})~,
\end{align}
where the second equality again follows from \citet[p.27]{Pfanzagl1985expansion} but now using the condition (1.2.11'') on p.25 in \citet{Pfanzagl1985expansion}. Combining results \eqref{Eqn: equivalence local, aux1}, \eqref{Eqn: equivalence local, aux2}, \eqref{Eqn: equivalence local, aux3}, \eqref{Eqn: equivalence local, aux4}, \eqref{Eqn: equivalence local, aux5}, and \eqref{Eqn: equivalence local, aux6}, we finally obtain that
\begin{align}
|\{F_n(x_n)-F(x_n)\}-\{F_n(y_n)-F(y_n)\}|=O(((\log n)/n)^{\frac{\tau+2}{2\tau+3}})~,
\end{align}
and hence, in view of \eqref{Eqn: equivalence local, aux0},
\begin{align}
D_n=O_{a.s.}(((\log n)/n)^{\frac{\tau+2}{2\tau+3}})~.
\end{align}
This completes the proof of the theorem.\qed

\noindent{\sc Proof of Lemma \ref{Lem: weak limit}:} By Assumption \ref{Ass: setup}(i), we obtain by Example 19.6 in \citet{Vaart1998} that $\sqrt n\{\mathbb F_n-F\}\convl\mathbb G$ in $\ell^\infty(\mathbf R_+)$, where $\mathbb G\equiv\mathbb B\circ F$---see also \citet[p.266]{Vaart1998} for a brief discussion of the limit $\mathbb G$. Define a map $\psi: \ell^\infty(\mathbf R_+)\to\ell^\infty(\mathbf R_+)$ by $\psi(g)=\mathcal Mg-g$ for any $g\in\ell^\infty(\mathbf R_+)$. Then by definition $\phi(g)=\|\psi(g)\|_p$. Since $F$ is concave under $\mathrm H_0$, Proposition 2.1 in \citet{Beare_Fang2016Grenander} implies that $\psi$ is Hadamard directionally differentiable tangentially to $C_0(\mathbf R_+)$ with the derivative $\psi_F': C_0(\mathbf R_+)\to\mathbf R$ given by $\phi_F'(h)=\mathcal M_F'h-h $. for all $h\in C_0(\mathbf R_+)$. Moreover, $\|\cdot\|_p$ is Hadamard directionally differentiable at 0 since $(\|0+t_nh_n\|_p-\|0\|_p)/t_n=\|h_n\|_p\to \|h\|_p$ whenever $t_n\downarrow 0$ and $h_n\to h$ in $\ell^\infty(\mathbf R_+)$. By the chain rule as in Proposition 3.6 in \citet{Shapiro1990}, we may therefore conclude that $\phi: \ell^\infty(\mathbf R_+)\to\mathbf R$ is Hadamard directionally differentiable tangentially to $C_0(\mathbf R_+)$ with the derivative $\phi_F': C_0(\mathbf R_+)\to\mathbf R$ given by $\phi_F'(h)=\|\mathcal M_F'h-h\|_p$. This, in view of Theorem 2.1 in \citet{FangSantos2018HDD}, implies that: under $\mathrm H_0$ (and so $\phi(F)=0$),
\begin{align}
\sqrt n\phi(\mathbb F_n)=\sqrt n\{\phi(\mathbb F_n)-\phi(F)\}\convl\phi_F'(\mathbb G)=\|\mathcal{M}_F'(\mathbb G)-\mathbb G\|_p~.
\end{align}
This completes the proof of the lemma.\qed

\noindent{\sc Proof of Theorem \ref{Thm: test concavity, pointwise}:} Consider first the case when $F$ is strictly concave on $\mathbf R_+$. By Assumption \ref{Ass: kn}(ii), $n^{-2/3}(\log n)^{2/3}/(n^{-1/2}\kappa_n)\to 0$ as $n\to\infty$. In turn, we may thus have by Theorem \ref{Thm: Equivalence} that
\begin{align}
\limsup_{n\to\infty}P(\sqrt{n}\phi(\mathbb F_n)> \hat c_{1-\alpha})\le &\limsup_{n\to\infty} P(\sqrt{n}\phi(\mathbb F_n)> \kappa_n)\notag\\
&=\limsup_{n\to\infty}P(\frac{\|\hat{\mathbb F}_n-\mathbb F_n\|_\infty}{n^{-1/2}\kappa_n}> 1)=0~.
\end{align}
Now suppose that $F$ is non-strictly concave. In this case, consistency of the rescaled bootstrap is justified by Proposition 2 in \citet{Dumbgen1993}. Alternatively, the bootstrap consistency follows by Lemma S.3.8 (for consistency of $\hat{\mathcal M}_n'$ as an derivative estimator of $\mathcal M_F'$) and Theorem 3.2 in \citet{FangSantos2018HDD}. In any case, $\hat c_{1-\alpha}^*\convp c_{1-\alpha}^*$ by Assumption \ref{Ass: quantile regularity} and Corollary 3.2 in \citet{FangSantos2015HDDarxiv}. Therefore, it follows that
\begin{align}
\limsup_{n\to\infty}P(\sqrt{n}\phi(\mathbb F_n)> \hat c_{1-\alpha})\le &\limsup_{n\to\infty} P(\sqrt{n}\phi(\mathbb F_n)> \hat c_{1-\alpha}^*)\notag\\
&\le P(\|\mathcal M_F'(\mathbb G)-\mathbb G\|_p\ge c_{1-\alpha})=\alpha~,
\end{align}
where the second inequality follows from Slutsky's lemma (so that $\sqrt{n}\phi(\mathbb F_n)- \hat c_{1-\alpha}^*\convl \|\mathcal M_F'(\mathbb G)-\mathbb G\|_p- c_{1-\alpha}$ ) and the portmanteau theorem, and the last step is due to Assumption \ref{Ass: quantile regularity}. This proves the first claim of the theorem.

We now turn to the second part of the theorem and suppose that $F$ is not concave. First, we show that $\hat c_{1-\alpha}^*=O_{p}(1)$ regardless of whether $F$ is concave or not, where $O_p(1)$ means ``bounded in probability.'' To this end, note that by Lemma 2.2 in \citet{Durot_Tocquet2003}, for any $h\in\ell^\infty(\mathbf R_+)$,
\begin{align}\label{Eqn: pointwise thm, aux1}
\|\hat{\mathcal M}_n'(h))\|_p\le \|\frac{\mathcal M(\mathbb F_n+t_nh)-\mathcal M(\mathbb F_n)}{t_n}\|_p\lesssim \|h\|_\infty~,
\end{align}
where $\lesssim$ means `` smaller than or equal to up to a universal constant''. Here, the constant only depends on the (known) weighting function $g$ if $p\in[1,\infty)$. In turn, result \eqref{Eqn: pointwise thm, aux1} implies that
\begin{align}\label{Eqn: pointwise thm, aux2}
\|\hat{\mathcal M}_n'(\sqrt n\{\mathbb F_n^*-&\mathbb F_n\})-\sqrt n\{\mathbb F_n^*-\mathbb F_n\}\|_p\notag\\
&\le \|\hat{\mathcal M}_n'(\sqrt n\{\mathbb F_n^*-\mathbb F_n\})\|_p+\|\sqrt n\{\mathbb F_n^*-\mathbb F_n\}\|_p\notag\\
&\lesssim 2\|\sqrt n\{\mathbb F_n^*-\mathbb F_n\}\|_\infty~.
\end{align}
Moreover, Theorem 3.1 in \citet{GineandZin1990} and Proposition 10.7 in \citet{Kosorok2008} yields $\|\sqrt n\{\mathbb F_n^*-\mathbb F_n\}\|_\infty=O_{p}(1)$ outer almost surely. This, together with result \eqref{Eqn: pointwise thm, aux2} and Lemma 3 in \citet{Cheng_Huang2010bootstrap}, implies that:
\begin{align}\label{Eqn: pointwise thm, aux3}
\|\hat{\mathcal M}_n'(\sqrt n\{\mathbb F_n^*-\mathbb F_n\})-\sqrt n\{\mathbb F_n^*-\mathbb F_n\}\|_p=O_{p}(1)~,
\end{align}
unconditionally. Fix $\epsilon\in(0,1)$. Then we may choose some $M>0$ such that
\begin{align}\label{Eqn: pointwise thm, aux4}
P(\|\hat{\mathcal M}_n'(\sqrt n\{\mathbb F_n^*-\mathbb F_n\})-\sqrt n\{\mathbb F_n^*-\mathbb F_n\}\|_p> M)\le \alpha\epsilon~.
\end{align}
By the definition of $\hat c_{1-\alpha}^*$, if $\hat c_{1-\alpha}^*>M$, then we must have
\begin{align}\label{Eqn: pointwise thm, aux5}
P_W(\|\hat{\mathcal M}_n'(\sqrt n\{\mathbb F_n^*-\mathbb F_n\})-\sqrt n\{\mathbb F_n^*-\mathbb F_n\}\|_p> M)> \alpha~.
\end{align}
We may then conclude from the implication of \eqref{Eqn: pointwise thm, aux5} that
\begin{align}\label{Eqn: pointwise thm, aux5a}
P(\hat c_{1-\alpha}^* & >M)\le P_X(P_W(\|\hat{\mathcal M}_n'(\sqrt n\{\mathbb F_n^*-\mathbb F_n\})-\sqrt n\{\mathbb F_n^*-\mathbb F_n\}\|_p> M)> \alpha)\notag\\
&\le \alpha^{-1} E_X [P_W(\|\hat{\mathcal M}_n'(\sqrt n\{\mathbb F_n^*-\mathbb F_n\})-\sqrt n\{\mathbb F_n^*-\mathbb F_n\}\|_p> M)]\notag\\
&\le \alpha^{-1} P(\|\hat{\mathcal M}_n'(\sqrt n\{\mathbb F_n^*-\mathbb F_n\})-\sqrt n\{\mathbb F_n^*-\mathbb F_n\}\|_p> M)\le\epsilon~,
\end{align}
where the second inequality follows by the Markov's inequality, the third inequality by Lemma 1.2.6 in \citet{Vaart1996}, and the last step by the result \eqref{Eqn: pointwise thm, aux4}. This shows that $\hat c_{1-\alpha}^*=O_{p}(1)$ and hence $\hat c_\alpha=\max\{\kappa_n,\hat c_{1-\alpha}^*\}=O_p(1)$ in view of $\kappa_n=o(1)$ by Assumption \ref{Ass: kn}(ii).

Next, by the triangle inequality and Lemma 2.2 in \citet{Durot_Tocquet2003},
\begin{align}\label{Eqn: pointwise thm, aux6}
\sqrt{n}\phi(\mathbb F_n)& =\sqrt n \|\mathcal M\mathbb F_n-\mathcal M F -\{\mathbb F_n-F\}+\mathcal M F-F\|_p\notag\\
&\ge \sqrt n\|\mathcal M F-F\|_p- \{\sqrt n \|\mathcal M\mathbb F_n-\mathcal M F\|_p +\sqrt n \|\mathbb F_n-F\|_p\}\notag\\
&\ge \sqrt n\|\mathcal M F-F\|_p-C \|\sqrt n\{\mathbb F_n-F\}\|_\infty~,
\end{align}
where $C>0$ is a constant depending on the weighting function $g$. Note that if $F$ is noncacave, then $\mathcal M F$ and $F$ must differ from each other on a set of positive Lebesgue measure, implying that $\|\mathcal M F-F\|_p>0$. Together with
$\|\sqrt n\{\mathbb F_n-F\}\|_\infty=O_{p}(1)$, we may then conclude from result \eqref{Eqn: pointwise thm, aux6} that $\sqrt{n}\phi(\mathbb F_n)\to\infty$ in probability. Combining with $\hat c_\alpha=O_p(1)$, it follows that
\begin{align}
P(\sqrt n \phi(\mathbb F_n)>\hat c_\alpha)\to 1~.
\end{align}
This completes the proof for the second claim of the theorem. \qed

\noindent{\sc Proof of Theorem \ref{Thm: test concavity, local}:} First, for convenience of the reader, we introduce the concept of contiguity. For each $n\in\mathbf N$, let $P_n$ and $Q_n$ be two generic probability measures on some measurable space $(\Omega_n,\mathcal A_n)$. Then the sequence $\{Q_n\}$ is said to be contiguous with respect to $\{P_n\}$ if, for any statistic $T_n:\Omega_n\to\mathbf R$, one has $T_n\convp 0$ under $Q_n$ whenever $T_n\convp 0$ under $P_n$. Heuristically, this means that any statistic $T_n$ that is asymptotically negligible under $P_n$ remains so under $Q_n$. We say that $\{Q_n\}$ and $\{P_n\}$ are mutually contiguous if $\{Q_n\}$ is contiguous with respect to $\{P_n\}$ and vice versa. We refer the reader to Chapter 6 in \citet{Vaart1998}---in particular, Lemma 6.4 there---for more details.

Since $\{P_t\}$ is a differentiable path, it follows by Theorem 12.2.3 and Corollary 12.3.1 in \citet{TSH2005} that $P_n^n$ and $P^n$ are mutually contiguous. By Theorem 2.1 in \citet{FangSantos2018HDD} and $\phi(F(P))=0$, we then have:
\begin{multline}\label{Eqn: local size, aux0}
\sqrt{n}\phi(\mathbb F_n)=\sqrt{n}\{\phi(\mathbb F_n)-\phi(F(P))\}\\
=\phi_F'(\sqrt{n}\{\mathbb F_n-F(P)\})+o_p(1)\text{ under }P^n~,
\end{multline}
where $\phi_F'(h)=\|\mathcal M_F'h-h\|_p$. In turn, we obtain by \eqref{Eqn: local size, aux0} and mutual continuity of $P_n^n$ and $P^n$ that, under $P_n^n$,
\begin{align}\label{Eqn: local size, aux1}
\sqrt{n}\phi(\mathbb F_n)=\phi_F'(\sqrt{n}\{\mathbb F_n-F(P)\})+o_p(1)~.
\end{align}
As is well known in the literature (see, for example, \citet[p.192]{BKRW993Efficient}), $\mathbb F_n$ is a regular estimator of $F$ so that, under $P_n^n$,
\begin{align}\label{Eqn: local size, aux2}
\sqrt{n}\{\mathbb F_n-F(P_{\eta/\sqrt{n}})\}\convl \mathbb G \text{ in }\ell^\infty(\mathbf R_+)~.
\end{align}
By result \eqref{Eqn: F diff}, we also have, as a deterministic result,
\begin{align}\label{Eqn: local size, aux2a}
\sqrt n \{F(P_{\eta/\sqrt{n}})-F(P)\}\to \eta\dot F(h) \text{ in }\ell^\infty(\mathbf R_+)~.
\end{align}
By Slutsky's theorem, we obtain from results \eqref{Eqn: local size, aux2} and \eqref{Eqn: local size, aux2a} that
\begin{align}\label{Eqn: local size, aux2b}
\sqrt{n}\{\mathbb F_n-F(P)\} \convl \mathbb G + \eta\dot F(h) \text{ in }\ell^\infty(\mathbf R_+)~,
\end{align}
under $P_n^n$. Combining \eqref{Eqn: local size, aux1} and \eqref{Eqn: local size, aux2b}, together with the continuous mapping theorem and Slutsky's theorem, we may thus conclude that, under $P_n^n$,
\begin{align}\label{Eqn: local size, aux2c}
\sqrt{n}\phi(\mathbb F_n)\convl \phi_F'(\mathbb G + \eta\dot F(h))\equiv \|\mathcal M_F'(\mathbb G+\eta\dot F(h))-\{\mathbb G+\eta\dot F(h)\}\|_p~.
\end{align}

For the first part of the theorem, we consider the two cases separately.

\noindent\underline{\sc Case I:} $F$ is strictly concave. By Assumption \ref{Ass: kn}(ii), $n^{-2/3}(\log n)^{2/3}/(n^{-1/2}\kappa_n)\to 0$ as $n\to\infty$. By Assumption \ref{Ass: kn}(ii) and Theorem \ref{Thm: Equivalence, local} we thus have
\begin{multline}\label{Eqn: local size, aux3}
\liminf_{n\to\infty}\pi_n(P_{1/\sqrt{n}})=\liminf_{n\to\infty}P_n^n ( \sqrt n\|\hat{\mathbb F}_n-\mathbb F_n\|_p >\max\{\kappa_n,\hat c_{1-\alpha}\})\\
\le \liminf_{n\to\infty}P_n^n (\frac{\sqrt n\|\hat{\mathbb F}_n-\mathbb F_n\|_p}{\kappa_n}>1)=0~.
\end{multline}

\noindent\underline{\sc Case II:} $F$ is non-strictly concave. Mutual contiguity of $P_n^n$ and $P^n$, together with $\hat c_{1-\alpha}^*\convp c_{1-\alpha}^*$ under $P^n$ from the proof of Theorem \ref{Thm: test concavity, pointwise}, implies that whenever $F$ is non-strictly concave,
\begin{align}\label{Eqn: local size, aux3a}
\hat c_{1-\alpha}^*\convp c_{1-\alpha}^* \text{ under }P_n^n~.
\end{align}
Since $c_{1-\alpha}^*$ is nonnegative as a quantile of $\|\mathcal M_F'(\mathbb G)-\mathbb G\|_p$ and $\kappa_n=o(1)$ by Assumption \ref{Ass: kn}(ii), we obtain from \eqref{Eqn: local size, aux3a} and the continuous mapping theorem that $\hat c_{1-\alpha}\convp c_{1-\alpha}^*$ under $P_n^n$. This, together with result \eqref{Eqn: local size, aux2c}, Slutsky's theorem and the portmanteau theorem, allows us to conclude that
\begin{align}\label{Eqn: local size, aux4}
\liminf_{n\to\infty}\pi_n(P_{1/\sqrt{n}})&=\liminf_{n\to\infty}P_n^n (\sqrt{n}\phi(\mathbb F_n)>\hat c_{1-\alpha})\notag\\
&=P(\|\mathcal M_F'(\mathbb G+\eta\dot F(h))-\{\mathbb G+\eta\dot F(h)\}\|_p>c_{1-\alpha}^*)~.
\end{align}
If in addition the cdf of $\|\mathcal M_F' (\mathbb G+\eta\dot F(h))-\{\mathbb G+\eta\dot F(h)\}\|_p$ is continuous at $c_{1-\alpha}^*$, then \eqref{Eqn: local size, aux4} holds with equality, again by the portmanteau theorem. This completes the proof for the first part of the theorem.

As for the second part, we again consider the following two cases separately.

\noindent\underline{\sc Case I:} $F$ is strictly concave. By Assumption \ref{Ass: kn}(ii) and Theorem \ref{Thm: Equivalence, local} we have
\begin{align}\label{Eqn: local size, aux5}
\limsup_{n\to\infty}\pi_n(P_{1/\sqrt{n}})\le \limsup_{n\to\infty}P_n^n (\sqrt n\|\hat{\mathbb F}_n-\mathbb F_n\|_p>\kappa_n)=0~.
\end{align}
This is in fact also implied by result \eqref{Eqn: local size, aux3} which holds for any $\eta\in\mathbf R$.

\noindent\underline{\sc Case II:} $F$ is non-strictly concave. We follow the arguments of Theorem 3.3 in \citet{FangSantos2018HDD}. To begin with, note that result \eqref{Eqn: local size, aux2c}, $\hat c_{1-\alpha}\convp c_{1-\alpha}^*$ under $P_n^n$ (as argued previously), Slutsky's theorem and the portmanteau theorem together imply that
\begin{align}\label{Eqn: local size, aux6}
\limsup_{n\to\infty}\pi_n(P_{1/\sqrt{n}})&\le \limsup_{n\to\infty}P_n^n (\sqrt{n}\phi(\mathbb F_n)>\hat c_{1-\alpha}^*)\notag\\
&\le P(\|\mathcal M_F'(\mathbb G+\eta\dot F(h))-\{\mathbb G+\eta\dot F(h)\}\|_p\ge c_{1-\alpha}^*)~.
\end{align}

Next, we show that the map $g\mapsto\phi_F'(g)=\|\mathcal M_F'(g)-g\|_p$ (defined on $\ell^\infty(\mathbf R_+)$) is subadditive. Note that (i) $\mathcal M_F'(g)-g\ge 0$ for all $g\in\ell^\infty(\mathbf R_+)$ because $\mathcal M_F'(g)$ majorizes $g$, and (ii) $\ell^\infty(\mathbf R_+)$ can be identified as a subspace of $L^p(\mathbf R_+)$ equipped with norm $\|\cdot\|_p$ if $p\in[1,\infty)$---of course, if $p=\infty$, then $\ell^\infty(\mathbf R_+)$ is a subspace of itself. Thus, we may view $\psi_F'$ as a real-valued map defined on some (normed) subspace of $L^p(\mathbf R_+)$ if $p\in[1,\infty)$ or of $\ell^\infty(\mathbf R_+)$ if $p=\infty$. Since $\phi_F'$ is defined in term of the $\|\cdot\|_p$ norm, such an embedding (of the domain $\ell^\infty(\mathbf R_+)$ of $\phi_F'$ into $L^p(\mathbf R_+)$) allows us to rewrite $\phi_F'$ according to Lemma \ref{Lem: Riesz}.

Specifically, let $\ell^\infty(\mathbf R_+)^{p*}$ be the topological dual space of $\ell^\infty(\mathbf R_+)$ viewed as a subspace of $L^p(\mathbf R_+)$ if $p\in[1,\infty)$, and be the topological dual space of $\ell^\infty(\mathbf R_+)$ if $p=\infty$.\footnote{Let $\mathbb D$ be a normed space with norm $\|\cdot\|_{\mathbb D}$ such as $L^p(\mathbf R_+)$ and $\ell^\infty(\mathbf R_+)$. The topological dual space $\mathbb D^*$ of $\mathbb D$ is the space of all continuous and linear functions $\varphi:\mathbb D\to\mathbf R$, equipped with the operator norm $\|\varphi\|_{op}\equiv\sup_{\|g\|_{\mathbb D}\le 1}|\varphi(g)|$.} Define $\mathbb S_+^p=\{\varphi\in \ell^\infty(\mathbf R_+)^{p*}: \|\varphi\|_{op}=1,\varphi\ge 0\}$ where $\|\cdot\|_{op}$ is the operator norm, and $\varphi\ge 0$ means that $\varphi(g)\ge 0$ whenever $g(x)\ge 0$ for all $x\in\mathbf R_+$. By Lemma \ref{Lem: Riesz} we now have: for any $g_1,g_2\in\ell^\infty(\mathbf R_+)$,
\begin{multline}\label{Eqn: local size, aux7}
\phi_F'(g_1 + g_2)=\sup_{\varphi\in \mathbb S_+^p}\varphi(\mathcal M_F'(g_1+g_2)-(g_1+g_2)) \\
=\sup_{\varphi\in \mathbb S_+^p}\{\varphi(\mathcal M_F'(g_1+g_2))-[\varphi(g_1)+\varphi(g_2)]\} ~,
\end{multline}
where the second equality follows by linearity of each $\varphi\in \mathbb S_+^p$. Since $\mathcal M_F'$ is convex by Proposition 2.2 in \citet{Beare_Fang2016Grenander} and positively homogeneous of degree one as a Hadamard directional derivative \citep{Shapiro1990}, it follows that $\mathcal M_F'$ must be subadditive, i.e.,
\begin{align}\label{Eqn: local size, aux7a}
\mathcal M_F'(g_1+g_2)\le \mathcal M_F'(g_1)+\mathcal M_F'(g_2)~.
\end{align}
Since each $\varphi\in \mathbb S_+^p$ is linear and satisfies $\varphi(g)\ge 0$ whenever $g(x)\ge 0$ for all $x\in\mathbf R_+$, we obtain from \eqref{Eqn: local size, aux7a} that, for any $\varphi\in \mathbb S_+^p$,
\begin{align}\label{Eqn: local size, aux7b}
\varphi(\mathcal M_F'(g_1)+\mathcal M_F'(g_2))\le \varphi(\mathcal M_F'(g_1)) + \varphi(\mathcal M_F'(g_2))~.
\end{align}
Combining results \eqref{Eqn: local size, aux7} and \eqref{Eqn: local size, aux7b} then leads to
\begin{align}\label{Eqn: local size, aux7c}
\phi_F'(g_1 +g_2) &\le \sup_{\varphi\in \mathbb S_+^p}\{\varphi(\mathcal M_F'(g_1))+\varphi(\mathcal M_F'(g_2))-[\varphi(g_1)+\varphi(g_2)]\}\notag\\
&\le \sup_{\varphi\in \mathbb S_+^p}\{\varphi(\mathcal M_F'(g_1))-\varphi(g_1)\}+\sup_{\varphi\in \mathbb S_+^p}\{\varphi(\mathcal M_F'(g_2))-\varphi(g_2)\}\notag\\
&=\sup_{\varphi\in \mathbb S_+^p}\{\varphi(\mathcal M_F'(g_1)-g_1)\}+\sup_{\varphi\in \mathbb S_+^p}\{\varphi(\mathcal M_F'(g_2)-g_2)\}~,
\end{align}
where the equality again follows by linearity of each $\varphi\in \mathbb S_+^p$. By Lemma \ref{Lem: Riesz}, we in turn have from \eqref{Eqn: local size, aux7c} that, for any $g_1,g_2\in\ell^\infty(\mathbf R_+)$,
\begin{multline}\label{Eqn: local size, aux7d}
\phi_F'(g_1 +g_2) \le  \|\mathcal M_F'(g_1)-g_1\|_p+\|\mathcal M_F'(g_1)-g_1\|_p\\
=\phi_F'(g_1)+\phi_F'(g_2)~,
\end{multline}

In light of result \eqref{Eqn: local size, aux7}, we immediately obtain that
\begin{align}\label{Eqn: local size, aux8}
P(\|\mathcal M_F'(\mathbb G&+\eta\dot F(h))-\{\mathbb G+\eta\dot F(h)\}\|_p\ge c_{1-\alpha}^*)\notag\\
&\le P(\|\mathcal M_F'(\mathbb G)-\mathbb G\|_p+\|\mathcal M_F'(\eta\dot F(h))-\eta\dot F(h)\|_p\ge c_{1-\alpha}^*)~.
\end{align}
Moreover, since $\{P_{\eta/\sqrt n}\}$ is a differentiable path under the null when $\eta\le 0$, it follows that $\phi(F(P_{\eta/\sqrt n}))=\phi(F(P))=0$ for all $n\in\mathbf N$ and $\eta\le 0$. Hence,
\begin{multline}\label{Eqn: local size, aux9}
0=\lim_{n\to\infty}\frac{\phi(F(P_{\eta/\sqrt n}))-\phi(F(P))}{n^{-1/2}}\\
=\phi_F'(\eta\dot F(h))\equiv \|\mathcal M_F'(\eta\dot F(h))-\eta\dot F(h)\|_p~.
\end{multline}
Results \eqref{Eqn: local size, aux6}, \eqref{Eqn: local size, aux8} and \eqref{Eqn: local size, aux9} then yield that
\[
\limsup_{n\to\infty}\pi_n(P_{1/\sqrt{n}})\le P(\|\mathcal M_F'(\mathbb G)-\mathbb G\|_p\ge c_{1-\alpha}^*)=\alpha~,
\]
where we also exploited that $c_{1-\alpha}^*$ is a continuity point of the cdf of $\|\mathcal M_F'(\mathbb G)-\mathbb G\|_p$ by Assumption \ref{Ass: quantile regularity}. This completes the proof of the second part.

Finally, for the third part of the theorem, we note by Lemma 2.2 in \citet{Durot_Tocquet2003} that
\begin{align}\label{Eqn: local size, aux10}
|\sqrt n \phi(\mathbb F_n)-\sqrt n\phi(F_n)| \lesssim \|\sqrt n \{\mathbb F_n-F_n\}\|_p
\lesssim  \sqrt n\| \mathbb F_n-F_n\|_\infty~,
\end{align}
where the second inequality follows from the definition of $\|\cdot\|_p$. Fix $M>0$. By Dvoretzky, Kiefer and Wolfowitz's inequality (see, for example, Theorem 11.2.18 in \citet{TSH2005}), we have that
\begin{align}\label{Eqn: local size, aux11}
P_n^n(\sqrt n\| \mathbb F_n-F_n\|_\infty>M)= P_n(\| \mathbb F_n-F_n\|_\infty>\frac{M}{\sqrt n})\lesssim \exp\{-2 M^2\}~,
\end{align}
which can be made arbitrarily small by choosing a sufficiently large $M$. Combining results \eqref{Eqn: local size, aux10} and \eqref{Eqn: local size, aux11} we thus obtain that
\begin{align}\label{Eqn: local size, aux12}
|\sqrt n \phi(\mathbb F_n)-\sqrt n\phi(F_n)| =O_p(1) \text{ under } P_n^n~.
\end{align}
Moreover, by Proposition 2.1 in \citet{Beare_Fang2016Grenander} and $\phi(F)=0$ (by the definition of $\mathbf H$), we obtain that, as $n\to\infty$,
\begin{align}\label{Eqn: local size, aux13}
\sqrt n \phi(F_n) =\sqrt n \{\phi(F_n)- \phi(F)\} \to \|\mathcal M_F'(\eta\dot F(h))-\eta\dot F(h)\|_p~.
\end{align}
By the positive homogeneity of degree one of $\mathcal M_F'$ (as a Hadamard directional derivative), we thus have from \eqref{Eqn: local size, aux13} that: for $\Delta\equiv \|\mathcal M_F'(\dot F(h))-\dot F(h)\|_p$,
\begin{align}\label{Eqn: local size, aux14}
\sqrt n \phi(F_n)= \eta \Delta+o(1)~,
\end{align}
for all $\eta>0$. It follows from \eqref{Eqn: local size, aux12} and \eqref{Eqn: local size, aux14} that, under $P_n^n$,
\begin{align}\label{Eqn: local size, aux15}
\sqrt n \phi(\mathbb F_n) = \sqrt n \phi(\mathbb F_n)-\sqrt n\phi(F_n) + \sqrt n \phi(F_n) = \eta \Delta +O_p(1)~.
\end{align}

Having derived the order of the test statistic as in \eqref{Eqn: local size, aux15}, we next evaluate the order of the critical value. For this, note that, by result \eqref{Eqn: pointwise thm, aux2},
\begin{align}\label{Eqn: local size, aux16}
\|\hat{\mathcal M}_n'(\sqrt n\{\mathbb F_n^*-\mathbb F_n\})-\sqrt n\{\mathbb F_n^*-\mathbb F_n\}\|_p\lesssim 2\sqrt n \|\mathbb F_n^*-\mathbb F_n\|_\infty~.
\end{align}
By Dvoretzky, Kiefer and Wolfowitz's inequality (see, for example, Theorem 11.2.18 in \citet{TSH2005}), we note that, for any $M>0$,
\begin{align}\label{Eqn: local size, aux17}
P_W(\sqrt n \|\mathbb F_n^*-\mathbb F_n\|_\infty>M) \lesssim \exp\{-2 M^2\}~,
\end{align}
which can be made arbitrarily small by choosing a sufficiently large $M$. Result \eqref{Eqn: local size, aux17}, together with arguments similar to those leading to \eqref{Eqn: pointwise thm, aux5a}, yields
\begin{align}\label{Eqn: local size, aux18}
\hat c_{1-\alpha^*} = O_p(1)\text{ under }P_n^n~.
\end{align}
Combining results \eqref{Eqn: local size, aux15} and \eqref{Eqn: local size, aux18}, we may thus conclude
\begin{align}
P_n^n(\sqrt n \phi(\mathbb F_n)>\hat c_{n,1-\alpha}) = P_n^n (\eta \Delta+O_p(1)>\max\{\kappa_n,\hat c_{1-\alpha}^*\}) = 1~,
\end{align}
by letting $n\to\infty$ followed by $\eta\uparrow\infty$, where we exploited $\kappa_n=o(1)$ by Assumption \ref{Ass: kn}(ii) and $\Delta>0$. This completes the proof of the theorem. \qed

Our final lemma entails some knowledge on Riesz space. For convenience of the reader, we introduce some relevant concepts in this regard. Just like normed spaces generalize to abstract spaces the operations of vector addition and scalar multiplication in Euclidean spaces, Riesz spaces generalize the binary relations $\ge$ and $\le$. Let $\mathbb E$ be a set. We say that $\mathbb E$ is partially ordered if there exists a binary relation $\ge$ that is transitive (i.e., $x\ge z$ whenever $x\ge y$ and $y\ge z$), reflexive (i.e., $x\ge x$) and antisymetric (i.e., $x\ge y$ and $y\ge x$ implies $x=y$). The notation $y\le x$ is equivalent to $x\ge y$. If $x\ge y$ but $x\neq y$, we also write $x>y$ or equivalently $y<x$. Let $\mathbb E$ be partially ordered by $\ge$. An element $z$ is the supremum of a pair $x,y\in\mathbb E$, denoted $x\vee y$, if $x\ge z$, $y\ge z$, and $z\le u$ whenever $x\le u$ and $y\le u$. The infimum of $x,y\in\mathbb E$, denoted $x\wedge y$, is defined similarly. Not every pair $x,y$ in $\mathbb E$ admits a supremum or infimum. But, if this is the case, then we call $\mathbb E$ a lattice.

If $\mathbb E$ is a vector space, then one may hope that the partial order $\ge$ is ``compatible'' with the algebraic structure of $\mathbb E$. Concretely, whenever $x\ge y$, one should have (i) $x+z\ge y+z$ for all $z\in\mathbb E$, and (ii) $\alpha x\ge \alpha y$ for all $\alpha\in\mathbf R_+$. In this case, we call $\mathbb E$ an (partially) ordered vector space. A partially ordered vector space that is also a lattice is called a Riesz space. For example, for $L^p(\mathbf R_+)$ with $p\in[1,\infty)$, we may define the partial order $\ge $ by: $f\ge g$ whenever $f(x)\ge g(x)$ for almost all $x\in\mathbf R_+$. For $\ell^\infty(\mathbf R_+)$, we define $\ge $ analogously: $f\ge g$ whenever $f(x)\ge g(x)$ for all $x\in\mathbf R_+$. Then both $L^p(\mathbf R^+)$ and $\ell^\infty(\mathbf R_+)$ are Riesz spaces. Given the notion of supremum, we may define the absolute value of $x\in\mathbb E$ by $|x|\equiv x^++x^-$ where $x^+\equiv x\vee 0$ and $x^-\equiv (-x)\vee 0$. If $\mathbb E$ is a Riesz space and is equipped with a norm $\|\cdot\|$ satisfying the property that $\|x\|\le \|y\|$ whenever $|x|\le |y|$, then $\mathbb E$ is called a normed Riesz space, and the norm $\|\cdot\|$ is called a lattice norm. We refer the reader to \citet{AliprantisandBorder2006} for more discussions. With these concepts in hand, we now have:

\begin{lem}\label{Lem: Riesz}
Let $\mathbb E$ be a normed Riesz space with partial order $\ge $ and lattice norm $\|\cdot\|_{\mathbb E}$. Then for any $x\ge 0$, we have
\[
\|x\|_{\mathbb E}=\sup_{\varphi\in \mathbb E^*: \|\varphi\|_{op}=1,\varphi\ge 0}\varphi(x)~,
\]
where $\mathbb E^*$ is the topological dual space of $\mathbb E$, $\|\cdot\|_{op}$ is the operator norm, and $\varphi\ge 0$ means that $\varphi(y)\ge 0$ whenever $y\ge 0$.
\end{lem}
\noindent{\sc Proof:} The conclusion trivially holds if $x=0$. Suppose that $x>0$. Then Theorem 39.3 in \citet{Zaanen1997} implies that there exists some $\varphi^*\in \mathbb E^*$ with $\|\varphi^*\|_{op}=1$ and $\varphi^*\ge 0$ such that $\varphi^*(x)=\|x\|_{\mathbb E}$. This in turn implies that
\[
\|x\|_{\mathbb E}=\varphi^*(x)\le \sup_{\varphi\in \mathbb E^*: \|\varphi\|_{op}=1,\varphi\ge 0}\varphi(x)~.
\]
On the other hand, we have
\[
\sup_{\varphi\in \mathbb E^*: \|\varphi\|_{op}=1,\varphi\ge 0}\varphi(x)\le \sup_{\varphi\in \mathbb E^*: \|\varphi\|_{op}=1,\varphi\ge 0}\|\varphi\|_{op}\|x\|_{\mathbb E}=\|x\|_{\mathbb E}~.
\]
This completes the proof of the lemma. \qed

\section{Some Extensions}\label{App: extension}
\renewcommand{\theequation}{B.\arabic{equation}}
\setcounter{equation}{0}

In this section, we discuss extensions of our test results to a general setup that includes inference on density, regression and hazard functions as special cases. To this end, We shall assume that the data $\{X_i\}_{i=1}^n$ are i.i.d. with common distribution $P_{\eta/\sqrt n}\in \mathbf P$, which we also denote by $P_n$ with some abuse of notation. Here, $\mathbf P$ denotes the model as before. This allows us to consider the standard i.i.d.\ setup (i.e., $\eta=0$) as well as the setup for local analysis, in a unified way. In turn, wherever appropriate, we identify the parameter of interest $\theta$ as a map $\theta: \mathbf P\to\ell^\infty([a,b])$ such that, for any $P\in\mathbf P$ and $t\in[a,b]$,
\begin{align}
\theta(P)(t)=\int_a^t g(u)\,\mathrm du~,
\end{align}
where $a,b\in\mathbf R$ are known with $a<b$ throughout. The dependence of $\theta(P)$ on $P$ sometimes is also suppressed for simplicity. While it may be possible to consider settings with unbounded support such as $[0,\infty)$, we confine our attention to the bounded case for simplicity---it is not essential to what we intend to convey anyway. Our objective is to test the null that $\theta(P_{\eta/\sqrt n})$ is concave versus otherwise. In the main text, $g$ is the density function and $\theta$ is the cdf. Below we formalize the hazard and the regression examples.

\begin{ex}[Monotone Hazard Rate]\label{Ex: hazard}
Let $X\in \mathbf R_+$ be a random variable with cdf $F$ that admits a pdf $f$ with respect to the Lebesgue measure. Then the hazard rate $\lambda$ at $u$ is defined by:
\begin{align}
\lambda(u)\equiv \frac{f(u)}{1-F(u)}~.
\end{align}
Heuristically, $\lambda(u)$ measures the probability of instantaneous failure rate at time $u$, given the subject has functioned until $u$ \citep{GroeneboomJongbloed2014Shape}. A leading example of shape restrictions is that $\lambda$ be monotonically increasing on an interval $[0,b]$ for some $b\in(0,\infty)$. This can be inferred from convexity of the cumulative hazard function $\Lambda$ defined by: for any $t\in[0,b]$,
\begin{align}
\Lambda(t)\equiv\Lambda_F(t)\equiv \int_0^t \frac{1}{1-F_-(u)}\,\mathrm dF(u)~,
\end{align}
where $u\mapsto F_-(u)\equiv\lim_{t\uparrow u}F(t)$ is the left-continuous version of $F$. Here, $F(u)=F(u-)$ for all $u\in\mathbf R_+$ since $F$ is absolutely continuous, so $\Lambda(t)=\int_0^t \lambda(u)\,\mathrm du$. The definition above, however, allows us to construct estimators of $\Lambda$ in a straightforward way. In this example, $\theta=-\Lambda$.  \qed
\end{ex}

\begin{ex}[Isotonic/Monotone Regression]\label{Ex: regression}
Let $\{(Y_i,Z_i)\}_{i=1}^n$ be i.i.d.\ bivariate random vectors generated according to
\begin{align}
Y_i = m(Z_i) + \epsilon_i~,
\end{align}
where each $Z_i$ is supported on $[a,b]$, the regression function $m: [a,b]\to\mathbf R$ is unknown, and each $\epsilon_i$ is a centered random variable independent of $Z_i$. Let $Q:[0,1]\to[a,b]$ be the population quantile function of $Z_1$. Then monotonicity of $m$ can be inferred from convexity of the population cumulative regression function (cqf) $\theta\in\ell^\infty([0,1])$ defined by: for all $t\in[0,1]$,
\begin{align}\label{Eqn: cqf}
\Lambda(t)=\int_0^t m(Q(u))\,\mathrm du~.
\end{align}
One may also view $\Lambda$ as a population analog of the cumulative sum diagram---see, for example, \citet[p.3865]{Beare_Fang2016Grenander}. Here, $\theta=-\Lambda$. \qed
\end{ex}

Since a function $\theta:[a,b]\to\mathbf R$ is convex if and only if $-\theta$ is concave, both Examples \ref{Ex: regression} and \ref{Ex: hazard} fall within the scope of our framework. For a coherent treatment, we shall thus maintain the null hypothesis that $\theta$ is concave versus otherwise, and revisit Examples \ref{Ex: regression} and \ref{Ex: hazard} after the general theory is presented. Towards this end, let $\hat\theta_n:\{X_i\}_{i=1}^n\to\ell^\infty([a,b])$ be an unconstrained estimator of $\theta$. Note that $X_i$ is a generic notation for the $i$-th observation; e.g, $X_i=(Y_i,Z_i)'$ in Example \ref{Ex: regression}. In turn, we then employ the test statistic $\sqrt n \phi(\hat\theta_n)$ where $\phi: \ell^\infty([a,b])\to\mathbf R$ is defined by
\begin{align*}
\phi(\theta)=\|\mathcal M\theta-\theta\|_p=\begin{cases}
[\int_a^b(\mathcal M\theta(x)-\theta(x))^p g(x)\, \mathrm dx]^{1/p} & \text{ if }p\in[1,\infty)\\
\sup_{a\le x\le b}|\mathcal M\theta(x)-\theta(x)| & \text{ if }p=\infty
\end{cases}~,
\end{align*}
for some positive weighting function $g\in L^1([a,b])$. Throughout, we study the test functional $\phi$ with an arbitrarily fixed $p\in[1,\infty]$.

We now proceed by imposing the following assumptions, where $C([a,b])$ is the space of continuous functions on $[a,b]$ equipped with the uniform norm.

\begin{ass}\label{Ass: DGP}
(i) $\{X_i\}_{i=1}^n$ are i.i.d.\ with a common distribution $P_{\eta/\sqrt n}$ for some $\eta\in\mathbf R$; (ii) $\{P_t: |t|<\epsilon\}\subset\mathbf P$ is a differentiable path with score $h$ in the sense of Definition \eqref{Defn: differentiable path}; (iii) $\|\sqrt n \{\theta(P_{\eta/\sqrt n})-\theta(P_0)\}-\eta\dot\theta(h)\|_{\infty}\to 0$ for some $\eta\dot\theta(h)\in C([a,b])$; (iv) $\theta_0\equiv\theta(P_{0})\in \ell^\infty([a,b])$ is concave.
\end{ass}

\begin{ass}\label{Ass: estimator}
(i) $\hat\theta_n:\{X_i\}_{i=1}^n\to\ell^\infty([a,b])$ satisfies $\sqrt n \{\hat\theta_n-\theta_n\}\convl \mathbb G$ in $\ell^\infty([a,b])$ under $P_n^n\equiv \prod_{i=1}^{n} P_{\eta/\sqrt n}$, where $\theta_n\equiv \theta(P_{\eta/\sqrt n})$ and $\mathbb G\in C([a,b])$; (ii) $\|\mathcal M\hat\theta_n-\hat\theta_n\|_\infty=O_p(s_n)$ under $P_n^n$ if $\theta_0$ is strictly concave for some sequence $\{s_n\}$ of strictly positive scalars satisfying $\sqrt n s_n=o(1)$.
\end{ass}

Assumption \ref{Ass: DGP}(i)(ii) formalizes the data generating process. Assumption \ref{Ass: DGP}(iii) is a generalization of the property \eqref{Eqn: F diff}, and is fulfilled whenever $\theta: \mathbf P\to\ell^\infty([a,b])$ is a regular parameter in the sense of Definition 5.1 in \citet{BKRW993Efficient}. If $\theta_n$ and $\theta_0$ are in $C([a,b])$ for all $n$, then so is $\eta\dot\theta(h)$ as a uniform limit of continuous functions. Assumption \ref{Ass: DGP}(iv) implies that $\{P_{\eta/\sqrt n}\}$ is a sequence of local perturbations such that $\theta_n$ tends to a concave function $\theta_0$. We stress that the parameter of interest in truth is $\theta_n$ while $\theta_0$ is the limit of $\{\theta_n\}$; if $\eta=0$ (e.g., the i.i.d.\ setup with a fixed distribution $P_0$), then $\theta_n=\theta_0$ (and $\eta\dot\theta(h)=0$). Assumption \ref{Ass: estimator}(i) requires an estimator $\hat\theta_n$ of $\theta_n$ that admits a weak limit $\mathbb G$ at rate $\sqrt n$. In turn, Assumption \ref{Ass: estimator}(ii) is a high level condition describing the consequence of degeneracy when $\theta_0$ is strictly concave. For the density problem in the main text, we have $s_n=(n^{-1}\log n)^{2/3}$ as verified by Theorem \ref{Thm: Equivalence, local}. Note that the order in probability (i.e., $O_p$) instead of almost surely (i.e., $O_{a.s.}$) suffices for our inferential purpose.

By a simple generalization of Lemma 3.2 in \citet{BeareandMoon2015}, $\mathcal M:\ell^\infty([a,b])\to\ell^\infty([a,b])$ is Hadamard directionally differentiable at the concave $\theta_0$ tangentially to $C([a,b])$. Like $\mathcal M_F'$ in the main text, the derivative $\mathcal M_{\theta_0}'(h)$ majorizes $h$ by concave functions on regions over which $\theta_0$ is affine but acts like an identity map elsewhere. We omit the explicit expression of the derivative $\mathcal M_{\theta_0}'$ as it is not important to the implementation of the bootstrap, but refer the reader to \citet{BeareandMoon2015}. By Assumptions \ref{Ass: DGP} and \ref{Ass: estimator}(i), we in turn have by analogous arguments as in the proof of Theorem \ref{Thm: test concavity, local} that
\begin{align}
\sqrt n \phi(\hat\theta_n)\convl \|\mathcal M_{\theta_0}'(\mathbb G+\eta\dot\theta(h))-(\mathbb G+\eta\dot\theta(h))\|_p \text{ under }P_n^n~.
\end{align}
Therefore, letting $\eta=0$ yields the pointwise weak limit
\begin{align}\label{Eqn: weak limit}
\|\mathcal M_{\theta_0}'(\mathbb G)-(\mathbb G)\|_p~.
\end{align}
If $\theta_0$ is strictly concave, then $\mathcal M_{\theta_0}'$ is the identity map so that \eqref{Eqn: weak limit} is degenerate at $0$, consistent with Assumption \ref{Ass: estimator}(ii). Following the main text, we differentiate strict concavity from non-strict concavity by introducing a tuning parameter $\kappa_n>0$ such that $\kappa_n\to 0$ and $\sqrt n s_n/\kappa_n\to 0$ as $n\to\infty$.

Towards construction of the critical values, we estimate the law of \eqref{Eqn: weak limit} by the rescaled bootstrap. First, we estimate $\mathcal M_{\theta_0}'$ by: for any $h\in\ell^\infty([a,b])$,
\begin{align}\label{Eqn: rescaled bootsrap3}
\hat{\mathcal M}_n'(h)=\frac{\mathcal M(\hat\theta_n+t_nh)-\mathcal M(\hat\theta_n)}{t_n}~,
\end{align}
where $t_n\downarrow 0$ such that $t_n\sqrt n\to\infty$. Next, we need to bootstrap the law of $\mathbb G$. In order to accommodate flexible bootstrap schemes, we consider a general bootstrap quantity $\hat{\mathbb G}_n: \{X_i,W_{ni}\}_{i=1}^n\to\ell^\infty([a,b])$ where $\{W_{ni}\}_{i=1}^n$ are bootstrap weights independent of the data $\{X_i\}_{i=1}^n$. For example, for the standard empirical bootstrap, $W_n\equiv(W_{n1},\ldots,W_{nn})$ is a multinomial random vector over $n$ categories with probabilities $(1/n,\ldots,1/n)$. More generally, one may consider multiplier or exchangeable bootstrap, corresponding respectively to i.i.d.\ or exchangeable weights \citep{Praestgaard_Wellner1993,Kosorok2008}. To formalize the notion of bootstrap consistency, we employ the bounded Lipschitz metric \citep{Vaart1996}. For a generic normed space $\mathbb D$ equipped with norm $\|\cdot\|_{\mathbb D}$ (e.g., $\mathbb D=\ell^\infty([a,b])$), let
\begin{align*}
\mathrm{BL}_1(\mathbb D)\equiv\{f:\mathbb D\to\mathbf R: \sup_{x\in\mathbb D}|f(x)|\le 1\text{ and }\sup_{x,y\in\mathbb D: x\neq y}\frac{|f(x)-f(y)|}{\|x-y\|_{\mathbb D}}\le 1 \}~.
\end{align*}
If $\hat{\mathbb G}_n:\{X_i,W_{ni}\}_{i=1}^n\to\mathbb D$ is a bootstrap for the law of a random element $\mathbb G$ in $\mathbb D$, then $\hat{\mathbb G}_n$ is said to be a consistent bootstrap for $\mathbb G$ if
\begin{align}
\sup_{f\in\mathrm{BL}_1(\mathbb D)} |E[f(\hat{\mathbb G}_n)|\{X_i\}_{i=1}^n]-E[f(\mathbb G)]|=o_p(1)~.
\end{align}
Finally, for $\alpha\in(0,1)$, we set our critical value $\hat c_{1-\alpha}=\max\{\kappa_n,\hat c_{1-\alpha}^*\}$ where
\begin{align}\label{Eqn: critical value}
\hat c_{1-\alpha}^*\equiv\inf\{c\in\mathbf R: P(\|\hat{\mathcal M}_n'(\hat{\mathbb G}_n)-(\hat{\mathbb G}_n)\|_p\le c|\{X_i\}_{i=1}^n)\ge 1-\alpha\}~.
\end{align}
We then reject the null that $\theta_n$ is concave if $\sqrt n \phi(\hat\theta_n)>\hat c_{1-\alpha}$.

To ensure validity of $\hat c_{1-\alpha}$ as our critical value, we impose:

\begin{ass}\label{Ass: bootstrap}
(i) $\hat{\mathbb G}_n: \{X_i,W_{ni}\}_{i=1}^n\to\ell^\infty([a,b])$ where $\{W_{ni}\}_{i=1}^n$ are bootstrap weights independent of $\{X_i\}_{i=1}^n$; (ii) $\sup_{f\in\mathrm{BL}_1(\mathbb D)} |E[f(\hat{\mathbb G}_n)|\{X_i\}_{i=1}^n]-E[f(\mathbb G)]|=o_p(1)$; (iii) $\hat{\mathbb G}_n$ is asymptotically measurable in $\{X_i,W_{ni}\}_{i=1}^n$; (iv) $f(\hat{\mathbb G}_n)$ is a measurable function of $\{W_{ni}\}_{i=1}^n$ outer almost surely in $\{X_i\}_{i=1}^n$ for any continuous and uniformly bounded $f:\ell^\infty([a,b])\to\mathbf R$.
\end{ass}

\begin{ass}\label{Ass: tuning}
(i) $\{t_n\}$ is a sequence of positive scalars such that $t_n\to 0$ and $\sqrt n t_n \to \infty$ as $n\to\infty$; (ii) $\{\kappa_n\}$ is a sequence of positive scalars such that $\kappa_n\to 0$ and $\sqrt n s_n/\kappa_n\to 0$ as $n\to\infty$.
\end{ass}

\begin{ass}\label{Ass: quantile}
The cdf of $\|\mathcal M_{\theta_0}'(\mathbb G)-\mathbb G\|_p$ is continuous and strictly increasing at its $(1-\alpha)$th quantile $c_{1-\alpha}^*$ when $\theta_0$ is non-strictly concave.
\end{ass}

Assumption \ref{Ass: bootstrap}(i)(ii) formalize the bootstrap consistency of $\hat{\mathbb G}_n$ per our discussions above. Assumption \ref{Ass: bootstrap}(iii)(iv) are mild technical conditions. The precise meaning of Assumption \ref{Ass: bootstrap}(iii) is that, for any $f\in\mathrm{BL}_1(\mathbb D)$,
\[
E^*[f(\hat{\mathbb G}_n)]-E_*[f(\hat{\mathbb G}_n)]\to 0~,
\]
where $E^*$ and $E_*$ are respectively outer and inner expectations with respect to $\{X_i,W_{ni}\}_{i=1}^n$ (jointly)---see Chapter 1.2 in \citet{Vaart1996} for more discussions. Assumption \ref{Ass: bootstrap}(iii) can be verified by appealing to existing results directly; see Theorems 2.6 and 10.4 in \citet{Kosorok2008}, or Theorem 3.6.13 in \citet{Vaart1996} in conjunction with Lemma S.3.9 in \citet{FangSantos2018HDD}. Assumption \ref{Ass: bootstrap}(iv) is met whenever $f(\hat{\mathbb G}_n)$ is continuous in $\{W_{ni}\}_{i=1}^n$, as in common bootstrap schemes. In the main text, Assumption \ref{Ass: bootstrap} is automatically satisfied for $\hat{\mathbb G}_n$ that is constructed by the classical empirical bootstrap---see Theorem 3.6.1 in \citet{Vaart1996}. Assumption \ref{Ass: tuning} imposes suitable rate conditions on the tuning parameter $t_n$ (for the rescaled bootstrap) and on $\kappa_n$ (for the KW-selection). Finally, Assumption \ref{Ass: quantile} is an analog of Assumption \ref{Ass: quantile regularity}, which ensures that bootstrap consistency delivers consistency of the critical value $\hat c_{1-\alpha}^*$.

Given the above assumptions, we may obtain the following analog of Theorem \ref{Thm: test concavity, local} which subsumes the pointwise results in Theorem \ref{Thm: test concavity, pointwise}. As before, let
\[
\pi_n(P_{\eta/\sqrt{n}})\equiv P_n^n(\sqrt{n}\phi(\hat\theta_n)>\hat c_{1-\alpha})~.
\]

\begin{thm}\label{Thm: general}
Let Assumptions \ref{Ass: DGP}, \ref{Ass: estimator}, \ref{Ass: bootstrap}, \ref{Ass: tuning} and \ref{Ass: quantile} hold, and set $\hat c_{1-\alpha}=\max\{\kappa_n,\hat c_{1-\alpha}^*\}$ with $\hat c_{1-\alpha}^*$ defined as in \eqref{Eqn: critical value}. Then it follows that
\begin{enumerate}
\item For any $\eta\in\mathbf R$, (i) if $\theta_0$ is strictly concave, then $\liminf_{n\to\infty}\pi_n(P_{\eta/\sqrt{n}})=0$, and (ii) if $\theta_0$ is non-strictly concave, then
\begin{multline}\label{Thm: test concavity, limit local2}
\liminf_{n\to\infty}\pi_n(P_{\eta/\sqrt{n}})\\ \ge P(\|\mathcal M_{\theta_0}' (\mathbb G+\eta\dot\theta(h))-\{\mathbb G+\eta\dot\theta(h)\}\|_p>c_{1-\alpha}^*)~,
\end{multline}
where \eqref{Thm: test concavity, limit local2} holds with equality if, in addition, the cdf of $\|\mathcal M_{\theta_0}' (\mathbb G+\eta\dot\theta(h))-\{\mathbb G+\eta\dot\theta(h)\}\|_p$ is continuous at $c_{1-\alpha}^*$.
\item For any $\eta\le 0$, we have
\begin{align}\label{Thm: test concavity, local size2}
\limsup_{n\to\infty}\pi_n(P_{\eta/\sqrt{n}})\le\alpha~.
\end{align}
\item If $\|\mathcal M_{\theta_0}' (\eta\dot\theta(h))-\eta\dot\theta(h)\|_p>0$, then
\begin{align}\label{Thm: test concavity, local power2}
\liminf_{\eta\uparrow \infty} \liminf_{n\to\infty}\pi_n(P_{\eta/\sqrt{n}}) = 1~.
\end{align}
\end{enumerate}
\end{thm}

The proof of Theorem \ref{Thm: general} is in complete accord with the proof of Theorem \ref{Thm: test concavity, local} given our high level assumptions, and is thus omitted. In what follows, we therefore focus on verifying some of the main assumptions and in particular Assumption \ref{Ass: estimator}(ii) for Examples \ref{Ex: hazard} and \ref{Ex: regression}. Assumptions \ref{Ass: DGP}(i)(ii)(iv) and \ref{Ass: quantile} may be thought of as regulating the data generating process, while Assumption \ref{Ass: tuning} is simply concerned with choices of tuning parameters.

\begin{exctd}[\ref{Ex: hazard}]
Let $\mathbf P$ be the class of distributions on $\mathbf R_+$ dominated by the Lebesgue measure. By Example 5.3.5 in \citet{BKRW993Efficient}, $\Lambda:\mathbf P\to\ell^\infty([0,b])$ is a regular parameter which may be estimated efficiently by the empirical cumulative hazard function $\hat\Lambda_n$: for any $t\in[0,b]$,
\begin{align}
\hat\Lambda_n(t)\equiv\Lambda_{\mathbb F_n}(t)=\int_0^t \frac{1}{1-\mathbb F_{n,-}(u)}\,\mathrm d\mathbb F_n(u)~,
\end{align}
where $\mathbb F_n$ is the empirical cdf and $\mathbb F_{n,-}$ is the left-continuous version of $\mathbb F_n$. Since the map $F\mapsto\Lambda_F$ is Hadamard differentiable (under regularity conditions)---see, for example, Lemma 20.14 in \citet{Vaart1998}, by the Delta method in conjunction with weak convergence of $\mathbb F_n$, one can show that Assumption \ref{Ass: estimator}(i) is satisfied with $\mathbb G=\mathbb W\circ \chi\in C([a,b])$, where $\mathbb W$ is the standard Brownian motion and $\chi(t)\equiv F(t)/(1-F(t))$ for any $t\in[0,b]$. We refer the reader also to Chapter 6 in \citet{Shorack_Wellner1986empirical} for a detailed treatment of weak convergence of $\hat\Lambda_n$. If $\mathbb F_n^*$ is the bootstrap empirical distribution as described below \eqref{Eqn: rescaled bootsrap}, then the aforementioned Hadamard differentiability implies that $\hat{\mathbb G}_n=\sqrt n \{\Lambda_{\mathbb F_n^*}-\Lambda_{\mathbb F_n}\}$ satisfies Assumption \ref{Ass: bootstrap} by a combination of Theorems 3.6.1 and 3.9.11 in \citet{Vaart1996} and Lemma S.3.9 in \citet{FangSantos2018HDD}. The main condition to verify is now Assumption \ref{Ass: estimator}(ii). In this regard, we believe that, as in Theorem \ref{Thm: Equivalence, local}, it is possible to obtain the same asymptotic order under $\prod_{i=1}^{n}P_{\eta/\sqrt n}$ as under $\prod_{i=1}^{n}P_0$. Since such a development is beyond the scope of this paper, we provide a shortcut building upon existing results under $\prod_{i=1}^{n}P_0$, at the cost of slowing down the rate only a little bit. To illustrate, we first note that existing results---see, for example, \citet{MacgibbonLuYounes2002Failure}, \citet{GroeneboomJongbloed2012Isotonic} and \citet{Durot_Lopuhaa2014KW}---imply under regularity conditions that,
\begin{align}
\|\mathcal M(-\hat\Lambda_n)-(-\hat\Lambda_n)\|_\infty=O_p((n^{-1}\log n)^{2/3}) \text{ under } \prod_{i=1}^{n}P_0~.
\end{align}
Therefore, we must have that, for any $\ell_n\uparrow\infty$ sufficiently slow,
\begin{align}
\|\mathcal M(-\hat\Lambda_n)-(-\hat\Lambda_n)\|_\infty=o_p((n^{-1}\log n)^{2/3}\ell_n) \text{ under } \prod_{i=1}^{n}P_0~.
\end{align}
Since $\prod_{i=1}^{n}P_0$ and $\prod_{i=1}^{n}P_{\eta/\sqrt n}$ are mutually contiguous (see the proof of Theorem \ref{Thm: test concavity, local}), it follows by Le Cam's first lemma that
\begin{align}
\|\mathcal M(-\hat\Lambda_n)-(-\hat\Lambda_n)\|_\infty=o_p((n^{-1}\log n)^{2/3}\ell_n) \text{ under } \prod_{i=1}^{n}P_{\eta/\sqrt n}~.
\end{align}
Thus, Assumption \ref{Ass: estimator}(ii) is satisfied with $s_n=(n^{-1}\log n)^{2/3}\ell_n$ for a suitable $\ell_n\uparrow$. For example, if $\ell_n=(\log n)^{1/3}$, then $s_n=n^{-2/3}\log n$. \qed
\end{exctd}

\begin{exctd}[\ref{Ex: regression}]
In verifying Assumptions \ref{Ass: DGP}(iii) and \ref{Ass: estimator}(i), we focus on the pointwise asymptotics, as a through investigation of regular/efficient estimation of $\Lambda$ as in Chapter 5 in \citet{BKRW993Efficient} is beyond the scope of this paper---and we could not find existing results in this regard. Nonetheless, we provide a sketch here. Let $p_z$ and $p_\epsilon$ be the densities of $Z_1$ and $\epsilon_1$ with respect to the Lebesgue measure. Then the joint density $p_x$ of $X_1\equiv (Y_1,Z_1)'$ at
$x\equiv(y,z)'$ is $p_\epsilon(y-m(z))p_z(z)$. Thus, the differentiable path $\{P_{\eta/\sqrt n}\}$ may be constructed by perturbing the density $p_\epsilon$, which causes no changes in $\Lambda$ in view of \eqref{Eqn: cqf}, or perturbing $m$ and/or $p_z$. Here, by ``perturbing'' we mean a pathwise differentiable sequence as in Definition \ref{Defn: differentiable path} in the case of $p_z$, or, in the case of $m$, a scalar-parametrization $\eta\mapsto m_\eta\in\mathbf H$ for some suitable Hilbert space $\mathbf H$ such that $m_\eta=m+\eta g + o(\eta)$ as $\eta\to 0$ for some $g\in\mathbf H$. A thorough treatment in this regard may be found in \citet{Vaart1991differentibility}. Since pathwise differentiability is Hadamard differentiability---see, for example, \citet[p.456]{BKRW993Efficient}, it follows by Lemmas 3.9.23 and 3.9.27 in \citet{Vaart1996} in conjunction with Lemma 3.9.3 in \citet{Vaart1996} that $\eta\mapsto \Lambda_\eta$ is pathwise differentiable where, for any $t\in[0,1]$,
\begin{align}
\Lambda_\eta(t)\equiv\int_0^t m_\eta (Q_\eta(u))\,\mathrm du~,
\end{align}
with the quantile function $Q_\eta$ arising from perturbing $p_z$. This completes sketching the verification of Assumption \ref{Ass: DGP}(ii).

Next, we verify Assumption \ref{Ass: estimator}(i) under pointwise asymptotics (i.e., $\eta=0$); the general case may be handled with the help of Theorem 3.10.7 and Lemma 3.10.11 in \citet{Vaart1996}, in a manner similar to Theorem 3.10.12 in \citet{Vaart1996}. Let $\{Z_i\}_{i=1}^n$ be arranged in ascending order as $\{Z_{(i)}\}_{i=1}^n$, and let the corresponding $\{Y_i\}_{i=1}^n$ and $\{\epsilon_i\}_{i=1}^n$ be denoted by $\{Y_{(i)}\}_{i=1}^n$ and $\{\epsilon_{(i)}\}_{i=1}^n$. Then the sample cumulative sum diagram is given by $\{(k,S_n(k)):k=0,\ldots,n\}$ where $S_n(0)=0$ and, for $k=1,\ldots,n$,
\begin{equation*}
S_n(k)=\sum_{j=1}^kY_{(j)}~.
\end{equation*}
The cumulative sum diagram plays an important role in isotonic regression because the isotonic regression estimator of $m$ evaluated at $Z_{(i)}$ is precisely the left derivative at $i$ of the greatest convex minorant over $[0,n]$ of $\{Y_{(i)}\}_{i=1}^n$---see \citet{Brunk1958InequalityParameter} and \citet{Mukerjee1988Monotone}. Following \citet{Beare_Fang2016Grenander}, we identify $S_n$ as a random element $\hat\Lambda_n\in\ell^\infty([0,1])$ defined by
\begin{align}
\hat\Lambda_n(t)=\frac{1}{n}\sum_{i=0}^{[nt]}Y_{(i)}+\frac{nt-[nt]}{n}Y_{([nt]+1)}
\end{align}
for all $t\in[0,1]$, where we define $Y_{(0)}=Y_{(n+1)}=0$ and for $x\in\mathbf R$ we denote by $[x]$ the largest integer in $[x-1,x]$. Proposition A.1 in \citet{Beare_Fang2016Grenander} then shows that, under regularity conditions,
\begin{align}\label{Eqn: CSD}
\sqrt n \{\hat\Lambda_n-\Lambda_{P_0}\}\convl \mathbb G\text{ under }\prod_{i=1}^n P_0~,
\end{align}
for some Gaussian process $\mathbb G\in C([0,1])$.

For Assumption \ref{Ass: bootstrap}, we propose the paired bootstrap, i.e., resample from $\{Y_i,Z_i\}_{i=1}^n$ with replacement, and construct $\hat\Lambda_n^*$ in the same fashion as $\hat\Lambda_n$. The main arguments underlying Proposition A.1 in \citet{Beare_Fang2016Grenander} are the Delta method and weak convergence of the partial sum defined by the errors. Thus, the bootstrap consistency of $\hat{\mathbb G}_n\equiv\sqrt n \{\hat\Lambda_n^*-\hat\Lambda_n\}$ follows from Theorem 3.9.4 in \citet{Vaart1996}, and consistency of bootstrapping the partial sum process for which we refer the reader to \citet{Kinateder1992Invariance}, \citet{HolmesKojadinovicQuessy2013Change} and \citet{Calhoun2018Block} for related results. \qed


\end{exctd}

\section{Additional Simulation Results}\label{App: Simulation}

This section collects some simulation results from Section \ref{Sec: simulation}.

{\renewcommand{\arraystretch}{1.5}
\setlength{\tabcolsep}{5pt}
\begin{table}[!ht]
\caption{Comparisons with Existing Tests Based on \eqref{Eqn: cdf compare1}: $n=200$}\label{Tab: comparison, n200}
\begin{footnotesize}
\begin{center}
\begin{tabular}{ccccccccccc}
\hline\hline
 &  \multicolumn{10}{c}{$\lambda$}  \\
\cmidrule{2-11}
 & $-1$ & $-0.1$ & $-0.05$ & $0$ & $0.05$ & $0.1$ & $0.15$ & $0.2$ & $0.25$ & $1$ \\
\hline
R37    & $0.009$ & $0.064$ & $0.086$ & $0.120$ & $0.120$ & $0.168$ & $0.170$ & $0.252$ & $0.283$ & $0.978$ \\
R38    & $0.009$ & $0.066$ & $0.085$ & $0.119$ & $0.123$ & $0.169$ & $0.173$ & $0.254$ & $0.286$ & $0.978$ \\
R3L    & $0.008$ & $0.065$ & $0.085$ & $0.119$ & $0.122$ & $0.166$ & $0.169$ & $0.248$ & $0.284$ & $0.979$ \\
R47    & $0.006$ & $0.049$ & $0.064$ & $0.095$ & $0.105$ & $0.143$ & $0.145$ & $0.213$ & $0.258$ & $0.971$ \\
R48    & $0.006$ & $0.051$ & $0.063$ & $0.096$ & $0.104$ & $0.141$ & $0.145$ & $0.215$ & $0.256$ & $0.969$ \\
R4L    & $0.006$ & $0.050$ & $0.064$ & $0.097$ & $0.104$ & $0.142$ & $0.144$ & $0.213$ & $0.256$ & $0.971$ \\
R57    & $0.003$ & $0.042$ & $0.058$ & $0.082$ & $0.095$ & $0.136$ & $0.132$ & $0.203$ & $0.240$ & $0.968$ \\
R58    & $0.003$ & $0.041$ & $0.059$ & $0.082$ & $0.096$ & $0.136$ & $0.133$ & $0.204$ & $0.238$ & $0.970$ \\
R5L    & $0.003$ & $0.042$ & $0.058$ & $0.081$ & $0.097$ & $0.134$ & $0.134$ & $0.203$ & $0.239$ & $0.968$ \\
R67    & $0.003$ & $0.040$ & $0.055$ & $0.078$ & $0.094$ & $0.126$ & $0.127$ & $0.194$ & $0.233$ & $0.968$ \\
R68    & $0.003$ & $0.040$ & $0.056$ & $0.077$ & $0.093$ & $0.126$ & $0.131$ & $0.196$ & $0.236$ & $0.968$ \\
R6L    & $0.003$ & $0.039$ & $0.057$ & $0.077$ & $0.094$ & $0.128$ & $0.128$ & $0.194$ & $0.233$ & $0.968$ \\
R77    & $0.003$ & $0.038$ & $0.056$ & $0.069$ & $0.092$ & $0.127$ & $0.125$ & $0.187$ & $0.225$ & $0.970$ \\
R78    & $0.003$ & $0.038$ & $0.056$ & $0.072$ & $0.093$ & $0.125$ & $0.125$ & $0.189$ & $0.225$ & $0.966$ \\
R7L    & $0.003$ & $0.038$ & $0.055$ & $0.071$ & $0.091$ & $0.126$ & $0.124$ & $0.189$ & $0.224$ & $0.968$ \\ \hline
P-Test & $0.003$ & $0.023$ & $0.035$ & $0.051$ & $0.071$ & $0.082$ & $0.109$ & $0.152$ & $0.198$ & $0.973$ \\
D-Test & $0.000$ & $0.015$ & $0.029$ & $0.044$ & $0.067$ & $0.110$ & $0.115$ & $0.163$ & $0.220$ & $0.979$ \\
KS-Test& $0.000$ & $0.015$ & $0.022$ & $0.037$ & $0.045$ & $0.067$ & $0.075$ & $0.098$ & $0.127$ & $0.945$ \\
KL-Test& $0.000$ & $0.012$ & $0.025$ & $0.033$ & $0.051$ & $0.081$ & $0.096$ & $0.134$ & $0.187$ & $0.978$ \\
\hline\hline
\end{tabular}
\end{center}
\end{footnotesize}
\end{table}
}

{\renewcommand{\arraystretch}{1.5}
\setlength{\tabcolsep}{5pt}
\begin{table}[!ht]
\caption{Comparisons with Existing Tests Based on \eqref{Eqn: cdf compare1}: $n=300$}\label{Tab: comparison, n300}
\begin{footnotesize}
\begin{center}
\begin{tabular}{ccccccccccc}
\hline\hline
  & \multicolumn{10}{c}{$\lambda$}  \\
\cmidrule{2-11}
 & $-1$ & $-0.1$ & $-0.05$ & $0$ & $0.05$ & $0.1$ & $0.15$ & $0.2$ & $0.25$ & $1$ \\
\hline
R37    & $0.009$ & $0.059$ & $0.071$ & $0.111$ & $0.147$ & $0.163$ & $0.220$ & $0.275$ & $0.366$ & $0.999$ \\
R38    & $0.009$ & $0.059$ & $0.069$ & $0.111$ & $0.148$ & $0.160$ & $0.220$ & $0.275$ & $0.364$ & $0.999$ \\
R3L    & $0.009$ & $0.059$ & $0.071$ & $0.110$ & $0.147$ & $0.166$ & $0.217$ & $0.274$ & $0.360$ & $0.999$ \\
R47    & $0.002$ & $0.038$ & $0.053$ & $0.084$ & $0.119$ & $0.136$ & $0.187$ & $0.238$ & $0.330$ & $0.999$ \\
R48    & $0.002$ & $0.040$ & $0.053$ & $0.082$ & $0.120$ & $0.135$ & $0.187$ & $0.242$ & $0.332$ & $0.999$ \\
R4L    & $0.002$ & $0.039$ & $0.053$ & $0.084$ & $0.120$ & $0.135$ & $0.187$ & $0.239$ & $0.327$ & $0.999$ \\
R57    & $0.002$ & $0.031$ & $0.045$ & $0.077$ & $0.110$ & $0.116$ & $0.178$ & $0.220$ & $0.314$ & $0.999$ \\
R58    & $0.002$ & $0.032$ & $0.045$ & $0.077$ & $0.110$ & $0.115$ & $0.175$ & $0.217$ & $0.315$ & $0.999$ \\
R5L    & $0.002$ & $0.032$ & $0.045$ & $0.077$ & $0.109$ & $0.116$ & $0.175$ & $0.216$ & $0.316$ & $0.999$ \\
R67    & $0.002$ & $0.027$ & $0.041$ & $0.070$ & $0.102$ & $0.113$ & $0.169$ & $0.220$ & $0.309$ & $0.998$ \\
R68    & $0.003$ & $0.028$ & $0.040$ & $0.072$ & $0.103$ & $0.112$ & $0.171$ & $0.218$ & $0.308$ & $0.998$ \\
R6L    & $0.002$ & $0.028$ & $0.040$ & $0.071$ & $0.103$ & $0.113$ & $0.169$ & $0.219$ & $0.311$ & $0.998$ \\
R77    & $0.003$ & $0.027$ & $0.036$ & $0.071$ & $0.104$ & $0.109$ & $0.166$ & $0.217$ & $0.305$ & $0.998$ \\
R78    & $0.003$ & $0.027$ & $0.038$ & $0.071$ & $0.103$ & $0.106$ & $0.165$ & $0.220$ & $0.307$ & $0.998$ \\
R7L    & $0.003$ & $0.027$ & $0.035$ & $0.071$ & $0.104$ & $0.112$ & $0.162$ & $0.217$ & $0.306$ & $0.998$ \\ \hline
P-Test & $0.000$ & $0.027$ & $0.034$ & $0.054$ & $0.084$ & $0.096$ & $0.128$ & $0.190$ & $0.265$ & $0.997$ \\
D-Test & $0.000$ & $0.016$ & $0.017$ & $0.051$ & $0.083$ & $0.101$ & $0.155$ & $0.218$ & $0.320$ & $0.999$ \\
KS-Test& $0.000$ & $0.011$ & $0.015$ & $0.036$ & $0.048$ & $0.056$ & $0.100$ & $0.144$ & $0.201$ & $0.988$ \\
KL-Test& $0.000$ & $0.012$ & $0.013$ & $0.041$ & $0.059$ & $0.070$ & $0.128$ & $0.173$ & $0.269$ & $1.000$ \\
\hline\hline
\end{tabular}
\end{center}
\end{footnotesize}
\end{table}
}

{\renewcommand{\arraystretch}{1.5}
\setlength{\tabcolsep}{5pt}
\begin{table}[!ht]
\caption{Comparisons with Existing Tests Based on \eqref{Eqn: cdf compare1}: $n=600$}\label{Tab: comparison, n600}
\begin{footnotesize}
\begin{center}
\begin{tabular}{ccccccccccc}
\hline\hline
  & \multicolumn{10}{c}{$\lambda$}  \\
\cmidrule{2-11}
 & $-1$ & $-0.1$ & $-0.05$ & $0$ & $0.05$ & $0.1$ & $0.15$ & $0.2$ & $0.25$ & $1$ \\
\hline
R37  & $0.006$ & $0.045$ & $0.059$ & $0.097$ & $0.153$ & $0.218$ & $0.274$ & $0.416$ & $0.496$ & $1.000$ \\
R38  & $0.006$ & $0.045$ & $0.059$ & $0.097$ & $0.157$ & $0.216$ & $0.274$ & $0.416$ & $0.496$ & $1.000$ \\
R3L  & $0.005$ & $0.043$ & $0.060$ & $0.098$ & $0.155$ & $0.217$ & $0.270$ & $0.416$ & $0.495$ & $1.000$ \\
R47  & $0.000$ & $0.029$ & $0.041$ & $0.080$ & $0.125$ & $0.183$ & $0.238$ & $0.375$ & $0.459$ & $1.000$ \\
R48  & $0.000$ & $0.028$ & $0.042$ & $0.079$ & $0.124$ & $0.185$ & $0.237$ & $0.376$ & $0.457$ & $1.000$ \\
R4L  & $0.000$ & $0.030$ & $0.042$ & $0.080$ & $0.123$ & $0.184$ & $0.234$ & $0.380$ & $0.457$ & $1.000$ \\
R57  & $0.000$ & $0.023$ & $0.032$ & $0.072$ & $0.116$ & $0.168$ & $0.228$ & $0.355$ & $0.439$ & $1.000$ \\
R58  & $0.000$ & $0.026$ & $0.032$ & $0.073$ & $0.116$ & $0.167$ & $0.227$ & $0.355$ & $0.438$ & $1.000$ \\
R5L  & $0.000$ & $0.024$ & $0.032$ & $0.072$ & $0.115$ & $0.168$ & $0.226$ & $0.356$ & $0.441$ & $1.000$ \\
R67  & $0.000$ & $0.023$ & $0.032$ & $0.068$ & $0.110$ & $0.162$ & $0.218$ & $0.347$ & $0.429$ & $1.000$ \\
R68  & $0.000$ & $0.021$ & $0.030$ & $0.068$ & $0.112$ & $0.161$ & $0.221$ & $0.350$ & $0.428$ & $1.000$ \\
R6L  & $0.000$ & $0.021$ & $0.029$ & $0.071$ & $0.112$ & $0.159$ & $0.219$ & $0.348$ & $0.432$ & $1.000$ \\
R77  & $0.000$ & $0.019$ & $0.028$ & $0.067$ & $0.110$ & $0.156$ & $0.216$ & $0.349$ & $0.424$ & $1.000$ \\
R78  & $0.000$ & $0.019$ & $0.028$ & $0.068$ & $0.110$ & $0.155$ & $0.216$ & $0.343$ & $0.425$ & $1.000$ \\
R7L  & $0.000$ & $0.020$ & $0.029$ & $0.067$ & $0.105$ & $0.157$ & $0.217$ & $0.344$ & $0.426$ & $1.000$ \\
 \hline
P-Test  & $0.003$ & $0.022$ & $0.035$ & $0.049$ & $0.105$ & $0.134$ & $0.202$ & $0.310$ & $0.435$ & $1.000$ \\
D-Test  & $0.000$ & $0.013$ & $0.015$ & $0.050$ & $0.107$ & $0.165$ & $0.234$ & $0.387$ & $0.469$ & $1.000$ \\
KS-Test & $0.000$ & $0.012$ & $0.016$ & $0.041$ & $0.063$ & $0.109$ & $0.157$ & $0.269$ & $0.343$ & $1.000$ \\
KL-Test & $0.000$ & $0.009$ & $0.012$ & $0.039$ & $0.086$ & $0.148$ & $0.200$ & $0.345$ & $0.447$ & $1.000$ \\
\hline\hline
\end{tabular}
\end{center}
\end{footnotesize}
\end{table}
}

{\renewcommand{\arraystretch}{1.5}
\setlength{\tabcolsep}{5pt}
\begin{table}[!ht]
\caption{Comparisons with Existing Tests Based on \eqref{Eqn: cdf compare1}: $n=800$}\label{Tab: comparison, n800}
\begin{footnotesize}
\begin{center}
\begin{tabular}{ccccccccccc}
\hline\hline
  & \multicolumn{10}{c}{$\lambda$}  \\
\cmidrule{2-11}
 & $-1$ & $-0.1$ & $-0.05$ & $0$ & $0.05$ & $0.1$ & $0.15$ & $0.2$ & $0.25$ & $1$ \\
\hline
R37  & $0.004$ & $0.039$ & $0.060$ & $0.100$ & $0.153$ & $0.234$ & $0.335$ & $0.492$ & $0.592$ & $1.000$ \\
R38  & $0.004$ & $0.037$ & $0.060$ & $0.101$ & $0.150$ & $0.234$ & $0.331$ & $0.489$ & $0.593$ & $1.000$ \\
R3L  & $0.004$ & $0.040$ & $0.060$ & $0.101$ & $0.152$ & $0.234$ & $0.335$ & $0.487$ & $0.592$ & $1.000$ \\
R47  & $0.002$ & $0.027$ & $0.050$ & $0.078$ & $0.117$ & $0.206$ & $0.299$ & $0.447$ & $0.550$ & $1.000$ \\
R48  & $0.003$ & $0.027$ & $0.049$ & $0.077$ & $0.115$ & $0.207$ & $0.296$ & $0.445$ & $0.549$ & $1.000$ \\
R4L  & $0.000$ & $0.027$ & $0.049$ & $0.078$ & $0.117$ & $0.205$ & $0.295$ & $0.449$ & $0.550$ & $1.000$ \\
R57  & $0.000$ & $0.021$ & $0.044$ & $0.069$ & $0.104$ & $0.197$ & $0.279$ & $0.435$ & $0.536$ & $1.000$ \\
R58  & $0.000$ & $0.022$ & $0.044$ & $0.069$ & $0.104$ & $0.196$ & $0.278$ & $0.437$ & $0.538$ & $1.000$ \\
R5L  & $0.000$ & $0.020$ & $0.045$ & $0.070$ & $0.102$ & $0.194$ & $0.281$ & $0.435$ & $0.536$ & $1.000$ \\
R67  & $0.000$ & $0.021$ & $0.044$ & $0.066$ & $0.095$ & $0.189$ & $0.272$ & $0.421$ & $0.525$ & $1.000$ \\
R68  & $0.000$ & $0.020$ & $0.044$ & $0.064$ & $0.095$ & $0.189$ & $0.275$ & $0.419$ & $0.529$ & $1.000$ \\
R6L  & $0.000$ & $0.021$ & $0.044$ & $0.067$ & $0.098$ & $0.190$ & $0.271$ & $0.421$ & $0.528$ & $1.000$ \\
R77  & $0.000$ & $0.020$ & $0.042$ & $0.060$ & $0.093$ & $0.183$ & $0.271$ & $0.414$ & $0.524$ & $1.000$ \\
R78  & $0.000$ & $0.021$ & $0.042$ & $0.061$ & $0.097$ & $0.182$ & $0.272$ & $0.415$ & $0.525$ & $1.000$ \\
R7L  & $0.000$ & $0.021$ & $0.042$ & $0.061$ & $0.096$ & $0.180$ & $0.269$ & $0.415$ & $0.524$ & $1.000$ \\
 \hline
P-Test & $0.003$ & $0.021$ & $0.040$ & $0.055$ & $0.113$ & $0.179$ & $0.265$ & $0.391$ & $0.514$ & $1.000$ \\
D-Test & $0.000$ & $0.009$ & $0.022$ & $0.045$ & $0.099$ & $0.191$ & $0.319$ & $0.466$ & $0.600$ & $1.000$ \\
KS-Test & $0.000$ & $0.006$ & $0.025$ & $0.042$ & $0.061$ & $0.131$ & $0.214$ & $0.340$ & $0.446$ & $1.000$ \\
KL-Test & $0.000$ & $0.007$ & $0.017$ & $0.036$ & $0.087$ & $0.178$ & $0.263$ & $0.428$ & $0.565$ & $1.000$ \\
\hline\hline
\end{tabular}
\end{center}
\end{footnotesize}
\end{table}
}

{\renewcommand{\arraystretch}{1.5}
\setlength{\tabcolsep}{5pt}
\begin{table}[!h]
\caption{Comparisons with Existing Tests Based on \eqref{Eqn: cdf compare1}: $n=1000$}\label{Tab: comparison, n1000}
\begin{footnotesize}
\begin{center}
\begin{tabular}{cccccccccccc}
\hline\hline
  & \multicolumn{10}{c}{$\lambda$}  \\
\cmidrule{2-11}
 & $-1$ & $-0.1$ & $-0.05$ & $0$ & $0.05$ & $0.1$ & $0.15$ & $0.2$ & $0.25$ & $1$ \\
\hline
R37  & $0.003$ & $0.029$ & $0.069$ & $0.096$ & $0.163$ & $0.230$ & $0.395$ & $0.517$ & $0.658$ & $1.000$ \\
R38  & $0.003$ & $0.030$ & $0.070$ & $0.095$ & $0.162$ & $0.231$ & $0.396$ & $0.518$ & $0.656$ & $1.000$ \\
R3L  & $0.004$ & $0.040$ & $0.060$ & $0.101$ & $0.152$ & $0.234$ & $0.335$ & $0.487$ & $0.592$ & $1.000$ \\
R47  & $0.000$ & $0.019$ & $0.054$ & $0.074$ & $0.127$ & $0.196$ & $0.349$ & $0.480$ & $0.627$ & $1.000$ \\
R48  & $0.000$ & $0.019$ & $0.053$ & $0.074$ & $0.125$ & $0.194$ & $0.349$ & $0.481$ & $0.627$ & $1.000$ \\
R4L  & $0.000$ & $0.019$ & $0.052$ & $0.074$ & $0.126$ & $0.195$ & $0.347$ & $0.479$ & $0.629$ & $1.000$ \\
R57  & $0.000$ & $0.017$ & $0.049$ & $0.068$ & $0.116$ & $0.187$ & $0.332$ & $0.470$ & $0.616$ & $1.000$ \\
R58  & $0.000$ & $0.017$ & $0.049$ & $0.068$ & $0.115$ & $0.187$ & $0.336$ & $0.471$ & $0.613$ & $1.000$ \\
R5L  & $0.000$ & $0.017$ & $0.049$ & $0.068$ & $0.115$ & $0.184$ & $0.332$ & $0.470$ & $0.613$ & $1.000$ \\
R67  & $0.000$ & $0.017$ & $0.047$ & $0.063$ & $0.109$ & $0.174$ & $0.328$ & $0.464$ & $0.609$ & $1.000$ \\
R68  & $0.000$ & $0.017$ & $0.046$ & $0.063$ & $0.108$ & $0.174$ & $0.327$ & $0.463$ & $0.608$ & $1.000$ \\
R6L  & $0.000$ & $0.017$ & $0.048$ & $0.062$ & $0.107$ & $0.172$ & $0.325$ & $0.459$ & $0.606$ & $1.000$ \\
R77  & $0.000$ & $0.017$ & $0.048$ & $0.061$ & $0.103$ & $0.171$ & $0.324$ & $0.454$ & $0.603$ & $1.000$ \\
R78  & $0.000$ & $0.017$ & $0.046$ & $0.061$ & $0.103$ & $0.172$ & $0.324$ & $0.457$ & $0.603$ & $1.000$ \\
R7L  & $0.000$ & $0.017$ & $0.047$ & $0.062$ & $0.104$ & $0.172$ & $0.324$ & $0.454$ & $0.605$ & $1.000$ \\
 \hline
P-Test & $0.002$ & $0.021$ & $0.038$ & $0.060$ & $0.112$ & $0.178$ & $0.302$ & $0.439$ & $0.604$ & $1.000$ \\
D-Test & $0.000$ & $0.004$ & $0.024$ & $0.056$ & $0.104$ & $0.187$ & $0.360$ & $0.509$ & $0.674$ & $1.000$ \\
KS-Test& $0.000$ & $0.009$ & $0.029$ & $0.043$ & $0.062$ & $0.123$ & $0.241$ & $0.380$ & $0.531$ & $1.000$ \\
KL-Test& $0.000$ & $0.006$ & $0.017$ & $0.043$ & $0.089$ & $0.167$ & $0.326$ & $0.487$ & $0.659$ & $1.000$ \\
\hline\hline
\end{tabular}
\end{center}
\end{footnotesize}
\end{table}
}

{\renewcommand{\arraystretch}{1.5}
\setlength{\tabcolsep}{5pt}
\begin{table}[!h]
\caption{Comparisons with Existing Tests Based on \eqref{Eqn: cdf compare2}: $n=200$}\label{Tab: comparison2, n200}
\begin{footnotesize}
\begin{center}
\begin{tabular}{cccccccccc}
\hline\hline
 &  \multicolumn{9}{c}{$\lambda$}  \\
\cmidrule{2-10}
 & $0.1$ & $0.2$ & $0.3$ & $0.4$ & $0.5$ & $0.6$ & $0.7$ & $0.8$ & $0.9$ \\
\hline
R37    & $0.142$ & $0.347$ & $0.651$ & $0.861$ & $0.944$ & $0.974$ & $0.961$ & $0.827$ & $0.214$  \\
R38    & $0.141$ & $0.346$ & $0.650$ & $0.861$ & $0.945$ & $0.974$ & $0.961$ & $0.826$ & $0.214$     \\
R3L    & $0.142$ & $0.342$ & $0.650$ & $0.860$ & $0.945$ & $0.974$ & $0.962$ & $0.825$ & $0.218$  \\
R47    & $0.114$ & $0.280$ & $0.582$ & $0.813$ & $0.906$ & $0.939$ & $0.915$ & $0.705$ & $0.128$     \\
R48    & $0.117$ & $0.280$ & $0.581$ & $0.811$ & $0.907$ & $0.939$ & $0.916$ & $0.709$ & $0.129$     \\
R4L    & $0.117$ & $0.279$ & $0.583$ & $0.811$ & $0.907$ & $0.941$ & $0.912$ & $0.704$ & $0.128$  \\
R57    & $0.109$ & $0.258$ & $0.542$ & $0.787$ & $0.881$ & $0.921$ & $0.881$ & $0.625$ & $0.093$     \\
R58    & $0.105$ & $0.260$ & $0.541$ & $0.779$ & $0.881$ & $0.923$ & $0.882$ & $0.624$ & $0.094$     \\
R5L    & $0.105$ & $0.259$ & $0.544$ & $0.779$ & $0.883$ & $0.921$ & $0.881$ & $0.627$ & $0.093$  \\
R67    & $0.100$ & $0.243$ & $0.511$ & $0.761$ & $0.861$ & $0.906$ & $0.846$ & $0.572$ & $0.075$     \\
R68    & $0.098$ & $0.243$ & $0.509$ & $0.759$ & $0.863$ & $0.908$ & $0.845$ & $0.577$ & $0.075$     \\
R6L    & $0.097$ & $0.246$ & $0.509$ & $0.759$ & $0.864$ & $0.907$ & $0.847$ & $0.575$ & $0.075$     \\
R77    & $0.095$ & $0.231$ & $0.494$ & $0.751$ & $0.844$ & $0.897$ & $0.834$ & $0.544$ & $0.066$     \\
R78    & $0.092$ & $0.232$ & $0.494$ & $0.753$ & $0.846$ & $0.895$ & $0.830$ & $0.546$ & $0.068$     \\
R7L    & $0.092$ & $0.234$ & $0.494$ & $0.752$ & $0.846$ & $0.896$ & $0.832$ & $0.545$ & $0.067$     \\ \hline
P-Test & $0.018$ & $0.016$ & $0.013$ & $0.009$ & $0.012$ & $0.005$ & $0.006$ & $0.004$ & $0.000$  \\
D-Test & $0.038$ & $0.053$ & $0.046$ & $0.044$ & $0.024$ & $0.019$ & $0.009$ & $0.008$ & $0.000$  \\
KS-Test& $0.038$ & $0.109$ & $0.255$ & $0.450$ & $0.522$ & $0.533$ & $0.398$ & $0.148$ & $0.005$ \\
KL-Test& $0.062$ & $0.093$ & $0.116$ & $0.138$ & $0.057$ & $0.014$ & $0.003$ & $0.000$ & $0.000$ \\
\hline\hline
\end{tabular}
\end{center}
\end{footnotesize}
\end{table}
}

{\renewcommand{\arraystretch}{1.5}
\setlength{\tabcolsep}{5pt}
\begin{table}[!h]
\caption{Comparisons with Existing Tests Based on \eqref{Eqn: cdf compare2}: $n=300$}\label{Tab: comparison2, n300}
\begin{footnotesize}
\begin{center}
\begin{tabular}{cccccccccc}
\hline\hline
  & \multicolumn{9}{c}{$\lambda$}  \\
\cmidrule{2-10}
 & $0.1$ & $0.2$ & $0.3$ & $0.4$ & $0.5$ & $0.6$ & $0.7$ & $0.8$ & $0.9$ \\
\hline
R37    & $0.166$ & $0.454$ & $0.851$ & $0.974$ & $0.997$ & $0.998$ & $0.997$ & $0.979$ & $0.487$  \\
R38    & $0.164$ & $0.460$ & $0.852$ & $0.975$ & $0.997$ & $0.998$ & $0.996$ & $0.980$ & $0.483$     \\
R3L    & $0.163$ & $0.454$ & $0.851$ & $0.974$ & $0.997$ & $0.998$ & $0.997$ & $0.979$ & $0.486$     \\
R47    & $0.137$ & $0.391$ & $0.782$ & $0.953$ & $0.992$ & $0.992$ & $0.993$ & $0.939$ & $0.342$     \\
R48    & $0.133$ & $0.387$ & $0.781$ & $0.952$ & $0.992$ & $0.993$ & $0.992$ & $0.940$ & $0.345$     \\
R4L    & $0.136$ & $0.391$ & $0.784$ & $0.952$ & $0.993$ & $0.992$ & $0.993$ & $0.940$ & $0.343$  \\
R57    & $0.119$ & $0.363$ & $0.749$ & $0.932$ & $0.986$ & $0.990$ & $0.988$ & $0.919$ & $0.272$     \\
R58    & $0.120$ & $0.364$ & $0.750$ & $0.932$ & $0.986$ & $0.990$ & $0.987$ & $0.919$ & $0.270$     \\
R5L    & $0.119$ & $0.365$ & $0.750$ & $0.932$ & $0.986$ & $0.990$ & $0.985$ & $0.918$ & $0.273$     \\
R67    & $0.112$ & $0.339$ & $0.727$ & $0.921$ & $0.983$ & $0.988$ & $0.978$ & $0.893$ & $0.221$     \\
R68    & $0.111$ & $0.337$ & $0.725$ & $0.920$ & $0.981$ & $0.988$ & $0.978$ & $0.894$ & $0.219$     \\
R6L    & $0.114$ & $0.342$ & $0.724$ & $0.920$ & $0.982$ & $0.987$ & $0.978$ & $0.892$ & $0.220$     \\
R77    & $0.109$ & $0.335$ & $0.706$ & $0.909$ & $0.978$ & $0.986$ & $0.977$ & $0.869$ & $0.189$     \\
R78    & $0.107$ & $0.332$ & $0.709$ & $0.910$ & $0.979$ & $0.985$ & $0.976$ & $0.869$ & $0.191$     \\
R7L    & $0.108$ & $0.336$ & $0.706$ & $0.908$ & $0.978$ & $0.986$ & $0.976$ & $0.867$ & $0.190$     \\ \hline
P-Test & $0.023$ & $0.020$ & $0.014$ & $0.012$ & $0.008$ & $0.009$ & $0.005$ & $0.003$ & $0.004$  \\
D-Test & $0.037$ & $0.040$ & $0.052$ & $0.032$ & $0.023$ & $0.023$ & $0.011$ & $0.003$ & $0.000$  \\
KS-Test& $0.057$ & $0.174$ & $0.485$ & $0.688$ & $0.844$ & $0.861$ & $0.754$ & $0.434$ & $0.019$ \\
KL-Test& $0.070$ & $0.150$ & $0.239$ & $0.218$ & $0.134$ & $0.037$ & $0.003$ & $0.000$ & $0.000$ \\
\hline\hline
\end{tabular}
\end{center}
\end{footnotesize}
\end{table}
}

{\renewcommand{\arraystretch}{1.5}
\setlength{\tabcolsep}{5pt}
\begin{table}[!h]
\caption{Comparisons with Existing Tests Based on \eqref{Eqn: cdf compare2}: $n=600$}\label{Tab: comparison2, n600}
\begin{footnotesize}
\begin{center}
\begin{tabular}{cccccccccc}
\hline\hline
 &  \multicolumn{9}{c}{$\lambda$}  \\
\cmidrule{2-10}
 & $0.1$ & $0.2$ & $0.3$ & $0.4$ & $0.5$ & $0.6$ & $0.7$ & $0.8$ & $0.9$ \\
\hline
R37    & $0.259$ & $0.827$ & $0.991$ & $1.000$ & $1.000$ & $1.000$ & $1.000$ & $1.000$ & $0.964$  \\
R38    & $0.260$ & $0.830$ & $0.991$ & $1.000$ & $1.000$ & $1.000$ & $1.000$ & $1.000$ & $0.964$  \\
R3L    & $0.260$ & $0.826$ & $0.991$ & $1.000$ & $1.000$ & $1.000$ & $1.000$ & $1.000$ & $0.964$  \\
R47    & $0.204$ & $0.758$ & $0.982$ & $1.000$ & $1.000$ & $1.000$ & $1.000$ & $1.000$ & $0.901$     \\
R48    & $0.204$ & $0.760$ & $0.983$ & $1.000$ & $1.000$ & $1.000$ & $1.000$ & $1.000$ & $0.904$     \\
R4L    & $0.204$ & $0.758$ & $0.982$ & $1.000$ & $1.000$ & $1.000$ & $1.000$ & $1.000$ & $0.904$     \\
R57    & $0.188$ & $0.724$ & $0.977$ & $1.000$ & $1.000$ & $1.000$ & $1.000$ & $1.000$ & $0.835$     \\
R58    & $0.190$ & $0.722$ & $0.977$ & $1.000$ & $1.000$ & $1.000$ & $1.000$ & $1.000$ & $0.832$     \\
R5L    & $0.189$ & $0.724$ & $0.978$ & $1.000$ & $1.000$ & $1.000$ & $1.000$ & $1.000$ & $0.839$     \\
R67    & $0.180$ & $0.703$ & $0.974$ & $1.000$ & $1.000$ & $1.000$ & $1.000$ & $1.000$ & $0.789$     \\
R68    & $0.179$ & $0.704$ & $0.973$ & $1.000$ & $1.000$ & $1.000$ & $1.000$ & $1.000$ & $0.787$     \\
R6L    & $0.178$ & $0.702$ & $0.973$ & $1.000$ & $1.000$ & $1.000$ & $1.000$ & $1.000$ & $0.789$     \\
R77    & $0.175$ & $0.694$ & $0.969$ & $1.000$ & $1.000$ & $1.000$ & $1.000$ & $1.000$ & $0.740$     \\
R78    & $0.173$ & $0.699$ & $0.969$ & $1.000$ & $1.000$ & $1.000$ & $1.000$ & $1.000$ & $0.741$     \\
R7L    & $0.174$ & $0.694$ & $0.969$ & $1.000$ & $1.000$ & $1.000$ & $1.000$ & $1.000$ & $0.742$     \\ \hline
P-Test & $0.025$ & $0.016$ & $0.009$ & $0.007$ & $0.008$ & $0.007$ & $0.007$ & $0.003$ & $0.003$  \\
D-Test & $0.052$ & $0.052$ & $0.053$ & $0.041$ & $0.028$ & $0.024$ & $0.010$ & $0.004$ & $0.000$  \\
KS-Test& $0.109$ & $0.526$ & $0.923$ & $0.997$ & $0.999$ & $1.000$ & $1.000$ & $0.980$ & $0.249$ \\
KL-Test& $0.140$ & $0.377$ & $0.560$ & $0.631$ & $0.604$ & $0.307$ & $0.015$ & $0.000$ & $0.000$ \\
\hline\hline
\end{tabular}
\end{center}
\end{footnotesize}
\end{table}
}

{\renewcommand{\arraystretch}{1.5}
\setlength{\tabcolsep}{5pt}
\begin{table}[!h]
\caption{Comparisons with Existing Tests Based on \eqref{Eqn: cdf compare2}: $n=800$}\label{Tab: comparison2, n800}
\begin{footnotesize}
\begin{center}
\begin{tabular}{cccccccccc}
\hline\hline
 &  \multicolumn{9}{c}{$\lambda$}  \\
\cmidrule{2-10}
 & $0.1$ & $0.2$ & $0.3$ & $0.4$ & $0.5$ & $0.6$ & $0.7$ & $0.8$ & $0.9$ \\
\hline
R37    & $0.311$ & $0.936$ & $1.000$ & $1.000$ & $1.000$ & $1.000$ & $1.000$ & $1.000$ & $0.995$  \\
R38    & $0.313$ & $0.934$ & $1.000$ & $1.000$ & $1.000$ & $1.000$ & $1.000$ & $1.000$ & $0.995$  \\
R3L    & $0.313$ & $0.934$ & $1.000$ & $1.000$ & $1.000$ & $1.000$ & $1.000$ & $1.000$ & $0.995$  \\
R47    & $0.255$ & $0.905$ & $1.000$ & $1.000$ & $1.000$ & $1.000$ & $1.000$ & $1.000$ & $0.986$  \\
R48    & $0.256$ & $0.908$ & $1.000$ & $1.000$ & $1.000$ & $1.000$ & $1.000$ & $1.000$ & $0.985$     \\
R4L    & $0.254$ & $0.907$ & $1.000$ & $1.000$ & $1.000$ & $1.000$ & $1.000$ & $1.000$ & $0.984$     \\
R57    & $0.238$ & $0.878$ & $1.000$ & $1.000$ & $1.000$ & $1.000$ & $1.000$ & $1.000$ & $0.964$     \\
R58    & $0.238$ & $0.880$ & $1.000$ & $1.000$ & $1.000$ & $1.000$ & $1.000$ & $1.000$ & $0.966$     \\
R5L    & $0.238$ & $0.879$ & $1.000$ & $1.000$ & $1.000$ & $1.000$ & $1.000$ & $1.000$ & $0.965$     \\
R67    & $0.233$ & $0.864$ & $1.000$ & $1.000$ & $1.000$ & $1.000$ & $1.000$ & $1.000$ & $0.949$     \\
R68    & $0.232$ & $0.864$ & $1.000$ & $1.000$ & $1.000$ & $1.000$ & $1.000$ & $1.000$ & $0.948$     \\
R6L    & $0.231$ & $0.864$ & $1.000$ & $1.000$ & $1.000$ & $1.000$ & $1.000$ & $1.000$ & $0.949$     \\
R77    & $0.225$ & $0.858$ & $1.000$ & $1.000$ & $1.000$ & $1.000$ & $1.000$ & $1.000$ & $0.934$     \\
R78    & $0.224$ & $0.858$ & $1.000$ & $1.000$ & $1.000$ & $1.000$ & $1.000$ & $1.000$ & $0.936$     \\
R7L    & $0.224$ & $0.856$ & $1.000$ & $1.000$ & $1.000$ & $1.000$ & $1.000$ & $1.000$ & $0.936$     \\ \hline
P-Test & $0.025$ & $0.021$ & $0.019$ & $0.013$ & $0.018$ & $0.010$ & $0.002$ & $0.003$ & $0.003$  \\
D-Test & $0.045$ & $0.045$ & $0.047$ & $0.035$ & $0.040$ & $0.022$ & $0.008$ & $0.002$ & $0.000$  \\
KS-Test& $0.143$ & $0.734$ & $0.996$ & $1.000$ & $1.000$ & $1.000$ & $1.000$ & $0.999$ & $0.581$ \\
KL-Test& $0.190$ & $0.512$ & $0.753$ & $0.837$ & $0.820$ & $0.531$ & $0.060$ & $0.000$ & $0.000$ \\
\hline\hline
\end{tabular}
\end{center}
\end{footnotesize}
\end{table}
}

{\renewcommand{\arraystretch}{1.5}
\setlength{\tabcolsep}{5pt}
\begin{table}[!h]
\caption{Comparisons with Existing Tests Based on \eqref{Eqn: cdf compare2}: $n=1000$}\label{Tab: comparison2, n1000}
\begin{footnotesize}
\begin{center}
\begin{tabular}{cccccccccc}
\hline\hline
  & \multicolumn{9}{c}{$\lambda$}  \\
\cmidrule{2-10}
 & $0.1$ & $0.2$ & $0.3$ & $0.4$ & $0.5$ & $0.6$ & $0.7$ & $0.8$ & $0.9$ \\
\hline
R37    & $0.362$ & $0.985$ & $1.000$ & $1.000$ & $1.000$ & $1.000$ & $1.000$ & $1.000$ & $1.000$  \\
R38    & $0.360$ & $0.985$ & $1.000$ & $1.000$ & $1.000$ & $1.000$ & $1.000$ & $1.000$ & $1.000$  \\
R3L    & $0.360$ & $0.983$ & $1.000$ & $1.000$ & $1.000$ & $1.000$ & $1.000$ & $1.000$ & $1.000$  \\
R47    & $0.309$ & $0.972$ & $1.000$ & $1.000$ & $1.000$ & $1.000$ & $1.000$ & $1.000$ & $0.999$  \\
R48    & $0.309$ & $0.974$ & $1.000$ & $1.000$ & $1.000$ & $1.000$ & $1.000$ & $1.000$ & $0.999$  \\
R4L    & $0.310$ & $0.973$ & $1.000$ & $1.000$ & $1.000$ & $1.000$ & $1.000$ & $1.000$ & $0.999$  \\
R57    & $0.280$ & $0.963$ & $1.000$ & $1.000$ & $1.000$ & $1.000$ & $1.000$ & $1.000$ & $0.997$  \\
R58    & $0.284$ & $0.962$ & $1.000$ & $1.000$ & $1.000$ & $1.000$ & $1.000$ & $1.000$ & $0.996$  \\
R5L    & $0.284$ & $0.964$ & $1.000$ & $1.000$ & $1.000$ & $1.000$ & $1.000$ & $1.000$ & $0.996$  \\
R67    & $0.270$ & $0.958$ & $1.000$ & $1.000$ & $1.000$ & $1.000$ & $1.000$ & $1.000$ & $0.993$  \\
R68    & $0.269$ & $0.955$ & $1.000$ & $1.000$ & $1.000$ & $1.000$ & $1.000$ & $1.000$ & $0.993$  \\
R6L    & $0.270$ & $0.959$ & $1.000$ & $1.000$ & $1.000$ & $1.000$ & $1.000$ & $1.000$ & $0.993$  \\
R77    & $0.253$ & $0.954$ & $1.000$ & $1.000$ & $1.000$ & $1.000$ & $1.000$ & $1.000$ & $0.993$     \\
R78    & $0.253$ & $0.953$ & $1.000$ & $1.000$ & $1.000$ & $1.000$ & $1.000$ & $1.000$ & $0.993$     \\
R7L    & $0.254$ & $0.951$ & $1.000$ & $1.000$ & $1.000$ & $1.000$ & $1.000$ & $1.000$ & $0.993$     \\ \hline
P-Test & $0.021$ & $0.014$ & $0.012$ & $0.016$ & $0.006$ & $0.013$ & $0.008$ & $0.007$ & $0.000$  \\
D-Test & $0.048$ & $0.039$ & $0.052$ & $0.046$ & $0.025$ & $0.024$ & $0.011$ & $0.003$ & $0.000$  \\
KS-Test& $0.177$ & $0.887$ & $1.000$ & $1.000$ & $1.000$ & $1.000$ & $1.000$ & $1.000$ & $0.852$ \\
KL-Test& $0.221$ & $0.637$ & $0.869$ & $0.939$ & $0.932$ & $0.775$ & $0.149$ & $0.000$ & $0.000$ \\
\hline\hline
\end{tabular}
\end{center}
\end{footnotesize}
\end{table}
}

\clearpage
\phantomsection
\bibliographystyle{ims}
\bibliography{D:/Dropbox/Common/Latex/mybibliography}

\end{document}